\documentclass[11pt]{article}
\usepackage{amsmath,amsfonts,amssymb,psfrag,epsf,epsfig}
\usepackage{exscale,makeidx,bbm,color}
\usepackage{graphicx}
 \DeclareGraphicsExtensions{.eps, .ps}
\usepackage{stmaryrd}

\newtheorem{prop}{Proposition}
\newtheorem{defn}[prop]{Definition}
\newtheorem{lm}[prop]{Lemma}
\newtheorem{ex}[prop]{Example}
\newtheorem{thm}[prop]{Theorem}
\newtheorem{cor}[prop]{Corollary}

\newtheorem{rem}[prop]{Remark}
\newenvironment{pf}{\begin{trivlist}\item[] \textbf{Proof.}}
                     {\hspace*{\fill} $\square$\end{trivlist}}

\newenvironment{pfofthmasc}{\begin{trivlist}\item[] \textbf{Proof of Theorem \ref{thmasc}.}}
                     {\hspace*{\fill} $\square$\end{trivlist}}
\newenvironment{pfofprop21}{\begin{trivlist}\item[] \textbf{Proof of Proposition \ref{prop21}.}}
                     {\hspace*{\fill} $\square$\end{trivlist}}
\newenvironment{pfofprop2729}{\begin{trivlist}\item[] \textbf{Proof of Propositions \ref{dfprunbranch} and \ref{dfprunbranch2}.}}
                     {\hspace*{\fill} $\square$\end{trivlist}}

\newcommand{\ed}{\mbox{$ \ \stackrel{d}{=}$ }}

\newcommand{\convd}{\overset{d}{\underset{{n\rightarrow \infty}}{\longrightarrow}}}
\newcommand{\eq}{\begin{equation}}
\newcommand{\en}{\end{equation}}

\newcommand{\beq}{\begin{eqnarray*}}
\newcommand{\eeq}{\end{eqnarray*}}

\newcommand{\ts}{\textstyle}
\def\build#1_#2^#3{\mathrel{\mathop{\kern 0pt#1}\limits_{#2}^{#3}}}

\newcommand{\bP}{\mathbb{P}}

\newcommand{\bN}{\mathbb{N}}
\newcommand{\bR}{\mathbb{R}}
\newcommand{\bD}{\mathbb{D}}

\newcommand{\bE}{\mathbb{E}}

\newcommand{\bT}{\mathbb{T}}

\newcommand{\cA}{\mathcal{A}}

\newcommand{\cF}{\mathcal{F}}

\newcommand{\cL}{\mathcal{L}}

\newcommand{\cP}{\mathcal{P}}

\newcommand{\cS}{\mathcal{S}}
\newcommand{\cT}{\mathcal{T}}

\newcommand{\cX}{\mathcal{X}}
\newcommand{\cY}{\mathcal{Y}}

\newcommand{\fk}{\mathbf{k}}
\newcommand{\bZ}{\mathbb{Z}}

\newcommand{\ft}{\mathbf{t}}

\def\beqlb{\begin{eqnarray}}\def\eeqlb{\end{eqnarray}}
\def\beqnn{\begin{eqnarray*}}\def\eeqnn{\end{eqnarray*}}

\begin{document}

\title{\vspace{-0.7cm}
Invariance principles for pruning processes of Galton-Watson trees}
\author{
Hui He%
\thanks{%
Laboratory of Mathematics and Complex Systems, School of Mathematical Sciences, Beijing Normal University, Beijing 100875, P.R.China; email hehui@bnu.edu.cn}
\and
Matthias Winkel\thanks{%
Department of Statistics, University of Oxford, 1 South Parks Road, Oxford OX1 3TG, UK; email winkel@stats.ox.ac.uk}}

\maketitle

\vspace{-0.4cm}

\begin{abstract} Pruning processes $(\cF(\theta),\theta\ge 0)$ have been studied separately for Galton-Watson trees and for
  L\'evy trees/forests. We establish here a limit theory that strongly connects the two studies. This solves an open problem
  by Abraham and Delmas, also formulated as a conjecture by L\"ohr, Voisin and Winter. Specifically, we show that for any  
  sequence of Galton-Watson forests $\cF_n$, $n\ge 1$, in the domain of attraction of a L\'evy forest $\cF$, suitably scaled
  pruning processes $(\cF_n(\theta),\theta\ge 0)$ converge in the Skorohod topology on cadlag functions with values in the space of (isometry
  classes of) locally compact real trees to limiting pruning processes. We separately treat pruning at branch points and pruning at edges. We apply our results to study ascension times and Kesten trees and forests. 

\emph{AMS 2010 subject classifications: Primary 60J80; Secondary 60J25, 60F17.\newline
}
\end{abstract}

\section{Introduction}

Consider a rooted combinatorial tree $(\ft,\rho)$, i.e.\ a connected acyclic graph $\ft$ with vertex set $V(\ft)$, edge set $E(\ft)$ and a special vertex $\rho\in V(\ft)$ called the \em root\em. Given a subset $A\subseteq E(\ft)$, we define the \em pruned subtree \em $(\ft_A,\rho)$ as the connected component $\ft_A$ of $\ft\setminus A$ containing $\rho$. Here $\ft\setminus A$ is the subgraph of $\ft$ with vertex set $V(\ft)$ and edge set $E(\ft)\setminus A$. Given an increasing family $(A(\theta),\theta\ge 0)$ of subsets of $E(\ft)$ with $A(0)=\varnothing$%
, we obtain a \em pruning process \em $\ft(\theta)=\ft_{A(\theta)}$, $\theta\ge 0$.   

In this paper, we establish a limit theory for certain random pruning processes associated with Galton-Watson trees. A Galton-Watson tree with offspring distribution $\xi$ on $\bN=\{0,1,2,\ldots\}$, or a \em ${\it GW}(\xi)$-tree \em for short, is the family tree $(\tau,\rho)$ of a population, in which, beginning with a \em progenitor \em $\rho$, each individual has an independent $\xi$-distributed number of children. We represent individuals by vertices $v\in V(\tau)$ and the parent-child relation by the edge set $E(\tau)$. For each $v\in V(\tau)$, let $E_v(\tau)$ be the set of edges from $v$ to its children (excluding the edge to its parent). If $E_v(\tau)=\varnothing$ then $v$ is a \em leaf \em of $\tau$. If $\#E_v(\tau)\ge 2$ then $v$ is a \em branch point\em. We define the sets ${\rm Lf}(\tau)$ of leaves and ${\rm Br}(\tau)$ of branch points.
Then $E(\tau)=\bigcup_{v\in V(\tau)\setminus{\rm Lf}(\tau)}E_v(\tau)$ is a disjoint union.   

Several pruning processes have appeared in the literature. Aldous and Pitman \cite{AP98} studied \em pruning at edges \em of a Galton-Watson tree $(\tau,\rho)$, where each edge $e\in E(\tau)$ has an independent exponentially distributed (${\rm Exp}(1)$) pruning time $M_e$ so that the set of edges pruned by time $\theta$ is 
$$A^E(\theta)=\{e\in E(\tau)\colon M_e\le\theta\},\qquad \theta\ge 0.$$ 
Abraham, Delmas and He \cite{ADH12} introduced pruning processes that exhibit \em pruning at branch points \em (also called \em pruning at nodes\em), where each branch point $v\in {\rm Br}(\tau)$ has an independent ${\rm Exp}(\#E_v(\tau)-1)$ pruning time $M_v$ that turns the branch point $v$ into a leaf. We obtain this by setting 
  $$A^B(\theta)=\bigcup_{v\in{\rm Br}(\tau)\colon M_v\le\theta}E_v(\tau),\qquad \theta\ge 0.$$
We denote the two pruning processes by $\tau^E(\theta)=\tau_{A^E(\theta)}$ and $\tau^B(\theta)=\tau_{A^B(\theta)}$, $\theta\ge 0$, \pagebreak respectively.  
The literature on invariance principles for discrete Galton-Watson processes goes back a long time, see e.g.\ Grimvall \cite{Gri74}. The starting point for a limit theory for pruning processes in the present paper is a recent extension to include the richer structure of their genealogical forests of Galton-Watson trees \cite{DuLG02,DuWi12,HL13}. In particular, it was shown that the only possible limits are L\'evy forests. L\'evy forests are parametrised by an initial distribution $\varrho$ on $[0,\infty)$ and a \em branching mechanism\em
\eq\label{brmech} \psi(u)=\alpha u+\beta u^2+\int_{(0,\infty)}\left(e^{-ur}-1+ur1_{(0,1)}(x)\right)\pi(dr)
\en
for some $\alpha\in\bR$, $\beta\ge 0$ and $\pi$ with $\int_{(0,\infty)}(1\wedge r^2)\pi(dr)<\infty$, also satisfying two further conditions: 
\eq\label{grey}\mbox{[Grey]}\quad\int^\infty\frac{du}{\psi(u)}<\infty\qquad\mbox{and [conservativity]}\quad
\int_{0+}\frac{du}{|\psi(u)|}=\infty.
\en
In view of the Grey condition, we can define $\eta\colon(0,\infty)\!\rightarrow\!(0,\infty)$ such that $\int_{\eta(h)}^\infty du/\psi(u)\!=\!h$. It is well-known (e.g.\ \cite[(69)]{DuWi12}) that $\int_{[0,\infty)}e^{-x\eta(h)}\varrho(dx)$ is the probability that a $(\psi;\varrho)$-L\'evy forest has height less than $h$.
The setting for the invariance principle is, for each $n\ge 1$,
\begin{itemize}\item an offspring distribution $\xi_n=(\xi_n(k),k\ge 0)$ with $\xi_n(1)<1$,
  \item an associated \em step distribution \em $\nu_n=(\nu_n(k),k\ge -1)$ given by $\nu_n(k)=\xi_n(k+1)$,
  \item an initial distribution $\mu_n=(\mu_n(k),k\ge 0)$ with $\mu_n(0)<1$, and
  \item a Galton-Watson real forest $\cF_n$ of a $\mu_n$-distributed number of independent ${\rm GW}(\xi_n)$-trees.
\end{itemize}
Here, \em real forests \em or \em forests of real trees \em are representations of forests of rooted combinatorial trees in the space $\bT$ of (isometry classes) of rooted locally compact metric space trees equipped with the Gromov-Hausdorff topology, see Section \ref{sectech} for a summary, \cite{Gro,EPW06,DuWi07,DuWi12} for details and \cite{EW06,GPW,Mie09,LVW} for related developments. 

We denote by $e^{-\eta_n(h)}$ the probability that a ${\rm GW}(\xi_n)$-tree has height less than $h\in\bN$, and by $\lfloor r\rfloor$ the integer part of $r\in[0,\infty)$, i.e.\ $\lfloor r\rfloor=k$, where $k\in\bN$ and $k\le r<k+1$.

\begin{thm}[Invariance principle for trees, Theorem 4.15 of \cite{DuWi12}]\label{invprinc}$\;\!\!\!$ Suppose that there is a positive sequence $\gamma_n\rightarrow\infty$ such that, as $n\rightarrow\infty$,
\eq\nu_n(\ts n\,\cdot\,)^{*\lfloor n\gamma_n\rfloor}\rightarrow\nu\quad\mbox{and}
\quad\mu_n(\ts n\,\cdot\,)\rightarrow\varrho\quad\mbox{weakly,}\qquad\mbox{and}\quad n\eta_n(\lfloor\ts\gamma_n\,\cdot\,\rfloor)\rightarrow \eta\quad\mbox{pointwise,}\label{invass}
\en
where $\nu$ is such that $\int_{\bR}e^{-rx}\nu(dx)=e^{\psi(r)}$ for a branching mechanism {\rm(\ref{brmech})} satisfying {\rm(\ref{grey})}. 
Then\vspace{-0.1cm} $$\cF_n/\gamma_n\convd\cF\qquad\mbox{in }\bT,\vspace{-0.2cm}$$ for a $(\psi;\varrho)$-L\'evy forest $\cF$, where $\convd$ denotes convergence in distribution, as $n\rightarrow\infty$.
\end{thm}
As $\bT$-valued random variables, (representatives of) L\'evy forests $\cF$ are equipped with a $\sigma$-finite length measure $\ell$ supported by non-leaf vertices. See Section \ref{sectech}. Aldous and Pitman \cite[Section 2.2]{AP98b} considered a fragmentation process, which in a setting with a root 
gives rise to a pruning process $(\cT^{\rm AP}(\theta),\theta\ge 0)$ for the Brownian Continuum Random Tree (CRT) \linebreak[2] $\cT$ of \cite{Ald91}. Specifically, $\cT$ is equipped with the set $\cA^{\rm AP}(\theta)$ of atoms of a Poisson random measure with intensity measure $\theta\ell$ in such a way that $(\cA^{\rm AP}(\theta),\theta\ge 0)$ is an increasing family, and $\cT^{\rm AP}(\theta)$ is the connected component of $\cT\setminus\cA^{\rm AP}(\theta)$ containing the root. The same construction applies to any L\'evy forest to give an \em Aldous-Pitman pruning process \em $(\cF^{\rm AP}(\theta),\theta\ge 0)$. This is an analogue of pruning at edges, because the (countable) set of branch points of degree $\ge 3$ has zero $\ell$-measure. See Section \ref{secprunedge}. Aldous and Pitman also establish a convergence of a discrete model based on uniform trees to a continuum limit, at the level of component sizes rather than trees, and they study a time-reversal of this process, the standard additive coalescent. 

The Aldous-Pitman pruning process for the Brownian CRT was generalised differently in \cite{AbS,Mie05,AD08,AD12} and placed in the tree-valued framework of pruning processes $(\cT^{\rm AD}(\theta),\theta\ge 0)$ for L\'evy trees \cite{LGLJ,DuLG02,DuLG05}. Abraham and Delmas \cite{AD12} pointed out the analogy between the Galton-Watson and L\'evy tree pruning models, but left open the question of a limit theory. Their generalisation is based on a measure $\omega$ constructed in \cite{Mie05,DuLG05}, which is supported by the branch points of $\cT$ of infinite degree. See Section \ref{sectpruneGW}. Specifically, cut points are placed in the set $\cA^{\rm AD}(\theta)$ of atoms of a Poisson random measure with intensity measure $\theta\omega+2\beta\theta\ell$, where $\beta$ is the quadratic coefficient in (\ref{brmech}). By considering forests of L\'evy trees, we can construct \em Abraham-Delmas pruning processes \em $(\cF^{\rm AD}(\theta),\theta\ge 0)$ for any L\'evy forest. This is an analogue of pruning at branch points. More precisely, $\omega$ is based on (suitably rescaled limiting) sizes of branch points and provides rates proportional to size, just as in the combinatorial pruning at branch points. 

More general pruning operations and pruning processes were introduced in \cite{ADV,He2014}, while \cite{ADH13} studied a two-parameter process that combines pruning of \cite{AD12} with growth of L\'evy trees of \cite{DuWi07}.
L\"ohr, Voisin and Winter \cite{LVW} started a systematic study of pruning processes as instances of a Markov process on a new space of bi-measure $\bR$-trees. In \cite[Section 4]{LVW}, examples of a limit theory for Aldous-Pitman pruning processes of Brownian and stable CRTs are obtained (see also their Remark 4.5 on a possible generalisation to compact L\'evy trees), but the general case of Abraham-Delmas pruning is left as a conjecture (in their Example 4.6). Their notion of convergence is different from ours. As their ``sampling and pruning measures'' only depend on the metric structure of the trees we consider here, we use the usual Gromov-Hausdorff metric. We offer a careful discussion in Section \ref{lvw}, after introducing some technical details. In our framework, we make precise and prove their conjecture (this is our Theorem \ref{thmprocconv}).

We denote by $\bD([0,\infty),\bT)$ the space of $\bT$-valued cadlag functions, equipped with the Skorohod topology. The main result of the present paper is the following.
\begin{thm}[Invariance principle for pruning at branch points]\label{thmprocconv} In the setting of Theorem \ref{invprinc}, the associated pruning processes  $(\cF_n^B(\theta),\theta\ge 0)$ with pruning at branch points converge: $$(\cF_n^B(\theta/n)/\gamma_n,\theta\ge 0)\convd(\cF^{\rm AD}(\theta),\theta\ge 0)\qquad \text{in }\bD([0,\infty),\bT),$$
where the limit is the Abraham-Delmas pruning process associated with a $(\psi;\varrho)$-L\'evy forest $\cF$. 
\end{thm}
We also establish a corresponding general result for Aldous-Pitman pruning at edges, as follows.

\begin{thm}[Invariance principle for pruning at edges]\label{thmprocconvedge} In the setting of Theorem \ref{invprinc}, the associated pruning processes $(\cF_n^E(\theta),\theta\ge 0)$ with pruning at edges converge: 
$$(\cF_n^E(\theta/\gamma_n)/\gamma_n,\theta\ge 0)\convd(\cF^{\rm AP}(\theta),\theta\ge 0)\qquad \text{in }\bD([0,\infty),\bT),$$
where the limit is the Aldous-Pitman pruning process associated with a $(\psi;\varrho)$-L\'evy forest $\cF$.
\end{thm}
The limits $(\cF^{\rm AD}(\theta),\theta\ge 0)$ and $(\cF^{\rm AP}(\theta),\theta\ge 0)$ coincide if and only if $\alpha=\pi=0$. In this ``Brownian case'', we have two convergence results for the same limiting process, with pre-limiting processes that only exhibit either pruning at branch points or pruning at edges. This can be explained by the prevailance of binary branch points in this case. More generally, while in the case $\beta>0$ the $\cF^{\rm AD}$ process includes features of pruning at edges, this feature is not needed for the pre-limiting processes in Theorem \ref{thmprocconv}, contrary to the conjecture of \cite{LVW}. In the case $\beta=0$ on the other hand, we typically have $\gamma_n/n\rightarrow 0$, see e.g.\ Lemma \ref{lm52} where $\gamma_n=\psi^\prime(n)$ with $\gamma_n/n\rightarrow 2\beta$, so the scaling of the pruning parameter is different in the two theorems. \pagebreak[2]

Let us briefly discuss our strategy to prove Theorems \ref{thmprocconv} and \ref{thmprocconvedge}. The first step is to reduce to statements about suitably $h$-erased pruning processes (Corollary \ref{cor25}), generalising the powerful notion of $h$-erasure  \cite{Nev86b,EPW06,DuWi12} from $\bT$ to decreasing $\bT$-valued functions. The second step is to compute the distributions of $h$-erased pruning processes
(Propositions \ref{edgeerase}, \ref{edgeerasediscr}, \ref{prop21} and \ref{propRhGWprun}). This leads outside the framework of pruning processes considered in \cite{LVW}, since pruning times will no longer be exponentially distributed (but mixed exponential). However, pre-limiting and limiting forests are now discrete, and the main step is to establish a new general convergence result (Theorem \ref{thmconvH}) for pruning processes in the framework of an invariance principle from \cite{DuWi12} for Galton-Watson real trees that converge to Galton-Watson real trees with exponentially distributed edge lengths. The final step is to apply 
Theorem \ref{thmconvH} to complete the proof of Theorem \ref{thmprocconv} in Section \ref{secthm2} and to adapt the proof to also prove
Theorem \ref{thmprocconvedge} in Section \ref{secthm3}.

These methods are very general an, in principle, apply to any sequence of discrete tree-valued pruning or tree growth processes, see articles from our bibliography and references therein. 

As an application of the results, we study Kesten(-L\'evy) trees, i.e.\ critical Galton-Watson (and L\'evy trees) suitably conditioned to have infinite height, following \cite{Kes86,Duq-08}. We derive new invariance principles (Theorems \ref{kesconv} and \ref{kesconvedge}) for pruning processes of Kesten-L\'evy trees and Kesten-L\'evy forests $(\cF_*^{\rm AD}(\theta),\theta\!\ge\! 0)$ in Section \ref{ascsec1}, while Section \ref{ascsec2} studies extended pruning processes $(\cF^{\rm AD}(\theta),\theta\!>\!-\theta_0)$ from their ascension time $A\!=\!\inf\{a\!\ge\! 0\colon\cF^{\rm AD}(-a)\mbox{ infinite}\}$. Specifically, we deduce from our invariance principles and discrete results of \cite{ADH12} that (Theorem \ref{thmasc})\vspace{-0.3cm}
  $$(\cF^{\rm AD}(\theta),\theta\ge -A)\ed(\cF_*^{\rm AD}(W/x+\theta),\theta\ge -\Theta),$$
for suitable $(W,\Theta)$, a new result for forests related to \cite{AD12}, who used different methods to establish similar results for single L\'evy trees, which we could now also deduce by limiting considerations. 

The structure of the paper is, as follows. 
Section \ref{sectech} gives an introduction to the Gromov-Hausdorff topology and Section \ref{secskor} to Skorohod's topology, and we also derive a general convergence criterion based on suitably $h$-erased pruning processes. Section \ref{lvw} discusses the topology of \cite{LVW}. In Section \ref{gwsec}, we introduce Galton-Watson real trees and L\'evy forests. Section \ref{secprunedge} discusses pruning at edges for Galton-Watson real trees, leading up to the Aldous-Pitman pruning processes for L\'evy forests. Section \ref{sectpruneGW} discusses pruning at branch points for Galton-Watson real trees, leading up to the Abraham-Delmas pruning processes for L\'evy forests.


In Section \ref{secH} we state and prove Theorem \ref{thmconvH}. In Section \ref{secaux} we obtain some auxiliary results that are used in Sections \ref{secthm2} and \ref{secthm3} to complete the proofs of Theorems \ref{thmprocconv} and \ref{thmprocconvedge}, respectively. In Section \ref{secthmlast}, we establish an invariance principle closely related to Theorem \ref{thmprocconvedge} but based on pruning at branch points with degree-independent rates. Applications to ascension times and Kesten(-L\'evy) trees and forests are discussed in Section \ref{ascsec}.

\section{Preliminaries on topologies for tree-valued processes}\label{prel}

\subsection{Real trees and the Gromov-Hausdorff topology on $\bT$}\label{sectech}

Following \cite{EPW06}, a rooted real tree $(T,d,\rho)$ is a metric space $(T,d)$ with a \em root \em $\rho\in T$, such that any two points $v,w\in T$ are connected by a unique injective path $[[v,w]]$, which furthermore has length $d(v,w)$. We denote by $\bT$ the set of all 
root-preserving isometry classes of complete separable locally compact rooted real trees. For two rooted real trees $(T,d,\rho)$ and $(T^\prime,d^\prime,\rho^\prime)$, we consider\vspace{-0.1cm}
$$\delta((T,d,\rho),(T^\prime,d^\prime,\rho^\prime))=\inf_{\phi,\phi^\prime}\overline{\Delta}_X^{\rm Haus}(\phi(T),\phi^\prime(T^\prime)),$$
where the infimum is taken over all pointed metric spaces $(X,\Delta_X,\rho_X)$ and all isometric embeddings $\phi\colon T\rightarrow X$ and $\phi^\prime\colon T^\prime\rightarrow X$ with $\phi(\rho)=\phi^\prime(\rho^\prime)=\rho_X$. Here, 
$$\overline{\Delta}_X^{\rm Haus}=\int_0^\infty\Delta_X^{\rm Haus}(\,\cdot\,\cap B(\rho_X,r),\,\cdot\,\cap B(\rho_X,r))e^{-r}dr$$ 
is a localised version of the Hausdorff distance $\Delta^{\rm Haus}_X$, based on restrictions to balls $B(\rho_X,r)=\{x\in X\colon\Delta_X(x,\rho_X)\le r\}$, $r\ge 0$. The distance function $\delta$ induces a metric on $\bT$, the Gromov-Hausdorff metric, which is also denoted by $\delta$. We abuse notation and write $T\in\bT$ to denote an isometry class. Occasionally, it is convenient to work with representatives. Every $T\in\bT$ can be represented as a metric subspace of $X=\ell_1(\bN)=\{x\!\in\![0,\infty)^\bN\colon||x||_1<\infty\}$ with metric induced by the $\ell_1$-norm $||x||_1=||(x_n)_{n\ge 0}||_1=\sum_{n\ge 0}|x_n|$ and root $0\in\ell_1(\bN)$. For $(X,\Delta_X,\rho_X)$, we denote the space of complete locally compact real trees in $X$ by $\bT_X$. See e.g.\ \cite{DuWi12} for details. 

For a rooted real tree $(T,d,\rho)$, we consider the height $\Gamma(T)=\sup\{d(\rho,v)\colon v\in T\}$, for any vertex $v\in T$ the subtree $T_v=\{w\in T\colon v\in[[\rho,w]]\}$ above $v$, and for any $h>0$ the $h$-erasure operation which sets $R^h(T)=\{\rho\}$ if $\Gamma(T)\le h$ and $R^h(T)=\{v\in T\colon\Gamma(T_v)\ge h\}$. Then $\Gamma$ and $R^h$ induce corresponding
functions $\Gamma\colon\bT\rightarrow[0,\infty]$ and $R^h\colon\bT\rightarrow\bT$, see \cite{EPW06,DuWi12}.

For any $v\in T$, let ${\rm n}(v,T)\in\bN\cup\{\infty\}$ be the \em degree of $v$ in $T$\em, i.e. the number of connected components of $T\setminus\{v\}$. We say 
$v\neq\rho$ is a \em branch point \em if ${\rm n}(v,T)\ge 3$ and a \em leaf \em if ${\rm n}(v,T)=1$. We denote the set of branch points by ${\rm Br}(T)$, the set of leaves by ${\rm Lf}(T)$. For $a>0$, consider ${\rm Blw}(T,a)=\{v\in T\!: d(\rho,v)\le a\}$ and the quotient space $({\rm Abv}(T,a),d_a,[\rho]_a)$ of $(T,d,\rho)$ by the equivalence relation $v\sim_a w$ iff $v=w$ or $v,w\in{\rm Blw}(T,a)$. We can represent ${\rm Abv}(T,a)=\{[\rho]_a\}\cup\bigcup_{i\in I(a)}T_i^\circ(a)$ as union of the  connected components of $\{v\in T\colon d(\rho,v)>a\}$, indeed as concatenation of trees $(T_i(a),d_i,\rho_i)$ at $\rho_i$, $i\in I(a)$, which we write as  $({\rm Abv}(T,a),d_a,[\rho]_a)=\circledast_{i\in I(a)}(T_i(a),d_i,\rho_i)$. Then
$({\rm Abv}(T,a),d_a,[\rho]_a)$ and $({\rm Blw}(T,a),d,\rho)$ are rooted real trees and induce ${\rm Abv},{\rm Blw}\colon\!\bT\!\times\![0,\infty)\!\rightarrow\!\bT$. 

It is a direct consequence of local compactness of $T$ that ${\rm Blw}(R^h(T),a))$ has at most finitely many leaves and branch points all with ${\rm n}(v,T)<\infty$, for all $h>0$ and $a>0$. In particular, there is a finite length measure $\ell$ that assigns length $d(v,w)$ to $[[v,w]]$ for all $v,w\in{\rm Blw}(R^h(T),a)$. While $T$ may have uncountable ${\rm Lf}(T)$ and countable dense ${\rm Br}(T)$ with ${\rm n}(v,T)=\infty$, $v\in{\rm Br}(T)$, this length measure consistently extends to a $\sigma$-finite measure $\ell$ on $T$, which is supported by $T\setminus{\rm Lf}(T)$. With trees with finite ${\rm Br}(T)\cup{\rm Lf}(T)$ in mind, we further define, for general $T$,
\begin{itemize}\item ${\rm n}(T)={\rm n}(\rho,T)\in\bN\cup\{\infty\}$, which we refer to as \em the number of trees in the forest $T$\em,
  \item $D(T)=\inf\{d(\rho,v)\colon v\in{\rm Lf}(T)\cup {\rm Br}(T)\setminus\{\rho\}\}\in[0,\infty]$, the height of the first branch point.
\end{itemize}
If ${\rm n}(T)=1$ and $D(T)\in(0,\infty)$, we furthermore define
\begin{itemize}\item $\vartheta(T)\!=\!{\rm Abv}(T,D(T))$, the (concatenation of) subtrees (if any) above the first branch point,
  \item ${\bf k}(T)\!=\!{\rm n}(\rho,\vartheta(T))\!\in\!\{0,2,3,\ldots\}\cup\{\infty\}$, the number of subtrees above the first branch point.
\end{itemize} 
If ${\rm n}(T)\neq 1$ or $D(T)\in\{0,\infty\}$, we define $\vartheta(T)=\{\rho\}$ and ${\bf k}(T)=0$. Then ${\rm n}$, $D$, $\vartheta$ and ${\bf k}$ induce functions on $\bT$. We collect some results from \cite{EPW06,DuWi07,DuWi12}.
\begin{prop}[\cite{EPW06,DuWi07,DuWi12}]\label{prop4}\begin{enumerate}
  \item[{\rm (i)}] $(\bT,\delta)$ is a Polish metric space.
  \item[{\rm (ii)}] $R^h\colon\bT\rightarrow\bT$ is continuous, $\delta(R^h(T),R^{h^\prime}(T))\le|h-h^\prime|$, $R^h\circ R^{h^\prime}=R^{h+h^\prime}$ for all $h,h^\prime\ge 0$.
  \item[{\rm (iii)}] $\Gamma\colon\bT\rightarrow[0,\infty]$ and ${\rm Abv},{\rm Blw}\colon\bT\times[0,\infty)\rightarrow\bT$ are continuous.
  \item[{\rm (iv)}] ${\rm n}\colon\bT\rightarrow\bN\cup\{\infty\}$, $D\colon\bT\rightarrow[0,\infty]$, $\vartheta\colon\bT\rightarrow\bT$, and ${\bf k}\colon\bT\rightarrow\bN\cup\{\infty\}$ are measurable.
  \item[{\rm (v)}] $D(R^h(T))>0$ and ${\bf k}(R^h(T))<\infty$ for all $T\in\bT$, $h>0$.
  \end{enumerate}\pagebreak[2]
\end{prop}

\subsection{Convergence criteria for Skorohod's topology}\label{secskor}

The convergence in Theorems \ref{thmprocconv} and \ref{thmprocconvedge} takes place in the space $\bD([0,\infty),\bT)$ of cadlag functions taking values in the space $\bT$ of isometry classes of complete separable locally compact rooted real trees. Since $\bT$ equipped with the (localised) Gromov-Hausdorff metric $\delta$ is a Polish metric space, the space $\bD([0,\infty),\bT)$ can be equipped with Skorohod's ($J_1$-)topology. We specialise from the higher generality of \cite[Proposition 3.5.3]{EK86} that for functions $x_n,x\in\bD([0,\infty),\bT)$, $n\ge 0$, we have $x_n\rightarrow x$ in the Skorohod sense, as $n\rightarrow\infty$,
if and only if there exists a sequence of continuous increasing bijections $\lambda_n\colon[0,\infty)\rightarrow[0,\infty)$ such that, as $n\rightarrow\infty$, 
$$\gamma(\lambda_n)=\sup_{\theta^\prime>\theta\ge 0}\left|\log\frac{\lambda_n(\theta^\prime)-\lambda_n(\theta)}{\theta^\prime-\theta}\right|\rightarrow 0\quad\mbox{and}\quad
  \sup_{0\le\theta\le\theta_0}\delta(x_n(\theta),x(\lambda_n(\theta)))\rightarrow 0\quad\mbox{for all $\theta_0>0$.}$$
Skorohod's topology is generated by the metric\vspace{-0.2cm}
$$d_{\rm Sk}(x,x^\prime)=\inf_{\lambda}\left(\gamma(\lambda)\vee\int_0^\infty e^{-u}\sup_{\theta\ge 0}\left(\delta(x(\theta\wedge u),x^\prime(\lambda(\theta)\wedge u))\wedge 1\right)du\right),\vspace{-0.1cm}$$
where the infimum is taken over all continuous increasing bijections $\lambda\colon[0,\infty)\rightarrow[0,\infty)$. With this definition,
$(\bD([0,\infty),\bT),d_{\rm Sk})$ is complete and separable. See \cite[Theorem 3.5.6]{EK86}.
\begin{lm} For all $x\in\bD([0,\infty),\bT)$ and $h>0$, we have $d_{\rm Sk}(R^h\circ x,x)\le h$.
  Furthermore, $d_{\rm Sk}(R^h\circ x_n,R^h\circ x)\rightarrow 0$ for all $h>0$, as $n\rightarrow\infty$, implies $d_{\rm Sk}(x_n,x)\rightarrow 0$ as $n\rightarrow\infty$.
\end{lm}
\begin{pf}  As $\delta(R^h(T),T)\le h$ for all $T\in\bT$, we have $d_{\rm Sk}(R^h\circ x,x)\le\sup_{\theta\ge 0}\delta(R^h(x(\theta)),x(\theta))\le h$. 
  Let $\varepsilon\!>\!0$ and set $h\!=\!\varepsilon/3$. Then there is $n_0\!\ge\! 0$ such that
   $d_{\rm Sk}(R^h\circ x_n,R^h\circ x)<\varepsilon/3$ for all $n\!\ge\! n_0$. Hence we find for all $n\ge n_0$\vspace{-0.1cm}
   $$d_{\rm Sk}(x_n,x)\le d_{\rm Sk}(x_n,R^h\circ x_n)+d_{\rm Sk}(R^h\circ x_n,R^h\circ x)+d_{\rm Sk}(R^h\circ x,x)<2h+\varepsilon/3=\varepsilon.\vspace{-0.65cm}$$
\end{pf}
We will be interested in pruning processes in the sense of the following general definition.

\begin{defn}[Pruning process]\label{prunproc}\rm Let $(X,\Delta_X,\rho_X)$ be a pointed metric space and $T\in\bT_X$. A right-continuous $\bT_X$-valued process $(\cT(\theta),\theta\ge 0)$ is called a \em pruning process of $T$ \em if it is decreasing for the inclusion partial order on the subsets of $X$ and if $\cT(0)=T$. We say that \em $(\cT(\theta),\theta\ge 0)$ is associated with point measure $\cP$ \em on $(0,\infty)\times T$, if $\cT(\theta)$ is the completion of the connected component of $T\setminus\cA(\theta)$ containing $\rho_X$, where $\cA(\theta)$ is the support of $\cP((0,\theta]\times\cdot)$, 
$\theta \ge 0$. We also call the $\bT$-valued process of isometry classes a pruning process. 
\end{defn}
In the next section, we will introduce families of pruning processes for which $T$ is a random tree and $\cP$ is a random point measure, often a Poisson random measure with some intensity measure of the form $d\theta\nu(dv)$. In this and similar settings, the following convergence criterion is useful.

\begin{cor}\label{cor25} If $\cX_n\!:=\!(\cT_n(\theta),\theta\!\ge\! 0)$, $n\!\ge\! 1$, and $\cX\!:=\!(\cT(\theta),\theta\!\ge\! 0)$ are pruning processes and if $\cX_n^h\!:=\!(\cT_n(\theta)\cap R^h(\cT_n(0)),\theta\!\ge\! 0)\convd(\cT(\theta)\cap R^h(\cT(0)),\theta\!\ge\! 0)\!=:\!\cX^h$ in the Skorohod sense \vspace{-0.1cm} for all $h>0$, then $(\cT_n(\theta),\theta\ge 0)=\cX_n\convd\cX=(\cT(\theta),\theta\ge 0)$ in the Skorohod sense. 

  The condition $\cX_n^h\convd\cX^h$ for all $h>0$ can be further relaxed: if for all $h>0$ there \vspace{-0.1cm} is a sequence $h_n\rightarrow h$
  for which $\cX^{h_n}_n\convd\cX^h$, then $\cX_n\convd\cX$, all in the Skorohod sense. 
\end{cor}
\begin{pf} Let $f\colon\bD([0,\infty),\bT_X)\rightarrow\bR$ be bounded and uniformly continuous, and let $\varepsilon>0$. Then there is
  $h\!>\!0$ such that for all $x,x^\prime\!\in\!\bD([0,\infty),\bT_X)$ with $d_{\rm Sk}(x,x^\prime)\!\le\! h$ we have 
  $|f(x)-f(x^\prime)|<\varepsilon/3$. Also, by hypothesis, there is $n_0\ge 1$ such that 
  $|\bE(f(\cX_n^h))-\bE(f(\cX^h))|<\varepsilon/3$ for all $n\ge n_0$. Now recall that ($\bT_X$-valued) pruning processes are decreasing (for the 
  inclusion partial order), so $R^h(\cT_n(\theta))\subseteq\cT_n(\theta)\cap R^h(\cT_n(0))\subseteq\cT_n(\theta)$, $\theta\ge 0$, 
  and by the previous lemma $d_{\rm Sk}(\cX_n^h,\cX_n)\le d_{\rm Sk}(R^h\circ\cX_n,\cX_n)\le h$ almost surely, so that for all $n\ge n_0$
  \begin{eqnarray*}\Big|\bE(f(\cX_n))\!-\!\bE(f(\cX))\Big|&\!\!\!\!\le\!\!\!\!&\left|\bE(f(\cX_n))\!-\!\bE(f(\cX_n^h))\right|\!+\!\left|\bE(f(\cX_n^h))\!-\!\bE(f(\cX^h))\right|\!+\!\left|\bE(f(\cX^h))\!-\!\bE(f(\cX))\right|\\
           &\!\!\!\!<\!\!\!\!&\varepsilon/3+\varepsilon/3+\varepsilon/3=\varepsilon.
  \end{eqnarray*}
  The required weak convergence follows by a suitable version of the Portmanteau theorem, see e.g.\ \cite[Theorem 2.1]{Bil68}. The relaxation of the condition to $\cX_n^{h_n}\convd\cX^h$ is straightforward.
\end{pf}

The reader may wonder why we consider $\cX^h\!=\!(\cT(\theta)\cap R^h(\cT(0)),\theta\!\ge\! 0)$, where we first $h$-erase then prune, instead of $R^h\circ\cX=(R^h(\cT(\theta)),\theta\ge 0)$, where we first prune then $h$-erase. The key advantage of $\cX^h$ is that it is a pruning process associated with a point measure $\cP^h$ that is just the restriction of the point measure $\cP$ of $\cX$ to $R^h(\cT(0))$. On the other hand, while $R^h\circ\cX$ is a pruning process, it is not associated with a natural point measure $\cP^h$, in general:\pagebreak[2]

\begin{ex}\label{ex8}\rm Consider a $Y$-shaped real tree $(T,d,\rho)$ with branch point $b\in T$ connecting three edges of unit length, namely a trunk $[[\rho,b]]$ and two branches $]]b,L_1]]$ and $]]b,L_2]]$ leading to two leaves $L_1$ and $L_2$. Let $h<2/3$. Then $R^h$ just shortens the two branches by $h$. If $\cP=\delta_{(\theta,v)}$ for some $\theta>0$ and $v\in]]b,L_2]]$ with $d(b,v)=h/2$, then the pruning process $\cX$ associated with $T$ and $\cP$ is such that the pruning event in $R^h\circ\cX$ prunes at $v$ and also $h$-erases $]]b,v]]$ entirely, but not the other branch at $b$, so any point measure $\cP^h$ on $(0,\infty)\times R^h(T)$ associated with $R^h\circ\cX$ will require infinitely many points on $]]b,v]]$ accumulating at $b$, which is not so useful. \vspace{-0.1cm}
\end{ex}

\subsection{Discussion of the topology and results by L\"ohr, Voisin and Winter \cite{LVW}}\label{lvw}

In \cite[Section 2]{LVW}, a topology on bi-measure $\bR$-trees is introduced. While we do not use their topology in the present paper, their results are closely related to ours, and we would like to discuss this in some detail both to clarify the connections and to justify our choice of topology. 

A $k$-pointed measure $\bR$-tree is a triplet $(T,(u_1,\ldots,u_k),\mu^s)$, where $(T,d,\rho)$ is a complete and separable rooted $\bR$-tree, $u_1,\ldots,u_k\in T$ and $\mu^s$ is a finite Borel measure on $T$. The measure $\mu^s$ is the \em sampling measure\em. Two $k$-pointed measure $\bR$-trees are called equivalent if the supports of the sampling measures (with the root added) are isometric by an isometry that preserves the roots, the $k$ points and the sampling measures. The space $\bT^w_k$ of equivalence classes of $k$-pointed measure $\bR$-trees is Polish when equipped with the $k$-pointed Gromov-Prohorov distance\vspace{-0.2cm}
$$\delta_k^w((T,(u_1,\ldots,u_k),\mu^s),(T^\prime\!,(u_1^\prime,\ldots,u_k^\prime),{\mu^s}^\prime))=\inf_{\phi,\phi^\prime}\left(\!\Delta^{\rm Pr}_X(\phi_*\mu^s,\phi^\prime_*{\mu^s}^\prime)+\sum_{i=1}^k\Delta_X(\phi(u_i),\phi^\prime(u_i^\prime))\!\right)\!,\vspace{-0.2cm}$$
where the infimum is over all metric measure spaces $(X,\Delta_X,\rho_X)$ and all isometric embeddings $\phi\colon T\!\rightarrow\! X$ and $\phi^\prime\colon T^\prime\!\!\rightarrow\! X$ with $\phi(\rho)\!=\!\phi^\prime(\rho^\prime)\!=\!\rho_X$, where $\Delta_X^{\rm Pr}$ is the Prohorov distance on $X$ and $\phi_*\mu^s=\mu^s\circ\phi^{-1}$ denotes the push-forward from $T$ to $X$ of the measure $\mu^s$ by the function $\phi$.

A \em bi-measure $\bR$-tree \em is a triplet $(T,\mu^s\!,\nu)$, where $(T,\mu^s)$ is a ($0$-pointed) measure $\bR$-tree and $\nu$ a Borel measure on $T$, which is $\sigma$-finite on and supported by $(T\setminus{\rm Lf}(T))\cup\{v\!\in\! T\colon\!\mu^s(\{v\})\!>\!0\}$, and which is finite on $[[\rho,v]]$ for all $v\!\in\! T$. The measure $\nu$ is called the \em pruning measure\em. Two bi-measure $\bR$-trees are equivalent if the measure $\bR$-trees are equivalent by an isometry that also preserves the pruning measures. We write the set of equivalence classes of bi-measure $\bR$-trees as $\bT^{\rm bi}$. A sequence in $\bT^{\rm bi}$ is said to \em converge in $\bT^{\rm bi}$ \em if random subtrees spanned by the root and $k$ points sampled from (normalised) sampling measures $\delta^w_k$-converge in distribution when equipped with the (finite) restrictions of the pruning measures, for all $k\!\ge\! 1$. This notion of convergence defines a separable metrisable topology on $\bT^{\rm bi}$, but completeness is not claimed. Lack of completeness would not be a problem for us, as limiting trees have already been constructed.

In \cite[Section 3]{LVW}, a $\bT^{\rm bi}$-valued pruning process is associated with each element of $\bT^{\rm bi}$: for a bi-measure $\bR$-tree $(T,\mu^s\!,\nu)$, this is a pruning process of $T$ associated with a Poisson point measure $\cP$ on $(0,\infty)\times T$ with intensity measure $d\theta\nu(dv)$, as in Definition \ref{prunproc}, but also equipped with the restrictions of $\mu^s$ and $\nu$. This pruning process is a stochastically continuous strong Markov process whose distribution on $\bD([0,\infty),\bT^{\rm bi})$ depends continuously on the initial condition. In \cite[Section 4]{LVW}, several examples are given that are relevant to us. Invariance principles in $\bD([0,\infty),\bT^{\rm bi})$ are obtained where trees are encoded in normalised excursions, i.e.\ for Brownian and stable CRTs, in the case of Aldous-Pitman pruning (cf. our Theorem \ref{thmprocconvedge}). 

Let us explore the framework of \cite{LVW} in the general setting of Theorems \ref{thmprocconv} and \ref{thmprocconvedge}. Sampling measures $\mu_n^s$ on $\cF_n$ and $\mu^s$ on $\cF$ do not feature at all. The topology on $\bT^{\rm bi}$ can find application if we can sample from a normalised counting or length measure $\mu_n^s$ on leaves (or vertices) or edges and, on the CRT side, from a normalised mass measure $\mu^s$ supported by the leaves of $\cF$. These measures exist (see \cite{DuLG02,DuWi07}) as finite measures if our locally compact trees are compact, i.e.\ precisely in the special case where the Galton-Watson and L\'evy forests are (sub)critical, and with further localisation we could prove the analogue of Theorem \ref{thmprocconvedge}, but not Theorem \ref{thmprocconv}. 

Pruning measures $\nu_n$ and $\nu$ on $\cF_n$ and $\cF$ are implicit in Theorems \ref{thmprocconv} and \ref{thmprocconvedge}. They capture the pruning mechanism as intensity measures $d\theta\nu_n(dv)$ or $d\theta\nu(dv)$ of Poisson random measures $\cP_n$ or $\cP$ of cut points in the sense of Definition \ref{prunproc}. In the case of pruning at branch points, the pruning measures are suitably rescaled \pagebreak[2] size measures $\nu_n(\{v\})={\rm n}(v,\cF_n)-2$ on branch points $v\in{\rm Br}(\cF_n)$, with $\nu=\omega+2\beta\ell$ on the CRT side, see Section \ref{sectpruneGW}. In the case of pruning at edges, the pruning measure on the CRT side is length measure $\nu=\ell$, see Section \ref{secprunedge}, while \cite{LVW} effectively made an (asymptotically negligible) modification to include pruning at edges. The authors take suitably rescaled counting or length measure $\nu_n$ on vertices or edges of $\cF_n$. 

Using counting measure on vertices corresponds (via the one-to-one correspondence between edges and non-root vertices) to counting measure on edges, when placing cut points at the top ends of edges, whereas removal of an edge more naturally means cutting at the bottom ends of edges (here ``top end'' and ``bottom end'' mean ``vertex further from the root'' and ``vertex closer to the root'', respectively). However, several edges share the same bottom end, so to prune only one edge at the bottom end vertex, we stop looking for a point process (or we would need infinitely many points in $\cP$ for each cut as in Example \ref{ex8}). Intuitively, the difference between top and bottom ends is negligible, as edge lengths tend to zero in the setting of Theorem \ref{thmprocconvedge}. Our methods can handle such more general pruning and prove negligibility, see Sections \ref{secthm3} and \ref{secthmlast}. 

On the other hand, the framework of \cite{LVW} is rather implicit about distances. Indeed, while the Prohorov metric is based on distances, Gromov-Prohorov convergence does not imply Gromov-Hausdorff convergence in general (see e.g.\ \cite{ALW}). In order to
include supercritical Galton-Watson and L\'evy forests, we use the (localised) Gromov-Hausdorff topology in Theorems \ref{thmprocconv} and \ref{thmprocconvedge}. While \cite{LVW} exploited that pruning measures are finite on sampled subtrees (and in $\bT^{\rm bi}$ they need to have such restrictions converge), we exploit that pruning measures are locally finite on $h$-erased subtrees (for length measures, and by Corollary \ref{coromega} also for the Abraham-Delmas pruning measures). As $h$-erased subtrees are discrete (by Proposition \ref{prop4}(v)), we effectively show the convergence of point processes branch by branch and in branch points (in Sections \ref{secthm2} and \ref{secthmlast}) and so establish the analogous convergence of pruning measures restricted to $h$-erased subtrees.

\section{Introduction to $\bT$-valued pruning processes}\label{prel2}

\subsection{Galton-Watson real trees and L\'evy forests}\label{gwsec}

To define Galton-Watson trees we specify the distribution of $(D,{\bf k},\vartheta)$ of Section \ref{sectech} recursively. Specifically, $\vartheta$ will be distributed as $\circledast_{i\in\{1,\ldots,\ell\}}\cT_i$ for independent and identically distributed $\cT_i$, $i\in\{1,\ldots,\ell\}$. If all $\cT_i$ have distribution $Q$, we denote this distribution of $\circledast_{i\in\{1,\ldots,\ell\}}\cT_i$ by $Q^{\circledast\ell}$. 

\begin{defn}[GW-real trees \cite{DuWi12}]\rm\begin{enumerate}\item[(i)] A Galton-Watson real tree with unit edge lengths and  offspring distribution $\xi$ on $\bN$ satisfying $\xi(1)<1$, a ${\rm GW}(\xi)$-real tree for short, is a $\bT$-valued random variable, whose distribution $Q_\xi$ is the unique distribution $Q$ on $\bT$ that satisfies\vspace{-0.3cm}
\begin{equation}\label{discrQ} Q(g(D)1_{\{\fk=i\}}G(\vartheta))=\sum_{m=1}^\infty g(m)\xi(1)^{m-1}\,\xi(i)\,Q^{\circledast i}(G)\vspace{-0.1cm}
\end{equation}
for all $i\in\bN$ and all nonnegative measurable functions $g$ on $[0,\infty)$ and $G$ on $\bT$.
\item[(ii)] Suppose $\xi$ satisfies $\xi(1)=0$  and [conservativity] $\int^{1-}|g_\xi(s)-s|^{-1}dr=\infty$, where $g_\xi(s)=\sum_{n\ge 0}\xi(n)s^n$. A ${\rm GW}(\xi,c)$-real tree, or a Galton-Watson real tree with exponentially distributed edge lengths with parameter $c\in(0,\infty)$ and offspring distribution $\xi$, is a $\bT$-valued random variable, whose distribution $Q_{\xi,c}$ is the unique distribution $Q$ on $\bT$ that satisfies\vspace{-0.2cm}
\begin{equation}\label{contQ} Q(g(D)1_{\{\fk=i\}}G(\vartheta))=\int_0^\infty g(x)ce^{-cx}dx\,\xi(i)\,Q^{\circledast i}(G)\vspace{-0.1cm}
\end{equation}
for all $i\in\bN$ and all nonnegative measurable functions $g$ on $[0,\infty)$ and $G$ on $\bT$.
\item[(iii)] A ${\rm GW}(\xi;\mu)$-real resp.\ ${\rm GW}(\xi,c;\mu)$-real forest is a $\bT$-valued random variable with distribution $P_{\xi}^\mu\!=\!\sum_{\ell\ge 0}\mu(\ell)Q_\xi^{\circledast\ell}$ resp.\ $P_{\xi,c}^\mu\!=\!\sum_{\ell\ge 0}\mu(\ell) Q_{\xi,c}^{\circledast\ell}$ so that $\mu=P_{\xi}^\mu({\rm n}\!\in\!\cdot)=P_{\xi,c}^\mu({\rm n}\!\in\!\cdot)$.
\end{enumerate}
\medskip

\noindent For any metric space $(X,\Delta_X,\rho_X)$, we also refer to a $\bT_X$-valued random variable as a Galton-Watson tree/forest, if its isometry class in $\bT$ has distribution $Q_\xi$, $Q_{\xi,c}$, $P_\xi^\mu$ or $P_{\xi,c}^\mu$, respectively.\pagebreak[2]
\end{defn}
Existence and uniqueness of $Q_\xi$ and $Q_{\xi,c}$ were shown in \cite[Lemma 2.15]{DuWi12}. We can rephrase the definitions of $Q_\xi$ and $Q_{\xi,c}$, as follows: under $Q_\xi$, respectively under $Q_{\xi,c}$,
\begin{enumerate}\item $D$ and $(\fk,\vartheta)$ are independent,
  \item $D\sim{\rm geom}(1-\xi(1))$ respectively $D\sim{\rm Exp}(c)$, where $\sim$ means ``has distribution'',
  \item $\fk\sim\widetilde{\xi}$ where $\widetilde{\xi}(i)=\xi(i)/(1-\xi(1))$, $i=0,2,3,\ldots$, and $\widetilde{\xi}(1)=0$, 
  \item and conditionally given $\{\fk=i\}$, we have $\vartheta\sim Q_\xi^{\circledast i}$ respectively $\vartheta\sim Q_{\xi,c}^{\circledast i}$.\pagebreak[2]
\end{enumerate}
The following theorem demonstrates how ${\rm GW}(\xi,c)$-real trees/forests appear as limits of ${\rm GW}(\xi_n)$-real trees/forests, summarising \cite[Lemma 3.22, Remark 3.23, Theorem 3.24, (89)]{DuWi12}. This result contains an invariance principle analogous to Theorem \ref{invprinc}, but in a discrete limit regime. 
\begin{thm}[\cite{DuWi12}]\label{invprincdiscr} Let $\xi$ be a conservative offspring distribution, $\xi(0)\!<\!1$, $\xi(1)\!=\!0$, $c\!>\!0$. Let $\mu$ be a distribution on $\bN$, $\mu(0)\!<\!1$. Let $\gamma_n\!>\!0$, $n\ge 1$, with $\gamma_n\!\rightarrow\!\infty$ as $n\rightarrow\infty$. Then the following convergences as $n\rightarrow\infty$ are equivalent:
  \begin{enumerate}
    \item[\rm(a)] $\widetilde{\xi}_n\rightarrow\xi$ and
      $\mu_n\rightarrow\mu$ weakly on $\bN$, and $\gamma_n(1-\xi_n(1))\rightarrow c$, where we write $\widetilde{\xi}_n$ for the
      distribution $\xi_n$ conditioned on $\bN\setminus\{1\}$, i.e.
      $\widetilde{\xi}_n(k)=\xi_n(k)/(1-\xi_n(1))$, $k\neq 1$, $\widetilde{\xi}_n(1)=0$.
    \item[\rm(b)] $(\nu_n)^{*\lfloor\gamma_n\rfloor}\rightarrow\overline{\nu}$ and $\mu_n\rightarrow\mu$ weakly on
      $\bZ$, where we write $\overline{\nu}$ for the law of $X_1$ for a compound Poisson process $(X_t,t\ge 0)$ with holding parameter
      $c$ and jump law $\nu(k)=\xi(k+1)$, $k\ge -1$.
    \item[\rm(c)] $P_{\xi_n}^{\mu_n}(\cdot/\gamma_n)\rightarrow P_{\xi,c}^{\mu}$ weakly on $\bT$.
    \item[\rm(d)] $Q_{\xi_n}(\cdot/\gamma_n)\rightarrow Q_{\xi,c}$ weakly on $\bT$, and $\mu_n\rightarrow\mu$ weakly on $\bN$. Furthermore, the joint distributions of $(D,\fk,\vartheta)$ under $Q_{\xi_n}$ converge weakly to those under $Q_{\xi,c}$.\pagebreak[2]
  \end{enumerate}
\end{thm}
Let us note that the class of Galton-Watson forests is closed under $h$-erasure:

\begin{lm}[cf. Kesten \cite{Kes86}, Neveu \cite{Nev86b}]\label{lmdiscrerasure}\begin{enumerate}\item[{\rm (i)}] For a ${\rm GW}(\xi;\mu)$-real forest $\cF$ and $h\!\in\!\bN$, the forest $R^h(\cF)$ is a ${\rm GW}(\xi^h;\mu^h)$-real
  forest, where $\xi^h$ and $\mu^h$ have generating functions
  \begin{equation}\label{xihmuh}g_{\xi^h}(s)=\frac{g_\xi(p+s(1-p))-p}{1-p}\quad\mbox{ and }\quad g_{\mu^h}(s)=g_\mu(p+s(1-p)),
  \end{equation}
  where $p=Q_\xi(\Gamma\le h)$.
  \item[{\rm (ii)}] For a ${\rm GW}(\xi,c;\mu)$-real forest $\cF$ and $h\in(0,\infty)$, the forest $R^h(\cF)$ is a ${\rm GW}(\xi^{h,c},c^h;\mu^h)$-real forest, where $c^h=c(1-g_\xi^\prime(p))$ and $\xi^{h,c}$ and $\mu^h$ have generating functions
  $$g_{\xi^{h,c}}(s)=s+\frac{g_\xi(s+(1-s)p)-(s+(1-s)p)}{(1-p)(1-g_\xi^\prime(p))}\quad\mbox{ and }\quad g_{\mu^h}(s)=g_\mu(p+(1-s)p),
  $$
  where $p=Q_{\xi,c}(\Gamma\le h)$.
  \end{enumerate}
\end{lm}
L\'evy forests have been introduced as genealogical forests of continuous-state branching processes \cite{LGLJ,DuLG02,DuLG05,DuWi07}, in the sense of a variety of limit theorems. Before introducing L\'evy forests and L\'evy trees rigorously, let us add that a $(\psi;\varrho)$-L\'evy forest $\cF$ consists of infinitely many L\'evy trees and can be written as concatenation $\circledast_{i\in I}(\cT_i(0),d_i,\rho_i)$, where $\sum_{i\in I}\delta_{(\cT_i,d_i,\rho_i)}$ is a Poisson random measure with intensity measure $\int_{[0,\infty)}x\bN_{\psi}\varrho(dx)$, where $\bN_\psi$ is the $\sigma$-finite measure on $\bT$ that describes a single L\'evy tree, see \cite{DuWi07,DuWi12}. L\'evy trees (and hence similarly L\'evy forests) can also be characterised by their branching property at fixed heights: roughly, under $\bN_\psi$, for each $a>0$, conditionally given ${\rm Blw}(\cT,a)$, the forest ${\rm Abv}(\cT,a)$ is the concatenation $\circledast_{i\in I(a)}(\cT_i(a),d_{i,a},\rho_{i,a})$ of a Poisson point process with intensity measure $\int_{[0,\infty)}x\bN_{\psi}\varrho_a(dx)$ for some distribution $\varrho_a$, see \cite{Weill,DuWi12}. For the purpose of this paper, it will be most convenient to introduce L\'evy forests and L\'evy trees via their $h$-erasures, as was established in \cite[Theorems 3.16, 3.18 and 3.20]{DuWi12}.

\begin{defn}[L\'evy forests \cite{DuWi12}]\label{Levyfor}\rm A $\bT$-valued random variable $\cF$ is a L\'evy forest
  if $R^h(\cF)$ is a Galton-Watson real forest for all $h>0$, and if ${\rm n}(\cF)=\infty$ with positive probability. 
  Specifically, $\cF$ is a $(\psi;\varrho)$-L\'evy forest if $R^h(\cF)$ is a ${\rm GW}(\xi^{h,\psi},c^{h,\psi};\mu^{h,\psi})$-real forest for all $h>0$, where 
  $c^{h,\psi}\!=\!\psi^\prime(\eta(h))$ and $\xi^{h,\psi}$ and $\mu^{h,\psi}$ have generating functions\vspace{-0.1cm}
  $$g_{\xi^{h,\psi}}(s)=s+\frac{\psi((1-s)\eta(h))}{\eta(h)\psi^\prime(\eta(h))}
	\qquad\mbox{and}\qquad 
    g_{\mu^{h,\psi}}(s)=\int_{[0,\infty)}e^{-x(1-s)\eta(h)}\varrho(dx),$$
  where $\int_{\eta(h)}^\infty du/\psi(u)=h$. The distribution of a $(\psi;\varrho)$-L\'evy forest on $\bT$ is denoted by 
  $P_\psi^\varrho$. For a metric space $(X,\Delta_X,\rho_X)$, we also refer to a $\bT_X$-valued random variable as a L\'evy forest if its isometry class has distribution $P_\psi^\varrho$. 
\end{defn}

\begin{prop}[L\'evy trees \cite{DuWi07,DuWi13}]\label{levytree} $\!\!\!\;$For every branching mechanism {\rm(\ref{brmech})} satisfying {\rm(\ref{grey})} there is a $\sigma$-finite
  measure $\bN_\psi$ on $\bT$ with the following properties. We have $\bN_\psi({\rm n}\!\neq\! 1)=0$, and the concatenation $\circledast_{i\in I}(T_i,d_i,\rho_i)$ of the points $(T_i,d_i,\rho_i)$, $i\!\in\!I$, of a Poisson random measure on $\bT$ with intensity measure $\int_{[0,\infty)}x\bN_\psi\varrho(dx)$ is a $(\psi;\varrho)$-L\'evy forest for any distribution $\varrho$ on $[0,\infty)$. We have $\bN_\psi(R^h\in\cdot\,|\Gamma\!>\!h)=Q_{\xi^{h,\psi},c^{h,\psi}}$, $\bN_\psi(\Gamma\!>\!h)\!=\!\eta(h)$, and there is a family of regular conditional probability measures $\bN_\psi(\,\cdot\,|\Gamma\!=\!h)$, $h\!>\!0$, such that $\frac{1}{\eta(h)}\int_{h^\prime}^\infty\bN_\psi(\,\cdot\,|\Gamma\!=\!h)|\eta^\prime(h)|dh=\bN_\psi(\,\cdot\,|\Gamma\!>\!h^\prime)$.  
\end{prop} 
By \cite[Theorem 4.7]{DuLG05}, the limit $\overline{\omega}(\{v\})=\lim_{h\downarrow 0}{\rm n}(v,R^h(\cF))/\eta(h)$ exists for all $v\in\cF$ a.s.\ for any $(\psi;\varrho)$-L\'evy forest. Since $\eta(h)\rightarrow\infty$ as $h\downarrow 0$, this limit can only be non-zero if ${\rm n}(v,\cF)=\infty$, hence $\overline{\omega}$ is supported by the root and branch points with infinite multiplicity. Following Miermont \cite{Mie05}, we refer to $\overline{\omega}(\{v\})$ as the \em width of $v$\em, and to the atomic measure $\overline{\omega}$ as the \em width measure \em of $\cF$. For a $(\psi;\varrho)$-L\'evy forest $\circledast_{i\in I}(T_i,d_i,\rho_i)$ as in Proposition \ref{levytree}, $\overline{\omega}(\{\rho\})=\lim_{h\downarrow 0}\frac{1}{\eta(h)}\#\{i\in I\colon\Gamma(T_i)>h\}$ has distribution $\varrho$, by the Strong Law of Large Numbers for Poisson processes. 

\begin{thm}[Reconstruction, \cite{DuWi13}]\label{reconstruction} For $(\psi;\varrho)$, $h\!>\!0$, $\xi^{h,\psi}, c^{h,\psi}, \mu^{h,\psi}$ as in Definition \ref{Levyfor}, consider a $\bT_X$-valued ${\rm GW}(\xi^{h,\psi}\!, c^{h,\psi}\!; \mu^{h,\psi})$-real forest $\cF_*$. Conditionally given \nolinebreak $\cF_*$,  
\begin{itemize}
\item let $\cP_1\!=\!\sum_{v\in I}\delta_{(v ,T_v )}$ be a Poisson random measure  
on $\cF_*\! \times\! \bT$ with intensity measure $\ell\!\times\!\bN^h_\psi$, where \vspace{-0.2cm}
$$ \bN^h_\psi = 2\beta\bN_\psi (\, \cdot\, ; \Gamma \le h) + \int_{(0, \infty)} x P_{\psi}^{\delta_x} (\, \cdot \, ; \Gamma \le h) \pi (dx),\vspace{-0.1cm}$$

\item independently 
for each $v\!\in\!{\rm Br}(\cF_* )$ with $m\!:=\!{\rm n} (v, \cF_*)\!-\!1$, let \vspace{-0.2cm}
$$T_v\sim\frac{1}{\lvert 
\psi^{(m)} (\eta(h))\rvert} \left( 
2\beta 1_{\{ j= 2 \}} \delta_{\{v\}}  + \int_{(0, \infty)} x^m P_\psi^{\delta_x} (\, \cdot \, ; \Gamma \leq h) \pi (dx)\right),\vspace{-0.1cm}$$ 
	
\item independently for each $v\!\in\!{\rm Lf}(\cF_*)$, let $T_v\!\sim\!\bN_\psi (\, \cdot \, |\Gamma = h)$, \vspace{-0.1cm}

\item and, independently, for the root $\rho$ of $\cF_*$ with $m:={\rm n}(\cF_*)$ subtrees in $\cF_*$, let\vspace{-0.1cm} 
$$T_\rho\sim\frac{1}{\lvert 
\cL^{(m)} (\eta(h))\rvert} \int_{[0, \infty)} x^m P_{\psi}^{\delta_x} (\, \cdot \, ;\Gamma \le h)\varrho (dx) 
\; ,\vspace{-0.1cm}$$
where $\cL^{(m)}$ is the $m$-th derivative of $\cL(\theta)= \int_{[0, \infty)} e^{-\theta x}\varrho (dx)  $, the Laplace transform of $\varrho$. 
\end{itemize}
We denote by $\cF$ the tree obtained by grafting on $\cF_*$ the tree $T_v$ at $v\in\cF_*$, for all $v\in I\cup{\rm Br}(\cF_*)\cup{\rm Lf}(\cF_*)\cup\{\rho\}$. 
Then, a.s.\ $R^{h} (\cF)=\cF_*$ and $\cF$ is a $(\psi;\varrho)$-L\'evy forest. 
\end{thm}
\begin{cor}\label{coromega} Let $\cF$ be a $(\psi;\varrho)$-L\'evy forest $\cF$ with width measure $\overline{\omega}$, and $h\!>\!0$. Then conditionally given
  $R^h(\cF)$, the restriction of $\overline{\omega}$ to $R^h(\cF)$ is a random measure. Its conditional distribution given $R^h(\cF)$ is, as follows. We have $\overline{\omega}(\{v\})\!=\!W_v$, $v\!\in\!I_h\cup{\rm Br}(R^h(\cF))\cup\{\rho\}$,  where independently, 
  \begin{itemize}
  \item $\cP=\sum_{v\in I_h}\delta_{(v,W_v)}$ is a Poisson random measure on $R^h(\cF)\times(0,\infty)$ with intensity measure $\ell(dv)\times xe^{-x\eta(h)}\pi(dx)$,
  \item $W_v$ has Laplace transform $\psi^{(m)}(\eta(h)\!+\!\theta)/\psi^{(m)}(\eta(h))$ if $m\!=\!{\rm n}(v,R^h(\cF))\!-\!1$,  $v\!\in\!{\rm Br}(R^h(\cF))$,
  \item $W_\rho$ has Laplace transform $\cL^{(m)}(\eta(h)\!+\!\theta)/\cL^{(m)}(\eta(h))$ if $m\!=\!{\rm n}(R^h(\cF))$. 
  \end{itemize}\pagebreak[2]
\end{cor}
\begin{pf} Instead of $(R^h(\cF),\cF)$, we may consider $(\cF_*,\cF)$ as in Theorem \ref{reconstruction}. As noted just above Theorem \ref{reconstruction}, $\frac{1}{\eta(h^\prime)}\#\{i\!\in\! I\colon\Gamma(T_i)\!>\!h^\prime\}\!\rightarrow\! x$ as $h^\prime\!\rightarrow\! 0$, $P^{\delta_x}_\psi$-a.s. As $P^{\delta_x}_\psi(\Gamma\le h)\!=\!e^{-x\eta(h)}$, the distribution of $T_\rho$ yields $W_\rho\sim x^me^{-x\eta(h)}\varrho(dx)/|\cL^{(m)}(\eta(h))|$, conditionally given ${\rm n}(\cF_*)=m$, with Laplace transform as claimed. The same argument, with $\varrho$ and $\cL$ replaced by $\pi$ and $\psi$, yields the distributions $(1_{\{m=2\}}2\beta\delta_0(dx)+x^me^{-x\eta(h)}\pi(dx))/|\psi^{(m)}(\eta(h))|$ of $W_v$, $v\in{\rm Br}(R^h(\cF))$, from the distributions of $T_v$, as well as the intensity measure of $\cP$, by standard mapping of the Poisson random measure $\cP_1$ under the map that projects $T_v$ onto the width of its root $v$.  
\end{pf}

\subsection{Pruning at edges of GW trees and Aldous-Pitman pruning of L\'evy forests}\label{secprunedge}

Consider a pruning process $(\tau^E(\theta),\theta\ge 0)$ as defined in the introduction following \cite{AP98}, constructed from a Galton-Watson tree $\tau=\tau^E(0)\sim{\rm GW}(\xi)$ and independent pruning times $M_e\sim{\rm Exp}(1)$, $e\in E(\tau)$. To study the convergence of such pruning processes, we will need to understand finite-dimensional marginals, which we will represent in the product space $\bT^k$. For times $0\!=\!\theta_1\!<\!\theta_2\!<\!\cdots\!<\!\theta_k\!<\!\theta_{k+1}\!=\!\infty$, the finite-dimensional vector $(\tau^E(\theta_j),1\le j\le k)$ is governed by discretised pruning times. For $1\le j\le k$, we write $J_e=j$ if $M_e\in(\theta_{j},\theta_{j+1}]$, i.e.\ if edge $e$ is pruned between times $\theta_j$ and $\theta_{j+1}$.
We denote the distribution of $J_e$ by $q_j=\bP(J_e\!=\!j)=e^{-\theta_j}\!-\!e^{-\theta_{j+1}}$, $1\!\le\! j\!\le\! k$. In particular, $q_k=e^{-\theta_k}$ is the probability of no pruning before time $\theta_k$. For a vertex $v\in V(\tau)$ with $\#E_v(\tau)=m$ children, the edges are pruned according to independent $J_e$, $e\in E_v(\tau)$. Therefore, the numbers  $N_j=\#E_v(\tau^E(\theta_j))-\#E_v(\tau^E(\theta_{j+1}))$ of edges present in $\tau^E(\theta_j)$, but not in $\tau^E(\theta_{j+1})$, and $N_k=\#E_v(\tau^E(\theta_k))$ present in $\tau^E(\theta_k)$, form a multinomially distributed random vector $(N_1,\ldots,N_k)$ with probability function \vspace{-0.1cm}
$$\bP(N_1=n_1,\ldots,N_k=n_k)=\frac{m!}{n_1!\cdots n_k!}\prod_{j=1}^kq_j^{n_j},\qquad n_1,\ldots,n_k\in\bN,\ n_1+\cdots+n_k=m.\vspace{-0.1cm}$$
As a consequence, $\#E_v(\tau^E(\theta_j))=N_j+\cdots+N_k$, $1\le j\le k$, have joint probability function \vspace{-0.1cm}
$$\bP(N_j+\cdots+N_k=i_j,1\le j\le k)=\frac{i_1!q_k^{i_k}}{i_k!}\prod_{j=1}^{k-1}\frac{q_j^{i_j-i_{j+1}}}{(i_j-i_{j+1})!},\qquad m=i_1\ge i_2\ge\cdots\ge i_k\ge 0.\vspace{-0.1cm}
$$
In a ${\rm GW}(\xi)$-real tree, the first branch point (or leaf) above the root is at a ${\rm geom}(1-\xi(1))$ height. Each edge is pruned with probability $1-q_k=1-e^{-\theta_k}$, so the minimum of the height of the first branch point of $\tau^E(0)$ and the first pruning height of $\tau^E(\theta_k)$ is ${\rm geom}(1-\xi(1)q_k)$. Furthermore, this first leaf or branch point is a leaf of $\tau^E(\theta_k)$ due to pruning below the first leaf or branch point of $\tau^E(0)$ with probability $(1-q_k)\xi(1)/(1-\xi(1)q_k)$, indeed pruning occurred between $\theta_j$ and $\theta_{j+1}$ with probability $q_j\xi(1)/(1-\xi(1)q_k)$, $1\le j\le k-1$. Finally, with probability $(1-\xi(1))/(1-\xi(1)q_k)$, there is no pruning of $\tau^E(\theta_k)$ below the first leaf or branch point of $\tau^E(0)$. More precisely, this first branch point (or leaf) is with $m\neq 1$ children with probability $\xi(m)/(1-\xi(1)q_k)$. 


For $1\le j\le k$, projections $\pi_j\colon\bT^k\rightarrow\bT$ are induced by $\pi_j(T_1,\ldots,T_k)=T_j$ and inclusions $\iota^{j,k}\colon\bT^{j}\!\rightarrow\!\bT^{k}$ by $\iota^{j,k}(T_1,\ldots,T_j)\!=\!(T_1,\ldots,T_j,\{\rho\},\ldots,\{\rho\})$. Then $\pi_j$ and $\iota^{j,k}$ are continuous. For a measure $Q$ on $\bT^{j}$, denote by $\iota^{j,k}_*Q$ the pushforward of $Q$ under $\iota^{j,k}$, so $\iota^{j,k}_*Q$ is the distribution on $\bT^k$ of a $Q$-distributed random variable in $\bT^j$ with trivial components $j\!+\!1,\ldots,k$ added.

We write $D_j=D\circ\pi_j$ and ${\rm n}_j={\rm n}\circ\pi_j$, $1\le j\le k$, and then introduce functions $\underline{{\rm n}},\underline{D},\underline{\vartheta},\underline{{\bf k}}$ on $\bT^k$ as $\underline{\rm n}=({\rm n}_1,\ldots,{\rm n}_k)$, $\underline{D}=\min\{D_1,\ldots,D_k\}$, $\underline{\vartheta}=({\rm Abv}(\pi_1,\underline{D}),\ldots,{\rm Abv}(\pi_k,\underline{D}))$, $\underline{\fk}=\underline{\rm n}\circ\underline{\vartheta}$. In particular, note that $\underline{\vartheta}$ and $\underline{\fk}$ refer to the subtrees above the same height $\underline{D}=\min\{D_1,\ldots,D_k\}$, not above the individual first branch points at heights $D_1,\ldots,D_k$. Note also that $\underline{\fk}$ takes values in vectors with not all entries equal to 1, but some may be equal to 1. It is reasonable to give the following definition for more general pruning distributions $(q_j,1\le j\le k)$.\pagebreak[2]

\begin{defn}[Pruning at edges, GW($\xi$)]\label{dfprunedge}\rm We fix an offspring distribution $\xi$ with $\xi(1)<1$ and consider a pruning distribution $q\!=\!(q_j,1\!\le\! j\!\le\! k)$. For $k\!=\!1$ we set $\widehat{Q}\!=\!Q_\xi$ as the distribution of an unpruned ${\rm GW}(\xi)$-real tree. For $k\!\ge\! 2$, a $q$-\em pruning at edges \em of a ${\rm GW}(\xi)$-real tree is a $\bT^{k}$-valued random variable $(\cT_1,\ldots,\cT_k)$, whose distribution $\widehat{Q}_{q_1,\ldots,q_{k-1}}$ is such that under $\widehat{Q}_{q_1,\ldots,q_{k-1}}$,
\begin{enumerate}\item[1.$\,\,$] $\underline{D}$ and $(\underline{\fk},\underline{\vartheta})$ are independent,
  \item[2.$^{\!\!\!\,E}$] $\underline{D}\sim{\rm geom}(1-\xi(1)q_k)$,\vspace{-0.7cm}
  \item[3.$^{\!\!\!\,E}$] $\underline{\fk}\sim\widetilde{\xi}$ where 
$\displaystyle\widetilde{\xi}(i_1,\ldots,i_k)\!=\!\!\left\{\!\!\begin{array}{lr}\displaystyle\frac{\xi(1)q_j}{1-\xi(1)q_k}&\hspace{-3.55cm}\mbox{if $i_1\!=\!\cdots\!=\!i_j\!=\!1\!>\!i_{j+1}\!=\!\cdots\!=\!i_k\!=\!0$, $1\!\le\! j\!\le\! k\!-\!1$,}\\
\displaystyle\frac{\xi(i_1)}{1-\xi(1)q_k}\frac{i_1!q_k^{i_k}}{i_k!}\prod_{j=1}^{k-1}\frac{q_j^{i_j-i_{j+1}}}{(i_j-i_{j+1})!}&\mbox{if $1\!\neq\! i_1\!\ge\!\cdots\!\ge\! i_k\!\ge\! 0$,}\end{array}\right.\hspace{-0.2cm}$
where the first line reflects pruning on an edge below the first branch point of $\cT_1$ and the second line pruning of a multinomially distributed number of the $i_1\ge 2$ edges just above the first branch point of $\cT_1$, also including $\widetilde{\xi}(0,\ldots,0)=\xi(0)/(1-\xi(1)q_k)$, 
  \item[4.$^{\!\!\!\,E}$] and conditionally given $\{\underline{\fk}=(i_1,\ldots,i_k)\}$, we have $\underline{\vartheta}$ as a concatenation of $i_j-i_{j+1}$ trees with distribution  $\iota^{j,k}_*\widehat{Q}_{q_1,\ldots,q_{j-1}}$, $1\le j\le k-1$, and $i_k$ trees with distribution $\widehat{Q}_{q_1,\ldots,q_{k-1}}$.
\end{enumerate}
\end{defn}
Note that the first formula in 3.$^{\!\!\!\,E}$ is actually a special case of the second formula here. Also, it is useful to leave the no-pruning probability $q_k=1-q_1-\cdots-q_{k-1}$ implicit in notation $\widehat{Q}_{q_1,\ldots,q_{k-1}}$ since in 4.$^{\!\!\!\,E}$ the no-pruning probability is not $q_j$, but $1-q_1-\cdots-q_{j-1}=q_j+\cdots+q_k$, as required.
\begin{defn}\rm A pruning process $(\cT^E(\theta),\theta\ge 0)$ in the sense of Definition \ref{prunproc} is called a ${\rm GW}(\xi)$-pruning process with \em pruning at edges \em if (the isometry classes of) $\cT^E(\theta_1),\ldots,\cT^E(\theta_k)$ have joint distribution $\widehat{Q}_{q_1,\ldots,q_{k-1}}$ for all $0\!=\!\theta_1\!<\!\theta_2\!<\!\cdots\!<\!\theta_k\!<\!\theta_{k+1}\!=\!\infty$, where $q_j=e^{-\theta_j}-e^{-\theta_{j+1}}$,
$1\!\le\!j\!\le\!k$.
\end{defn}
\begin{lm}\label{lm1} For each offspring distribution $\xi$, there is a unique family of distributions $\widehat{Q}_{q_1,\ldots,q_{k-1}}$, $q_i\ge 0$, $1\le i\le k-1$, $q_1+\cdots+q_{k-1}<1$, $k\ge 1$, that satisfies the recursive Definition \ref{dfprunedge}.
\end{lm}
\begin{pf} The proof of \cite[Lemma 2.15]{DuWi12} for the unpruned case $k=1$ of ${\rm GW}(\xi)$-real trees can be adapted, using induction on $k$, for any fixed sequence $q_i\ge 0$ with $q_1+\cdots+q_j<1$, $j\ge 1$.
\end{pf}
Definition \ref{dfprunedge} decomposes the distribution $\widehat{Q}_{q_1,\ldots,q_{k-1}}$ into quantities amenable to taking limits. The following definition applies to L\'evy forests, as well as to ${\rm GW}(\xi,c)$-real trees/forests:
\begin{defn}[Aldous-Pitman pruning, \cite{AP98b}]\rm A pruning process 
  $(\cT^{\rm AP}(\theta),\theta\!\ge\! 0)$ of $T\in\bT_X$ is called an \em Aldous-Pitman pruning process of $T$ \em if
  $\cP=\cP^{\rm AP}$ is a Poisson random measure on $(0,\infty)\times T$ with intensity measure $d\theta\ell(dv)$, where $\ell$ is
  the length measure on $T$.
\end{defn}\pagebreak[2]

\begin{prop}\label{dfprunedgec} For an Aldous-Pitman pruning process $(\cT^{\rm AP}(\theta),\theta\ge 0)$ of a ${\rm GW}(\xi,c)$-real tree, with $\xi(1)\!=\!0$, the distributions $\widehat{Q}_{\theta_2,\ldots,\theta_k}$ on $\bT^k$ of
$(\cT^{\rm AP}(\theta_1),\ldots,\cT^{\rm AP}(\theta_k))$, $0\!=\!\theta_1\!<\!\theta_2\!<\!\cdots\!<\!\theta_k$, are uniquely  determined by $\widehat{Q}\!=\!Q_{\xi,c}$ for $k\!=\!1$, and the following recursive rule for $k\!\ge\! 2$. Under $\widehat{Q}_{\theta_2,\ldots,\theta_k}$,
\begin{enumerate}\item[1.$\,\,$] $\underline{D}$ and $(\underline{\fk},\underline{\vartheta})$ are independent,
  \item[$\overline{\rm 2}.^{\!\!\!\,E}$] $\underline{D}\sim{\rm Exp}(c+\theta_k)$,\vspace{-0.6cm}
  \item[$\overline{\rm 3}.^{\!\!\!\,E}$] $\underline{\fk}\sim\widetilde{\xi}$ where 
$\displaystyle\widetilde{\xi}(i_1,\ldots,i_k)\!=\!\!\left\{\!\!\begin{array}{ll}\displaystyle\frac{\theta_{j+1}-\theta_j}{c+\theta_k}&\mbox{if $i_1\!=\!\cdots\!=\!i_j\!=\!1\!>\!i_{j+1}\!=\!\cdots\!=\!i_k\!=\!0$, $1\!\le\! j\!\le\! k\!-\!1$,}\\[0.4cm]
\displaystyle\frac{c\xi(i_1)}{c+\theta_k}&\mbox{if $1\!\neq\! i_1\!=\!\cdots\!=\! i_k$,}\end{array}\right.\hspace{-0.2cm}$
where the first line reflects pruning on the branch below the first branch point of $\cT^{\rm AP}(0)$ and the second line no pruning below the first branch point of $\cT^{\rm AP}(0)$,
  \item[$\overline{\rm 4}.^{\!\!\!\,E}$] and conditionally given $\{\underline{\fk}=(i_1,\ldots,i_k)\}$, we have $\underline{\vartheta}$ as a concatenation of $i_j-i_{j+1}$ trees with distribution  $\iota^{j,k}_*\widehat{Q}_{\theta_2,\ldots,\theta_j}$, $1\le j\le k-1$, and $i_k$ trees with distribution $\widehat{Q}_{\theta_2,\ldots,\theta_k}$.\pagebreak[2] 
\end{enumerate}
\end{prop}
\begin{pf} In an Aldous-Pitman pruning process up to time $\theta_k$, the Poisson random measure has intensity measure $\theta_k\ell$. For a ${\rm GW}(\xi,c)$-real tree, this holds branch by branch, cutting each ${\rm Exp}(c)$ branch into a 
${\rm geom}(c/(c+\theta_k))$ number of ${\rm Exp}(c+\theta_k)$ parts. In particular, the probability of pruning on the first branch is $\theta/(c+\theta_k)$. By standard thinning, this further splits into $(\theta_{j+1}-\theta_j)/(c+\theta_k)$, $1\le j\le k-1$, for pruning of $\cT^{\rm AP}(\theta_{j+1})$, but not $\cT^{\rm AP}(\theta_j)$. By the same reasoning as in the setup of Definition \ref{dfprunedge}, we deduce the distribution of $(\underline{D},\underline{\fk},\underline{\vartheta})$ under $\widehat{Q}_{\theta_2,\ldots,\theta_k}$ from the Definition of ${\rm GW}(\xi,c)$-real trees and the independence and identical distribution of the Poisson random measure on subtrees. Uniqueness is obtained as indicated in Lemma \ref{lm1}.
\end{pf}
\begin{cor} For an Aldous-Pitman pruning process $(\cF^{\rm AP}(\theta),\theta\ge 0)$ of a ${\rm GW}(\xi,c;\mu)$-real forest, $(\cF^{\rm AP}(\theta_1),\ldots,\cF^{\rm AP}(\theta_k))$ has joint distributions $\widehat{P}_{\theta_2,\ldots,\theta_k}^\mu:=\sum_{\ell\ge 0}\mu(\ell)\widehat{Q}_{\theta_2,\ldots,\theta_k}^{\circledast\ell}$.  
\end{cor}
Note that pruning on every edge only occurs at exponentially distributed heights, so no thinning occurs at the first branch point. Hence, either $\underline{\vartheta}\sim \widehat{Q}_{\theta_2,\ldots,\theta_j}$ or $\underline{\vartheta}\sim \widehat{Q}_{\theta_2,\ldots,\theta_k}^{\circledast i_1}$, no concatenations of different $\widehat{Q}_{\theta_2,\ldots,\theta_j}$, $1\!\le\! j\!\le\! k$, occur. 
Recall that L\'evy forests $\cF$ were defined via $R^h(\cF)$, $h>0$.
\begin{prop}\label{edgeerase} A $\bT_X$-valued random process $(\cF(\theta),\theta\!\ge\! 0)$ is an Aldous-Pitman pruning process of a $(\psi;\varrho)$-L\'evy forest $\cF\!=\!\cF(0)$ if and only if for all $h\!>\!0$, the process $(\cF(\theta)\cap R^h(\cF),\theta\!\ge\! 0)$ is an Aldous-Pitman pruning process of a ${\rm GW}(\xi^{h,\psi},c^{h,\psi};\mu^{h,\psi})$-real forest, with $\xi^{h,\psi},c^{h,\psi},\mu^{h,\psi}$ as in Definition 
  \ref{Levyfor}. In particular, $\cF(\theta)=\cF^{\rm AP}(\theta)$ is a $(\widehat{\psi}_\theta;\varrho)$-L\'evy forest, where $\widehat{\psi}_\theta(u)=\psi(u)+\theta u$.\pagebreak[2]
\end{prop}
\begin{pf} For the ``if'' part, note that $\cF$ is a $(\psi;\varrho)$-L\'evy forest by definition, and that the Aldous-Pitman 
  pruning processes $(\cF(\theta)\cap R^h(\cF),\theta\ge 0)$, $h>0$, yield a Poisson random measure $\cP$ with intensity measure
  $d\theta\ell(dv)$ on $(0,\infty)\times\bigcup_{h>0}R^h(\cF)=[0,\infty)\times\cF\setminus{\rm Lf}(\cF)$, but since $\ell$ does 
  not charge ${\rm Lf}(\cF)$ and adding further point masses on ${\rm Lf}(\cF)$ to $\cP$ would not change the distribution of the
  associated pruning process, this identifies $(\cF(\theta),\theta\ge 0)$ as an Aldous-Pitman pruning process. The 
  ``only if'' part is straightforward. Finally, note that we obtain from Proposition \ref{dfprunedgec} for $k=2$ that $\cF^{\rm AP}(\theta)\cap R^h(\cF^{\rm AP}(0))\sim{\rm GW}(\xi_\theta^{h,\psi},c_\theta^{h,\psi};\mu^{h,\psi})$, where
$c_\theta^{h,\psi}=c^{h,\psi}+\theta=\psi_\theta^\prime(\eta(h))$ and\vspace{-0.1cm}
$$g_{\xi_\theta^{h,\psi}}(s)=\frac{c^{h,\psi}}{c^{h,\psi}+\theta}g_{\xi^{h,\psi}}(s)+\frac{\theta}{c^{h,\psi}+\theta}=s+\frac{\widehat{\psi}_\theta((1-s)\eta(h))}{\eta(h)\widehat{\psi}_\theta^\prime(\eta(h))},$$
  and these are offspring distributions as they appear when erasing a $(\widehat{\psi}_\theta;\varrho)$-L\'evy forest. 
\end{pf}

We state a related result for Galton-Watson trees with pruning at edges. Its proof is easier and left to the reader.

\begin{prop}\label{edgeerasediscr} For a ${\rm GW}(\xi;\mu)$-pruning process $(\cF(\theta),\theta\!\ge\! 0)$ with pruning at edges and $\xi^h,\mu^h$ from {\rm(\ref{xihmuh})},  $(\cF(\theta)\cap R^h(\cF(0)),\theta\!\ge\! 0)$ is a ${\rm GW}(\xi^h;\mu^h)$-pruning process with pruning at edges.
\end{prop}


\subsection{Pruning at branch points of GW trees, and Abraham-Delmas pruning}\label{sectpruneGW}

Consider a pruning process $(\tau^B(\theta),\theta\!\ge\! 0)$ as defined in the introduction following \cite{ADH12}, constructed from a Galton-Watson tree $\tau\!=\!\tau^B(0)\!\sim\!{\rm GW}(\xi)$ and independent pruning times $M_v\!\sim\!{\rm Exp}(\#E_v(\tau)\!-\!1)$, $v\!\in\!{\rm Br}(\tau)$. Then no pruning occurs below the first branch point. We will later need more general pruning, which allows non-exponential pruning time distributions $H_m$, $m\!\ge\! 1$, for branch points with $\#E_v(\tau)\!=\!m$, and also pruning below the first branch point.

\begin{defn}[$H$-pruning, pruning at branch points, GW($\xi$)]\label{prbrpp}\rm $\;\!\!\!$Fix an offspring distribution $\xi$ with $\xi(1)<1$ and a family $H\!=\!(H_m,m\!\ge\! 1)$ of pruning time distributions on $(0,\infty]$. Given a $\bT_X$-valued ${\rm GW}(\xi)$-real tree $\cT$, an \em $H$-pruning process of a ${\rm GW}(\xi)$-real tree \em is a pruning process in the sense of Definition \ref{prunproc} that is associated with a point measure 
  $\cP\!=\!\sum_{v\in{\rm Br}^+(\cT)}\delta_{(\Theta_v,v)}$ with  
  ${\rm Br}^+(\cT)\!=\!\{v\!\in\!\cT\!\setminus\!(\{\rho\}\cup{\rm Lf}(\cT))\colon d(\rho,v)\!\in\!\bN\}$ as the set of
  branch points including further points between unit segments on branches and  
  $\Theta_v\!\sim\! H_{{\rm n}(v,\cT)-1}(d\theta)$, $v\!\in\!{\rm Br}^+(\cT)$. We write $(\cT^H(\theta),\theta\ge 0)$. If $H_m={\rm Exp}(m-1)$, $m\ge 2$, and 
  $H_1=\delta_\infty$ is the point mass at $\infty$, it is called a ${\rm GW}(\xi)$-pruning process with \em pruning at branch points \em and written $(\cF^B(\theta),\theta\ge 0)$. For ${\rm GW}(\xi;\mu)$-real forests, we similarly define pruning processes $(\cF^B(\theta),\theta\ge 0)$ with pruning at branch points and $H$-pruning processes $(\cF^H(\theta),\theta\ge 0)$.
\end{defn}

\begin{prop}\label{dfprunbranch} For an $H$-pruning process $(\cT^H(\theta),\theta\ge 0)$ of a ${\rm GW}(\xi)$-real tree, the distributions $\widetilde{Q}_{H;\theta_2,\ldots,\theta_k}$ of $(\cT^H(\theta_1),\ldots,\cT^H(\theta_k))$, $0\!=\!\theta_1\!<\!\theta_2\!<\!\cdots\!<\!\theta_k$, are uniquely determined by $\widetilde{Q}_H\!=\!Q_\xi$ for $k\!=\!1$ and the following recursive rule for $k\!\ge\! 2$. Under $\widetilde{Q}_{H;\theta_2,\ldots,\theta_k}$,
\begin{enumerate}\item[1.$\,\,$] $\underline{D}$ and $(\underline{\fk},\underline{\vartheta})$ are independent,\vspace{-0.1cm}
  \item[2.$^{\!\!\!\,B}$] $\underline{D}\sim{\rm geom}(1-\xi(1)H_1((\theta_k,\infty]))$,\vspace{-0.1cm}
  \item[3.$^{\!\!\!\,B}$] $\underline{\fk}\sim\widetilde{\xi}$ where 
$\displaystyle\widetilde{\xi}(i_1,\ldots,i_k)\!=\!\frac{\xi(i_1)H_{i_1}((\theta_j,\theta_{j+1}])}{1-\xi(1)H_1((\theta_k,\infty])}\vspace{0.1cm}$
if $i_1\!=\!\cdots\!=\!i_j\!>\!i_{j+1}\!=\!\cdots\!=\! i_k\!=\! 0$, $1\!\le\!j\!\le\!k$, but excluding the case $i_1\!=\!1$, $j\!=\!k$.
where the case $i_1=1$ reflects pruning at a single-child vertex below the first branch point of $\cT_1$ and the case $i_1\neq 1$ no pruning or pruning at the first branch point of $\cT_1$, and  $\widetilde{\xi}(0,\ldots,0)=\xi(0)/(1-\xi(1)H_1((\theta_k,\infty]))$, \vspace{-0.1cm}
  \item[4.$^{\!\!\!\,B}$] and conditionally given $\{\underline{\fk}=(i_1,\ldots,i_k)\}$, we have $\underline{\vartheta}$ as a concatenation of $i_j-i_{j+1}$ trees with distribution  $\iota^{j,k}_*\widetilde{Q}_{H;\theta_2,\ldots,\theta_j}$, $1\le j\le k-1$, and $i_k$ trees with distribution $\widetilde{Q}_{H;\theta_2,\ldots,\theta_k}$.\vspace{-0.1cm}\pagebreak[2]
\end{enumerate}
\end{prop}

\begin{defn}[$\overline{H}$-pruning, GW($\xi,c$)]\rm\label{prbrpp2} Let $\overline{H}_1$ be a measure on $(0,\infty)$, finite on bounded sets, $\overline{H}\!:=\!(\overline{H}_1;H_m,m\!\ge\! 2)$ and $c\!>\!0$. An \em $\overline{H}$-pruning process of a ${\rm GW}(\xi,c)$-real tree \em is a pruning \linebreak process of a $\bT_X$-valued ${\rm GW}(\xi,c)$-real tree $\cT$ associated with 
  $\cP=\cP_1+\sum_{v\in{\rm Br}(\cT)}\delta_{(\Theta_v,v)}$ where $\cP_1$ is a Poisson
  random measure on $(0,\infty)\times\cT$ with intensity measure $\overline{H}_1(d\theta)\ell(dv)$ and 
  $\Theta_v\sim H_{{\rm n}(v,\cT)-1}(d\theta)$, $v\in{\rm Br}(\cT)$. We write $(\cT^{\overline{H}}(\theta),\theta\ge 0)$, and  $(\cF^{\overline{H}}(\theta),\theta\ge 0)$ for ${\rm GW}(\xi,c;\mu)$-forests.
\end{defn}

\vspace{-0.3cm}

\begin{prop}\label{dfprunbranch2} For an $\overline{H}$-pruning process $(\cT^{\overline{H}}(\theta),\theta\ge 0)$ of a ${\rm GW}(\xi,c)$-real tree, with $\xi(1)=0$, the distributions $\widetilde{Q}_{\overline{H};\theta_2,\ldots,\theta_k}$ of $(\cT^{\overline{H}}(\theta_1),\ldots,\cT^{\overline{H}}(\theta_k))$, $0\!=\!\theta_1\!<\!\theta_2\!<\!\cdots\!<\!\theta_k$, are \vspace{-0.1cm} uniquely determined by $\widetilde{Q}_{\overline{H}}\!=\!Q_{\xi,c}$ for $k\!=\!1$ and the following recursive rule for $k\!\ge\! 2$. Under 
$\widetilde{Q}_{\overline{H};\theta_2,\ldots,\theta_k}$,
\begin{enumerate}\item[1.$\,\,$] $\underline{D}$ and $(\underline{\fk},\underline{\vartheta})$ are independent,\vspace{-0.1cm}
  \item[$\overline{\rm 2}.^{\!\!\!\,B}$] $\underline{D}\sim{\rm Exp}(c+\overline{H}_1((0,\theta_k]))$,\vspace{-0.7cm}
  \item[$\overline{\rm 3}.^{\!\!\!\,B}$] $\underline{\fk}\sim\widetilde{\xi}$ where 
$\displaystyle\widetilde{\xi}(i_1,\ldots,i_k)\!=\!\!\left\{\!\!\begin{array}{cl}\displaystyle\frac{\overline{H}_1((\theta_j,\theta_{j+1}])}{c+\overline{H}_1((0,\theta_k])}&\mbox{if $i_1\!=\!1$,}\\[0.4cm]
\displaystyle\!\frac{c\xi(i_1)H_{i_1}((\theta_j,\theta_{j+1}])}{c+\overline{H}_1((0	,\theta_k])}&\mbox{if $i_1\!\neq\!1$,}\end{array}\right.\!\!\!\!\vspace{0.1cm}i_1\!=\!\cdots\!=\!i_j\!>\!i_{j+1}\!=\cdots\!=\!i_k\!=\!0$, $1\!\le\!j\!\le\!k$, but excluding the case $i_1=1$, $j=k$, 
where the first line reflects pruning on the branch below the first branch point of $\cT_1$ and the second line no pruning or pruning at the first branch point of $\cT_1$, and $\widetilde{\xi}(0,\ldots,0)=c\xi(0)/(c+\overline{H}_1((0,\theta_k]))$,\vspace{-0.1cm}
  \item[$\overline{\rm 4}.^{\!\!\!\,B}$] and conditionally given $\{\underline{\fk}=(i_1,\ldots,i_k)\}$, we have $\underline{\vartheta}$ as a concatenation of $i_j-i_{j+1}$ trees with distribution  $\iota^{j,k}_*\widetilde{Q}_{\overline{H};\theta_2,\ldots,\theta_j}$, $1\le j\le k-1$, and $i_k$ trees with distribution $\widetilde{Q}_{\overline{H};\theta_2,\ldots,\theta_k}$.\pagebreak
\end{enumerate}   
\end{prop}
\begin{pfofprop2729} The recursive rules for $\widetilde{Q}_{H;\theta_2,\ldots,\theta_k}$ and $\widetilde{Q}_{\overline{H};\theta_2,\ldots,\theta_k}$ are straightforward. Uniqueness follows as indicated in Lemma \ref{lm1}.
\end{pfofprop2729}

\noindent The definition of Abraham-Delmas pruning processes depends on the width measure $\overline{\omega}$ supported by the branch points of a L\'evy forest $\cF$, see after Proposition \ref{levytree}. We denote by $\omega$ the restriction of $\overline{\omega}$ to $\cF\setminus\{\rho\}$. Following L\"ohr, Voisin and Winter \cite{LVW}, we define Abraham-Delmas' \cite{AD12} pruning processes as a special case of pruning processes driven by a more general $\sigma$-finite pruning measure $\nu$ on $\cF\setminus{\rm Lf}(\cF)$, i.e.\ $\nu$ is finite on compact subsets of $\cF\setminus{\rm Lf}(\cF)$ such as $[[\rho,v]]$, $v\!\in\!\cF\setminus{\rm Lf}(\cF)$.

\begin{defn}[$\nu$-pruning, Abraham-Delmas pruning]\rm Let $T\in\bT_X$ and $\nu$ a $\sigma$-finite measure on the skeleton 
  $T\setminus{\rm Lf}(T)$ of $T$. A pruning process $(\cT^\nu(\theta),\theta\ge 0)$ is called a \em $\nu$-pruning process of $T$ \em if it is associated with a Poisson random measure with intensity measure $d\theta\nu(dv)$. \pagebreak[2]
  
  Let $\cF$ be a $(\psi;\varrho)$-L\'evy forest and $\nu=\omega+2\beta\ell$, where $\omega$ and $\ell$ are the width and length 
  measures of $\cF$, and $\beta$ is the quadratic coefficient of $\psi$ in (\ref{brmech}). Then we refer to a $\nu$-pruning 
  process of $\cF$ as an \em Abraham-Delmas pruning process of $\cF$\em, and use notation $(\cF^{\rm AD}(\theta),\theta\ge 0)$.
\end{defn}

Note that Aldous-Pitman pruning processes of $T\in\bT_X$ are $\nu$-pruning processes for $\nu=\ell$, the length measure of $T$. We can obtain the analogue of Proposition \ref{edgeerase} and identify the pruning process $(\cF^{\rm AD}(\theta)\cap R^h(\cF^{\rm AD}(0)),\theta\ge 0)$ as $\nu^h$-pruning process of the ${\rm GW}(\xi^h,c^h;\mu^h)$, where $\nu^h$ is the restriction of $\nu=\omega+2\beta\ell$ to $R^h(\cF^{\rm AD}(0))$. However, this is less useful than the $\ell$-pruning process of Proposition \ref{edgeerase}, since unlike $\ell$, the measure $\nu^h$ is not an intrinsic measure that can be constructed from $R^h(\cF^{\rm AD}(0))$. Cf.\ \cite[Proposition 2.25]{LVW}. Instead, we conclude this section by providing an autonomous description of this pruning process as an $\overline{H}$-pruning process, and a similar result for the $h$-erasure of a ${\rm GW}(\xi)$-pruning process with pruning at branch points, which yields an $H$-pruning process. 

%
%
%
%

\begin{prop}\label{prop21} A $\bT_X$-valued random process $\cX:=(\cF(\theta),\theta\!\ge\! 0)$ is an Abraham-Delmas pruning process of a $(\psi;\varrho)$-L\'evy forest $\cF\!=\!\cF(0)$ if and only if the process $\cX^h:=(\cF(\theta)\cap R^h(\cF),\theta\!\ge\! 0)$ is an $\overline{H}$-pruning process of a ${\rm GW}(\xi^{h,\psi},c^{h,\psi};\mu^{h,\psi})$-real forest for all $h>0$, with $\xi^{h,\psi},c^{h,\psi},\mu^{h,\psi}$ as in Definition \ref{Levyfor} and\vspace{-0.3cm}  
    $$\overline{H}_1((0,\theta])=\psi^{\prime}(\eta(h)+\theta)-\psi^\prime(\eta(h)),\ \theta\!\ge\! 0,\quad\mbox{and}\  
      H_m((\theta,\infty])=\frac{\psi^{(m)}(\eta(h)+\theta)}{\psi^{(m)}(\eta(h))},\ \theta\!\ge\! 0,\ m\!\ge\! 2.$$
\end{prop}
\begin{rem}\rm
When $\eta(h)=\psi^{-1}(\lambda)$, the process $\cX^h=(\cF(\theta)\cap R^h(\cF),\theta\ge 0)$ has the same distribution as the pruning process $\cX_\lambda:=(\cF_\lambda(\theta),\theta\ge 0)$ based on Poisson sampling, studied in \cite[Section 6.1]{ADH13}. While convergence of $\cF_\lambda(\theta)$ to $\cF(\theta)$, in distribution as $\lambda\rightarrow\infty$, is easily obtained from \cite[Proposition 4.1 and Theorem 5.1]{ADH13}, Proposition \ref{prop21} here implies Skorohod convergence $\cX_\lambda\rightarrow\cX$, in distribution as $\lambda\rightarrow\infty$, since $d_{\rm Sk}(\cX^h,\cX)\le h$ as in the proof of Corollary \ref{cor25} yields $\cX^h\rightarrow\cX$ almost surely as $h\downarrow 0$, while 
$d_{\rm Sk}(\cX_\lambda,\cX)$ does not allow deterministic a.s.\ bounds.
\end{rem}
\begin{pfofprop21} For the ``only if'' part, Corollary \ref{coromega} provides the conditional distribution of $\omega$ given $R^h(\cF)$. For each realisation of $(R^h(\cF),\omega)$, the process $(\cF^{\rm AD}(\theta)\cap R^h(\cF),\theta\ge 0)$ is a $\nu^h$-pruning process, where $\nu^h$ is the restriction of $\nu=\omega+2\beta\ell$ to $R^h(\cF)$. Specifically, each  $v\in R^h(\cF)$ with $W_v=\omega(\{v\})>0$ has an independent pruning time $M_v\sim{\rm Exp}(W_v)$, and further pruning occurs according to an independent Poisson random measure $\sum_{v\in I_\ell}\delta_{(M_v,v)}$ with intensity measure $d\theta 2\beta\ell(dv)$. Given the distribution of $W_v$ in branch points of $R^h(\cF)$ from the proof of Corollary \ref{coromega}, the conditional distribution of $M_v$ given only $R^h(\cF)$ is a mixed exponential distribution with survival function\vspace{-0.1cm}
$$\frac{1_{\{m=2\}}2\beta}{|\psi^{(m)}(\eta(h))|}+\int_{(0,\infty)}e^{-x\theta} \frac{x^me^{-x\eta(h)}}{|\psi^{(m)}(\eta(h))|}\pi(dx)=\frac{\psi^{(m)}(\eta(h)+\theta)}{\psi^{(m)}(\eta(h))}\vspace{-0.1cm}$$
if $m={\rm n}(v,R^h(\cF))-1\ge 2$. For $m=1$, the atoms of $\omega$ of size $W_v$ at $v$ on the branches follow a Poisson random measure $\sum_{v\in I_\omega}\delta_{(W_v,v)}$, and the association of $M_v\sim{\rm Exp}(W_v)$ pruning times is a marking operation for this Poisson random  measure. By the mapping theorem for Poisson random measures, $\sum_{v\in I_\omega}\delta_{(M_v,v)}$ is a Poisson random measure with intensity measure $\ell(dv)\times\int_{(0,\infty)}xe^{-x\theta}d\theta xe^{-x\eta(h)}\pi(dx)$. By superposition of Poisson random measures, the
$\nu^h$-pruning process has pruning on branches according to a Poisson random measure with intensity measure\vspace{-0.1cm} 
$$\left(2\beta d\theta+\int_{(0,\infty)}x^2e^{-x(\eta(h)+\theta)}\pi(dx)d\theta\right)\ell(dv)=\psi^{\prime\prime}(\eta(h)+\theta)d\theta\ell(dv).\vspace{-0.1cm}$$
  For the ``if'' part, first note that in the ``only if'' setting, 
  $\cX^h:=(\cF(\theta)\cap R^h(\cF),\theta\ge 0)\rightarrow(\cF(\theta),\theta\ge 0)=:\cX$ a.s., since 
  $d_{\rm Sk}(\cX^h,\cX)\le h$, as in the proof of Corollary \ref{cor25}. Now consider \em any \em forest $\widetilde{\cF}$ such
  that $\widetilde{\cX}^h:=(\widetilde{\cF}(\theta)\cap R^h(\widetilde{\cF}),\theta\ge 0)\ed\cX^h$ for all $h>0$. 
  Then $\widetilde{\cX}^h\!\rightarrow\!\widetilde{\cX}\!:=\!(\widetilde{\cF}(\theta),\theta\!\ge\! 0)$ a.s.\ and 
  $\widetilde{\cX}^h\!\convd\!\cX$ in distribution. Hence $\widetilde{\cX}\!\ed\!\cX$, as required.
\end{pfofprop21}
We deduce from this new result marginal distributions obtained in \cite{AD12,ADV} using different methods.
\begin{cor}[\cite{AD12,ADV}] Let $(\cF^{\rm AD}(\theta),\theta\!\ge\! 0)$ be an Abraham-Delmas pruning process of a $(\psi;\varrho)$-L\'evy forest. Then $\cF^{\rm AD}(\theta)$ is a $(\widetilde{\psi}_\theta;\varrho)$-L\'evy forest, $\theta\!\ge\! 0$, where  $\widetilde{\psi}_\theta(u)=\psi(\theta\!+\!u)-\psi(\theta)$, $u\!\ge\! 0$. 
\end{cor}
\begin{pf} Recall from Definition \ref{Levyfor} that $\xi^{h,\psi}(m)=\eta(h)^{m-1}(-1)^m\psi^{(m)}(\eta(h))/\psi^\prime(\eta(h))m!$, $m\neq 1$. We see from the recursive definition of $\widetilde{Q}_{\overline{H};\theta}$ that $\cF^{\rm AD}(\theta)\cap R^h(\cF^{\rm AD}(0))\!\sim\!{\rm GW}(\xi_\theta^{h,\psi},c_\theta^{h,\psi};\mu^{h,\psi})$, where
  $c^{h,\psi}_\theta\!=\!c^{h,\psi}\!+\!\overline{H}_1((0,\theta])=\psi^\prime(\eta(h))\!+\!(\psi^\prime(\eta(h)\!+\!\theta)\!-\!\psi^\prime(\eta(h)))=\psi^\prime(\eta(h)\!+\!\theta)=\psi_\theta^\prime(\eta(h))$ and\vspace{-0.1cm}
  \begin{eqnarray*} \hspace{-0.0cm}g_{\xi^{h,\psi}_\theta}(s)\!&\!\!\!\!=\!\!\!\!&\!\frac{c^{h,\psi}}{c^{h,\psi}\!+\!\overline{H}_1((0,\theta])}\sum_{m\neq 1}\xi^{h,\psi}(m)\Big(H_m((\theta,\infty])s^m+H_m((0,\theta])\Big)\!+\!\left(1\!-\!\frac{c^{h,\psi}}{c^{h,\psi}\!+\!\overline{H}_1((0,\theta])}\right)\\
    &\!\!\!\!=\!\!\!\!&\!\frac{1}{\psi^\prime(\eta(h)\!+\!\theta)}\sum_{m\neq 1}\frac{(\eta(h))^{m-1}(-1)^m}{m!}\left(\psi^{(m)}(\eta(h)\!+\!\theta)s^m+\psi^{(m)}(\eta(h))-\psi^{(m)}(\eta(h)\!+\!\theta)\right)\\[-0.3cm]
	&\!\!\!\!\!\!\!\!&\!\hspace{9.4cm}+\left(1-\frac{\psi^\prime(\eta(h))}{\psi^\prime(\eta(h)+\theta)}\right)\\
    &\!\!\!\!=\!\!\!\!&\!s+\frac{\psi(\eta(h)\!+\!\theta\!-\!\eta(h)s)+\psi(\eta(h)\!-\!\eta(h))-\psi(\eta(h)\!+\!\theta\!-\!\eta(h))}{\eta(h)\psi^\prime(\eta(h)+\theta)}=s+\frac{\widetilde{\psi}_\theta((1\!-\!s)\eta(h))}{\eta(h)\widetilde{\psi}_\theta^\prime(\eta(h))}.
  \end{eqnarray*}
  The desired result follows from Definition \ref{Levyfor} and the fact that 
  $\cF^{\rm AD}(\theta)\cap R^h(\cF^{\rm AD}(0))\rightarrow\cF^{\rm AD}(\theta)$ almost surely, as $h\downarrow 0$. 
\end{pf}\vspace{-0.1cm}

\noindent The analogous result for Galton-Watson trees with pruning at branch points is as follows.
\begin{prop}\label{propRhGWprun} Let $(\cF^{B}(\theta),\theta\ge 0)$ be a pruning process of a ${\rm GW}(\xi;\mu)$-real forest $\cF^{B}(0)=\cF$ with pruning at branch points at exponential rates, and let $h\in\bN$. Then $(\cF^{B}(\theta)\cap R^h(\cF),\theta\ge 0)$ is an $H$-pruning process of $R^h(\cF)\sim{\rm GW}(\xi^h;\mu^h)$, with $\xi^h,\mu^h$ as in {\rm(\ref{xihmuh})} and \vspace{-0.1cm}
$$H_1((\theta,\infty])\!=\!\frac{g^{\prime}_\xi(pe^{-\theta})}{g^\prime_\xi(p)},\ \theta\ge 0,\quad\mbox{and}\  H_m((\theta,\infty])\!=\!\frac{e^{-(m-1)\theta}g_\xi^{(m)}(pe^{-\theta})}{g_\xi^{(m)}(p)},\ \theta\ge 0,\ m\!\ge\! 2,\vspace{-0.1cm}$$
  where $p=Q_\xi(\Gamma\le h)$. In particular, $\cF^B(\theta)\cap R^h(\cF)\sim{\rm GW}(\xi_\theta^h;\mu^h)$, where\vspace{-0.1cm}
  $$g_{\xi_\theta^h}(s)=\frac{g_\xi(pe^{-\theta}+s(1-p)e^{-\theta})-pe^{-\theta}-g_\xi(e^{-\theta})+e^{-\theta}}{(1-p)e^{-\theta}}.$$
\end{prop}
The proof is more elementary than for Abraham-Delmas pruning processes and left to the reader.

\section{Convergence of $H$-pruning processes}\label{secH}\label{nodesandedges}

In this section we establish the following general convergence result for $H$-pruning processes of Galton-Watson trees in the discrete limit regime of Theorem \ref{invprincdiscr}. We will later use this result to establish Theorem \ref{thmprocconv}, and 
we will use similar arguments for Theorem \ref{thmprocconvedge}.

\begin{thm}\label{thmconvH} In the setting of Theorem \ref{invprincdiscr}, consider a sequence $(\cF_n(\theta),\theta\ge 0)$ of $H^{(n)}$-pruning processes of ${\rm GW}(\xi_n;\mu_n)$-real forests, $n\ge 1$, and suppose that $H^{(n)}=(H_m^{(n)}, m\ge 1)$ is such that  
  $$\gamma_nH^{(n)}_1\rightarrow\overline{H}_1\ \mbox{vaguely on $[0,\infty)$,}\ \mbox{and}\  H_m^{(n)}\rightarrow H_m\ \mbox{weakly on $[0,\infty]$, }m\ge 2,\ \mbox{as $n\rightarrow\infty$.}$$
  Suppose furthermore that $H_m^{(n)}$, $m\ge 1$, $n\ge 1$, $\overline{H}_1$ and $H_m$, $m\ge 2$, are non-atomic. 
  Then $(\cF_n(\theta)/\gamma_n,\theta\ge 0)\convd(\cF(\theta),\theta\ge 0)$, where $(\cF(\theta),\theta\ge 0)$ is an $\overline{H}$-pruning process of a ${\rm GW}(\xi,c;\mu)$-real forest with $\overline{H}=(\overline{H}_1;H_m,m\ge 2)$.
\end{thm}

\begin{rem}\rm Theorem \ref{thmconvH}, as well as our proofs, remain valid when $H_m^{(n)}$ and/or $H_m$, $m\ge 2$, have an atom at $\infty$. 
  Finite-dimensional convergence, as well as our proofs, hold for general pruning time distributions -- we just need to exclude
  the countable number of $\theta$-values that appear as atoms of pruning time distributions. Showing tightness in 
  Lemma \ref{lm36} is straightforward in many special cases, e.g. when pruning times are integer-valued, but the general result 
  appears to require a less immediate extra argument to deal with multiple pruning events, which we do not attempt here, as we will not require this higher generality (see 
  \cite[Theorem VI.2.15]{JaS} for convergence criteria when processes are increasing but not simple counting functions). 
\end{rem}


\noindent The proof of Theorem \ref{thmconvH} is spread over the following three subsections. Specifically, we will establish finite-dimensional convergence in Proposition \ref{finiteH} and tightness in $\bD([0,\infty),\bT)$ in Corollary \ref{cor38}.

\subsection{One-dimensional convergence}

We start by a simple lemma, which follows easily from the definition of $H/\overline{H}$-pruning processes.

\begin{lm}\label{lm28}\begin{enumerate}\item[\rm(a)] Let $(\cT^H(\theta),\theta\ge 0)$ be an $H$-pruning process of a ${\rm GW}(\xi)$-real tree as in Definition \ref{prbrpp}. Let $\theta\ge 0$. Then $\cT^H(\theta)$ is a ${\rm GW}(\xi^\theta)$-real tree, where
  $$\xi^\theta(i)=\xi(i)H_i((\theta,\infty]),\ i\ge 1,\qquad \xi^\theta(0)=\xi(0)+\sum_{i\ge 1}\xi(i)H_i((0,\theta]).$$
  \item[\rm(b)] Let $(\cT^{\overline{H}}(\theta),\theta\ge 0)$ be an $\overline{H}$-pruning process of a ${\rm GW}(\xi,c)$-real tree as in Definition \ref{prbrpp2}. Let $\theta\ge 0$. Then $\cT^{\overline{H}}(\theta)$ is a ${\rm GW}(\xi^{\theta,c},c^\theta)$-real tree, where
  $c^\theta=c+\overline{H}_1((0,\theta])$ and
  $$\xi^{\theta,c}(i)=\frac{c\xi(i)H_i((\theta,\infty])}{c+\overline{H}_1((0,\theta])},\ i\ge 2,\quad \xi^{\theta,c}(0)=\xi(0)+\frac{1}{c+\overline{H}_1((0,\theta])}\sum_{i\ge 2}\xi(i)H_i((0,\theta]).$$
  \end{enumerate}
\end{lm}

\begin{prop}\label{oneH}\begin{enumerate}\item[\rm(i)] In the setting of Theorem \ref{thmconvH}, $\cF_n(\theta)/\gamma_n\convd\cF(\theta)$ for each $\theta\ge 0$.
\item[\rm(ii)] In the case $\mu_n=1$, $n\ge 1$, when the forest $\cF_n(\theta)$ consists of a single tree $\cT_n(\theta)$, $n\ge 1$, and hence 
$\cF(\theta)$ consists of a single tree $\cT(\theta)$, we have
  $$(\cT_n(\theta)/\gamma_n,D(\cT_n(\theta)/\gamma_n),\fk(\cT_n(\theta)/\gamma_n))\convd
    (\cT(\theta),D(\cT(\theta)),\fk(\cT(\theta))).$$
\end{enumerate}
\end{prop}
\begin{pf} First note that in forests of $H$-pruning processes, the number of trees does not depend on $\theta$, so we only need to consider the case of single trees, i.e.\ $\mu=\mu_n=\delta_1$. We will apply Theorem \ref{invprincdiscr} to the pruned trees. To this end, note that with convergence of offspring and pruning distributions as assumed in Theorem \ref{thmconvH}, we obtain, as $n\rightarrow\infty$,  $$\gamma_n(1\!-\!\xi_n^\theta(1))=\gamma_n(1\!-\!\xi_n(1)H_1^{(n)}((\theta,\infty]))=\gamma_n(1\!-\!\xi_n(1))+\xi_n(1)\gamma_nH_1^{(n)}((0,\theta])\rightarrow c+\overline{H}_1((0,\theta]),$$
where $\xi_n^\theta$ is associated with $\xi_n$ and $H$ as in Lemma \ref{lm28}. For $i\ge 2$, as $n\rightarrow\infty$,
$$\widetilde{\xi}^\theta_n(i):=\frac{\xi^\theta_n(i)}{1-\xi^\theta_n(1)}=\frac{\gamma_n(1-\xi_n(1))\widetilde{\xi}_n(i)H_i^{(n)}((\theta,\infty])}{\gamma_n(1-\xi^\theta_n(1))}\rightarrow\frac{c\xi(i)H_i((\theta,\infty])}{c+\overline{H}_1((0,\theta])}=\xi^{\theta,c}(i).$$
  This establishes criterion (a) of Theorem \ref{invprincdiscr}, and the equivalence with criterion (d) of Theorem \ref{invprincdiscr} completes this proof.
\end{pf}

\subsection{Finite-dimensional convergence}

We first note Lemmas \ref{lm30} and \ref{lm31}, which follow easily from the definitions of $H$- and $\overline{H}$-pruning processes.

\begin{lm}\label{lm30}\begin{enumerate}\item[\rm(a)] Let $(\cT^H(\theta),\theta\ge 0)$ be an $H$-pruning process of a ${\rm GW}(\xi)$-real tree as in Definition \ref{prbrpp}. Then for all $\theta>0$, we have $(D(\cT^H(0)),D(\cT^H(\theta)))\sim(A,A\wedge B)$, where $A\sim{\rm geom}(1-\xi(1))$ and $B\sim{\rm geom}(H_1(0,\theta])$ are independent.
  \item[\rm(b)] Let $(\cT^{\overline{H}}(\theta),\theta\ge 0)$ be an $\overline{H}$-pruning process of a ${\rm GW}(\xi,c)$-real tree as in Definition \ref{prbrpp2}. Then for all $\theta>0$, we have $(D(\cT^{\overline{H}}(0)),D(\cT^{\overline{H}}(\theta)))\sim(A,A\wedge B)$, where $A\sim{\rm Exp}(c)$ and $B\sim{\rm Exp}(\overline{H}_1(0,\theta])$ are independent.
  \end{enumerate}
\end{lm}

\begin{lm}\label{lm31}\begin{enumerate}\item[\rm(a)] Let $(\cT(\theta),\theta\ge 0)$ be an $H$-pruning process of a ${\rm GW}(\xi)$-real tree as in Definition \ref{prbrpp}, and let $\theta^\prime\ge 0$. Then the post-$\theta^\prime$-process  $(\cT^H(\theta^\prime+\theta),\theta\ge 0)$ is an $H^{\theta^\prime}$-pruning process of a ${\rm GW}(\xi^{\theta^\prime})$-real tree, where $\xi^{\theta^\prime}$ is as in Lemma \ref{lm28} and  $H^{\theta^\prime}_i((0,\theta])=H_i((\theta^\prime,\theta^\prime+\theta])/H_i((\theta^\prime,\infty])$, $\theta\ge 0$, if $H_j((\theta^\prime,\infty])>0$, $H^{\theta^\prime}_i=\delta_\infty$ otherwise, $i\ge 1$. 
  \item[\rm(b)] Let $(\cT^{\overline{H}}(\theta),\theta\ge 0)$ be an $\overline{H}$-pruning process of a ${\rm GW}(\xi,c)$-real tree as in Definition \ref{prbrpp2}, and let $\theta^\prime\ge 0$. Then  $(\cT^{\overline{H}}(\theta^\prime+\theta),\theta\ge 0)$ is an $\overline{H}^{\theta^\prime}$-pruning process of a ${\rm GW}(\xi^{\theta^\prime,c},c^{\theta^\prime})$-real tree, where $(\xi^{\theta^\prime,c},c^{\theta^\prime})$ is as in Lemma \ref{lm28} and  $H^{\theta^\prime}_i((0,\theta])=H_i((\theta^\prime,\theta^\prime+\theta])/H_i((\theta^\prime,\infty])$ if $H_i((\theta^\prime,\infty])>0$, $H^{\theta^\prime}_i=\delta_\infty$ otherwise, $i\ge 2$, and $\overline{H}^{\theta^\prime}_1((0,\theta])=\overline{H}_1((\theta^\prime,\theta^\prime+\theta])$, $\theta\ge 0$.
\end{enumerate}  
\end{lm}

\begin{cor}\label{cor32} In the setting of Proposition \ref{oneH}(ii), we have for all $\theta\ge\theta^\prime\ge 0$
 $$(D(\cT_n(\theta^\prime)/\gamma_n),D(\cT_n(\theta)/\gamma_n),1_{\{D(\cT_n(\theta^\prime))=D(\cT_n(\theta))\}})\convd
  (D(\cT(\theta^\prime)),D(\cT(\theta)),1_{\{D(\cT(\theta^\prime))=D(\cT(\theta))\}}).$$
\end{cor}
\begin{pf} By Lemma \ref{lm31}, it suffices to consider the case $\theta^\prime=0$. By Lemma \ref{lm30}, we may write the LHS as $(A_n,A_n\wedge B_n,1_{\{A_n\le B_n\}})$ and the RHS as $(A,A\wedge B,1_{\{A\le B\}})$, where $\gamma_nA_n\sim{\rm geom}(1-\xi_n(1))$ and $\gamma_nB_n\sim{\rm geom}(H^{(n)}((0,\theta]))$ are independent, and $A\sim{\rm Exp}(c)$ and $B\sim{\rm Exp}(\overline{H}_1((0,\theta]))$ are independent. Now $A_n\convd A$ and $B_n\convd B$ are instances of the usual geometric to exponential convergence as demonstrated in the proof of Proposition \ref{oneH}. Then $(A_n,A_n\wedge B_n,A_n-B_n)\convd(A,A\wedge B,B-A)$ as $(a,b)\mapsto(a,a\wedge b,b-a)$ is continuous. Hence,
\begin{eqnarray*}&&\bP(A_n\le a,A_n\wedge B_n\le b,1_{\{A_n\le B_n\}}=1)=\bP(A_n\le a,A_n\wedge B_n\le b,A_n-B_n\le 0)\\
\longrightarrow&&\bP(A\le a,A\wedge B\le b,A-B\le 0)=\bP(A\le a,A\wedge B\le b,1_{\{A\le B\}}=1).\\[-1.2cm]
\end{eqnarray*} 
\end{pf}

\begin{prop}\label{finiteH} In the setting of Theorem \ref{thmconvH}, we have finite-dimensional convergence, i.e.\  $(\cF_n(\theta_1)/\gamma_n,\ldots,\cF_n(\theta_k)/\gamma_n)\!\convd\!(\cF(\theta_1),\ldots,\cF(\theta_k))$ for $0\!=\!\theta_1\!<\!\theta_2\!<\!\cdots\!<\!\theta_k\!<\!\theta_{k+1}\!=\!\infty$.\vspace{-0.1cm}
\end{prop}
\begin{pf} As in Proposition \ref{oneH}, we may assume $\mu\!=\!\mu_n\!=\!\delta_1$. We fix an increasing sequence $(\theta_j,j\ge 1)$ with $\theta_1=0$ and proceed by induction on $k$. For $k=1$, convergence holds by Proposition \ref{oneH} (or by Theorem \ref{invprincdiscr}). For $k\ge 2$, we simplify notation and denote the scaled trees by $\cS_n^{(j)}:=\cT_n(\theta_j)/\gamma_n$, $1\!\le\! j\!\le\! k$, $n\!\ge\! 1$. The induction hypothesis states that \vspace{-0.1cm}
\begin{equation}(\cS_n^{(1)},\ldots,\cS_n^{(j)})\convd \widetilde{Q}_{\overline{H};\theta_2,\ldots,\theta_j}\quad\mbox{for all $j<k$.}\label{doublestar}\vspace{-0.2cm}
\end{equation}
  By Proposition \ref{oneH}, the sequence in $n\!\ge\! 1$ of
  distributions on $\bT\!\times\![0,\infty)\!\times\!\bN$ of $(\cS_n^{(j)},D(\cS_n^{(j)}),\fk(\cS_n^{(j)}))$ is tight for
  each $j\!=\!1,\ldots,k$. By Corollary \ref{cor32}, also the sequence in $n\!\ge\! 1$ of distributions on $[0,\infty)^2\!\times\!\{0,1\}$ of 
  $(D(\cS_n^{(j)}),D(\cS_n^{(k)}),1_{\{D(\cS_n^{(j)})=D(\cS_n^{(k)})\}})$ is tight for each 
  $j\!=\!1,\ldots,k\!-\!1$. As tightness of marginals implies tightness of joint distribution, the sequence of distributions on $(\bT\!\times\![0,\infty)\!\times\!\bN\!\times\!\{0,1\})^k$ of 
  $((\cS_n^{(j)},D(\cS_n^{(j)}),\fk(\cS_n^{(j)}),1_{\{D(\cS_n^{(j)})=D(\cS_n^{(k)})\}}),1\!\le\! j\!\le\! k)$ is 
  tight. Consider any subsequence $(n(r))_{r\ge 1}$ along which the distributions converge. By Skorohod's representation theorem,
  we may assume that they converge almost surely. Denote the limit by $((\cT_j,D_j,\fk_j,B_j),1\le j\le k)$. We deduce from 
  Proposition \ref{oneH} and Corollary \ref{cor32} that
  $$D_j=D(\cT_j),\quad\fk_j=\fk(\cT_j),\quad B_j=1_{\{D(\cT_j)=D(\cT_k)\}},\quad\mbox{for all $j=1,\ldots,k$, almost surely.}\vspace{-0.1cm}$$
  Now recall the definition of $(\underline{D},\underline{\fk},\underline{\vartheta})$ as functions on $\bT^k$ before Definition \ref{dfprunedge} and note that
  \begin{equation}\hspace{-0.2cm}\begin{array}{rcl}\underline{D}\left(\cS_{n(r)}^{(1)},\ldots,\cS_{n(r)}^{(k)}\right)\!\!\!\!&=&\!\!\!\!\min_{1\le j\le k}D\left(\cS_{n(r)}^{(j)}\right)=D\left(\cS_{n(r)}^{(k)}\right)\longrightarrow D\left(\cT_k\right)=\underline{D}\left(\cT_1,\ldots,\cT_k\right)\\[0.2cm]
  \underline{\vartheta}\left(\cS_{n(r)}^{(1)},\ldots,\cS_{n(r)}^{(k)}\right)\!\!\!\!&=&\!\!\!\!\left({\rm Abv}\left(\cS_{n(r)}^{(j)},D(\cS_{n(r)}^{(k)})\right),1\le j\le k\right)\longrightarrow\underline{\vartheta}\left(\cT_1,\ldots,\cT_k\right)\\[0.2cm]
  \underline{\fk}\left(\cS_{n(r)}^{(1)},\ldots,\cS_{n(r)}^{(k)}\right)\!\!\!\!&=&\!\!\!\!\left(\fk\left(\cS_{n(r)}^{(j)}\right)1_{\{D(\cS_{n(r)}^{(1)})=D(\cS_{n(r)}^{(k)})\}}+1_{\{D(\cS_{n(r)}^{(j)})>D(\cS_{n(r)}^{(k)})\}},1\le j\le k\right)\\[0.3cm]
  &&\hspace{-0.8cm}\longrightarrow\left(\fk\left(\cT_j\right)1_{\{D(\cT_1)=D(\cT_k)\}}+1_{\{D(\cT_j)>D(\cT_k)\}},1\!\le\! j\!\le\! k\right)=\underline{\fk}\left(\cT_1,\ldots,\cT_k\right) 
  \end{array} \label{star}\vspace{-0.1cm}
  \end{equation}
  almost surely, as $r\rightarrow\infty$. Denote by $Q_r$ the distribution on $\bT^k$ of 
  $(\cS_{n(r)}^{(1)},\ldots,\cS_{n(r)}^{(k)})$ and by $Q$ the distribution on $\bT^k$ of $(\cT_1,\ldots,\cT_k)$. For $1\le j\le k$, projections $\pi^{k,j}\colon\bT^k\rightarrow\bT^k$ are induced by $\pi^{k,j}(T_1,\ldots,T_k)=(T_1,\ldots,T_j,\{\rho\},\ldots,\{\rho\})$. These projections are continuous. Let $g\colon[0,\infty)\rightarrow[0,\infty)$ and $G\colon\bT\rightarrow[0,\infty)$ be bounded continuous functions. By Proposition  \ref{dfprunbranch}, $Q_r$ satisfies for all $i_1=\cdots=i_j>i_{j+1}=\cdots=0$, $1\le j\le k$, $$Q_r(g(\underline{D})1_{\{\underline{\fk}=(i_1,\ldots,i_k)\}}G(\underline{\vartheta}))=Q_r(g(\underline{D}))Q_r(\underline{\fk}=(i_1,\ldots,i_k))Q_r^{\circledast i_1}(G\circ\pi^{k,j}).$$
We have $Q_r\rightarrow Q$ weakly, as $r\rightarrow\infty$, and by (\ref{star}), we have shown that the joint distributions of $(\underline{D},\underline{\fk},\underline{\vartheta})$ under $Q_r$ converge to the joint distribution of $(\underline{D},\underline{\fk},\underline{\vartheta})$ under $Q$. Hence for all $g,G,i_1,\ldots,i_k$ as above
\begin{equation}\label{uniqQ}Q(g(\underline{D})1_{\{\underline{\fk}=(i_1,\ldots,i_k)\}}G(\underline{\vartheta}))=Q(g(\underline{D}))Q(\underline{\fk}=(i_1,\ldots,i_k))Q^{\circledast i_1}(G\circ\pi^{k,j}).
\end{equation}
By induction hypothesis (\ref{doublestar}), we can identify the limit of projections onto the first $j<k$ components, and hence $Q^{\circledast i_1}(G\circ\pi^{k,j})=\widetilde{Q}^{\circledast i_1}_{\overline{H};\theta_2,\ldots,\theta_j}(G)$. By Proposition \ref{oneH}(ii), $\underline{D}\sim{\rm Exp}(c+\overline{H}((0,\theta_k]))$ under $Q$, and we also check, using arguments as for Proposition \ref{oneH}, that for all $i_1,\ldots,i_k$ as above,\vspace{-0.4cm}
\begin{eqnarray*}Q(\underline{\fk}=(i_1,\ldots,i_k))&\!\!\!\!=\!\!\!\!&\lim_{r\rightarrow\infty}Q_r(\underline{\fk}=(i_1,\ldots,i_k))=\lim_{r\rightarrow\infty}\frac{\xi_{n(r)}(i_1)H_{i_1}^{(n(r))}((\theta_j,\theta_{j+1}])}{1-\xi_{n(r)}(1)H_1^{(n(r))}((\theta_k,\infty])}\\ &\!\!\!\!=\!\!\!\!&\lim_{r\rightarrow\infty}\frac{\gamma_{n(r)}(1\!-\!\xi_{n(r)}(1))\widetilde{\xi}_{n(r)}(i_1)H_{i_1}^{(n(r))}((\theta_j,\theta_{j+1}])}{\gamma_{n(r)}(1\!-\!\xi_{n(r)}(1)H_1^{(n(r))}((\theta_k,\infty]))}=\frac{c\xi(i_1)H_{i_1}((\theta_j,\theta_{j+1}])}{c\!+\!\overline{H}_1((0,\theta_k])},\vspace{-0.1cm}
\end{eqnarray*}
applying part 3.$^{\!\!B}$ of Proposition \ref{dfprunbranch} and Theorem \ref{invprincdiscr}(a).    
But then $Q=\widetilde{Q}_{\overline{H};\theta_2,\ldots,\theta_k}$ is uniquely identified by (\ref{uniqQ}). Since $Q$ does not depend on $(n(r))_{r\ge 1}$, the tight sequence $Q_r$, $r\ge 1$, of distributions on $\bT^k$ converges to $\widetilde{Q}_{\overline{H};\theta_2,\ldots,\theta_k}$. \pagebreak 
\end{pf}

\subsection{Tightness}

For $T\in\bT$, let $V_a(T)=\#{\rm Br}({\rm Blw}(T,a))+\#{\rm Lf}({\rm Blw}(T,a))$ be the number of branch points and leaves of ${\rm Blw}(T,a)$ including leaves at height $a$ that are due to the truncation at height $a$, but excluding the root. We need the following result only to demonstrate a method of proof that we then refine for the following result, which we do need.

\begin{lm} In the setting of Theorem \ref{invprincdiscr}, the distribution of $V_a$ under $P_{\xi_n}^{\mu_n}(\cdot/\gamma_n)$ converges weakly to the distribution of $V_a$ under $P_{\xi,c}^{\mu}$. 
\end{lm}
\begin{pf} It is straightforward to deduce the result for forests from the corresponding result for the case $\mu=\mu_n=\delta_1$ of single trees. Since $\xi$ is conservative, $V_a$ is $\bN$-valued under $Q_{\xi,c}$. Let $p_a^{(n)}(m)=Q_{\xi_n}(V_a(\cdot/\gamma_n)=m)$ and $p_a(m)=Q_{\xi,c}(V_a=m)$, $m\ge 1$. We will show that $|p_{a_n}^{(n)}(m)-p_{a_n}(m)|\rightarrow 0$ for all $a_n\rightarrow a\ge 0$ and $m\ge 1$, as $n\rightarrow\infty$, using strong induction on $m$. For $m=1$, this holds since $(D,\fk)$ under $Q_{\xi_n}(\cdot/\gamma_n)$ converges and since the distribution of $D$ under $Q_{\xi,c}$ is continuous, so, as $n\rightarrow\infty$,
$$p_a^{(n)}(1)=Q_{\xi_n}(D(\cdot/\gamma_n)\ge a\mbox{ or }\fk=0)\rightarrow Q_{\xi,c}(D\ge a\mbox{ or }\fk=0)=p_a(1).$$
This convergence is uniform, because $a\mapsto p_a^{(n)}(1)$ is decreasing for all $n\ge 1$ and since $a\mapsto p_a(1)$ is continuous. For $m=2$, we have
$p_a^{(n)}(2)=0=p_a(2)$, and for $m\ge 3$, we have
$$p_a^{(n)}(m)=\sum_{\ell=2}^{m-1} Q_{\xi_n}(D(\cdot/\gamma_n)<a,\fk=\ell,V_{a-D(\cdot/\gamma_n)}(\vartheta/\gamma_n)=m-1)$$
and for all $0\le b\le a$
$$Q_{\xi_n}^{\circledast\ell}(V_{a-b}(\cdot/\gamma_n)=m)=\sum_{\genfrac{}{}{0pt}{}{i_1\ge 1,\ldots,i_\ell\ge 1\colon} {i_1+\cdots+i_\ell=m-1}}\prod_{j=1}^\ell p_{a-b}^{(n)}(i_j).$$
We consider the bounded and continuous function $b\mapsto\prod_{j=1}^\ell p_{b}(i_j)=:f(b)$. Since the terms in the products are probabilities in $[0,1]$, we estimate
\begin{eqnarray*}&&\hspace{-0.5cm}\left|Q_{\xi_n}\left(1_{\{D(\cdot/\gamma_n)<a_n,\fk=\ell\}}\prod_{j=1}^\ell p_{a_n-D(\cdot/\gamma_n)}^{(n)}(i_j)\right)-Q_{\xi,c}\left(1_{\{D<a_n,\fk=\ell\}}\prod_{j=1}^\ell p_{a_n-D}(i_j)\right)\right|\\
  &&\le \left|Q_{\xi_n}\left(1_{\{a_n-D(\cdot/\gamma_n)>0,\fk=\ell\}}\left(\prod_{j=1}^\ell p_{a_n-D(\cdot/\gamma_n)}^{(n)}(i_j)-\prod_{j=1}^\ell p_{a_n-D(\cdot/\gamma_n)}(i_j)\right)\right)\right|\\
	&&\hspace{0.5cm}+\left|Q_{\xi_n}\left(1_{\{a_n-D(\cdot/\gamma_n)>0,\fk=\ell\}}\prod_{j=1}^\ell p_{a_n-D(\cdot/\gamma_n)}(i_j)\right)-Q_{\xi,c}\left(1_{\{a-D>0,\fk=\ell\}}\prod_{j=1}^\ell p_{a-D}(i_j)\right)\right|\\
	&&\hspace{0.5cm}+\left|Q_{\xi,c}\left(1_{\{a_n-D>0,\fk=\ell\}}\prod_{j=1}^\ell p_{a_n-D}(i_j)-1_{\{a-D>0,\fk=\ell\}}\prod_{j=1}^\ell p_{a-D}(i_j)\right)\right|\\
	&&\le\sum_{j=1}^\ell Q_{\xi_n}\left(\left|p_{a_n-D(\cdot/\gamma_n)}^{(n)}(i_j)-p_{a_n-D(\cdot/\gamma_n)}(i_j)\right|\right)\\
    &&\hspace{0.5cm}+\left|Q_{\xi_n}\left(1_{\{a_n-D(\cdot/\gamma_n)>0,\fk=\ell\}}f(a_n-D(\cdot/\gamma_n))\right)-Q_{\xi,c}\left(1_{\{a-D>0,\fk=\ell\}}f(a-D)\right)\right|\\
	&&\hspace{0.5cm}+\left|Q_{\xi,c}\left(1_{\{a_n-D>0,\fk=\ell\}}f(a_n-D)-1_{\{a-D>0,\fk=\ell\}}f(a-D)\right)\right|\\
    &&\hspace{-0.1cm}\longrightarrow 0,
\end{eqnarray*}
as $n\rightarrow\infty$ and $a_n\rightarrow a$, by the induction hypothesis, the Dominated Convergence Theorem, and since $(D,\fk)$ under $Q_{\xi_n}(\cdot/\gamma_n)$ converges. Since there are only finitely many $\ell=2,\ldots,m-1$ and $i_1\ge 1,\ldots,i_\ell\ge1$ with $i_1+\cdots+i_\ell=m-1$, we obtain 
$$\left|p_{a_n}^{(n)}(m)-p_{a_n}(m)\right|\rightarrow 0,$$
as $n\rightarrow\infty$ and $a_n\rightarrow a$, and the induction proceeds.\pagebreak[2]
\end{pf}

Denote by $L_a(T)\in[0,\infty]$ the total length of a real tree $T\in\bT$ truncated at height $a$, i.e.\ the total length of ${\rm Blw}(T,a)$. We can refine this result to obtain joint convergence of total lengths below height $a$ with the numbers $V_{a,i}$ of branch points $x$ below height $a$ that have degree $i+1$, $i\ge 2$.

\begin{lm}\label{lm29} In the setting of Theorem \ref{invprincdiscr}, the distribution of $(L_a,(V_{a,i},i\ge 2))$ under $P_{\xi_n}^{\mu_n}(\cdot/\gamma_n)$ converges weakly to the distribution of $(L_a,(V_{a,i},i\ge 2))$ under $P_{\xi,c}^\mu$.
\end{lm}
\begin{pf} We use the same method of proof as in the previous lemma to show that for all $\lambda\ge 0$ and for all ${\bf m}\in {\bf M}:=\{(m_i)_{i\ge 2}\colon m_i\in\bN,\sum_{i\ge 2}m_i<\infty\}$,
  $$p_a^{(n)}({\bf m}):=Q_{\xi_n}(e^{-\lambda L_a(\cdot/\gamma_n)}1_{\{V_{a,i}(\cdot/\gamma_n)=m_i,i\ge 2\}})\rightarrow 
                  Q_{\xi,c}(e^{-\lambda L_a}1_{\{V_{a,i}=m_i,i\ge 2\}})=:p_a({\bf m}),$$
  and indeed $|p_{a_n}^{(n)}({\bf m})-p_{a_n}({\bf m})|\rightarrow 0$, as $n\rightarrow\infty$ and $a_n\rightarrow a\ge 0$. Then
  $$p_a^{(n)}({\bf 0})=Q_{\xi_n}(e^{-\lambda L_a(\cdot/\gamma_n)}1_{\{D(\cdot/\gamma)\ge a\ \mathrm{or}\ \fk=0\}})\rightarrow Q_{\xi,c}(e^{-\lambda L_a}1_{\{D\ge a\ \mathrm{or}\ \fk=0\}})=p_a({\bf 0}),$$
  and as a continuous limit of decreasing functions, this holds uniformly in $a\ge 0$. Writing ${\bf e}_\ell$ for the $\ell$th unit 
  vector $(0,\ldots,0,1,0,\ldots)$ in ${\bf M}$, we note that 
  $$p_a^{(n)}({\bf m})=\sum_{\ell\ge 2\colon m_\ell\ge 1}\sum_{\genfrac{}{}{0pt}{}{({\bf m}^{(1)},\ldots,{\bf m}^{(\ell)})\in {\bf M}^\ell\colon} {{\bf m}^{(1)}+\cdots+{\bf m}^{(\ell)}={\bf m}-{\bf e}_\ell}} Q_{\xi_n}\left(e^{-\lambda D(\cdot/\gamma_n)}1_{\{D(\cdot/\gamma_n)<a,\fk=\ell\}}\prod_{i=1}^\ell p_{a-D(\cdot/\gamma_n)}^{(n)}({\bf m}^{(i)})\right).$$
  The remainder of the proof is now easily adapted from the proof of the previous lemma.
\end{pf}
For a pruning process $\cX\!:=\!(\cF(\theta),\theta\!\ge\! 0)$, we denote by $M_{a,\theta}(\cX)=\#\{\theta^\prime\!\in\!(0,\theta]\colon{\rm Blw}(\cF(\theta^\prime),a)\!\neq\!{\rm Blw}(\cF(\theta^\prime-),a)\}$ the number of pruning events that correspond to jump times of the pruning process during time interval $(0,\theta]$ and below height $a$. Using the point process $\cP$ from Definition \ref{prbrpp},  
we also consider the number that includes further pruning times for events in components when they have already been disconnected from the root; we denote this total number of pruning events in $(0,\theta]$ and below height $a$ by $N_{a,\theta}(\cX)=\cP([0,\theta]\times{\rm Blw}(\cF(0),a))$.

\begin{lm}\label{lm36} In the setting of Theorem \ref{thmconvH}, for a sequence $\cX_n=(\cF_n(\theta),\theta\ge 0)$ of $H^{(n)}$-pruning processes of
  ${\rm GW}(\xi_n,\mu_n)$-real forests, $n\ge 1$, we have, as $n\rightarrow\infty$,
  $$(N_{a,\theta}(\cX_n/\gamma_n),\theta\ge 0)\rightarrow(N_{a,\theta},\theta\ge 0),\qquad\mbox{in distribution in $\bD([0,\infty),\bN)$,}$$
  where the distribution of $(N_{a,\theta},\theta\ge 0)$ is, as follows. For a ${\rm GW}(\xi,c;\mu)$-real forest $\cF$, conditionally given $(L_a(\cF),(V_{a,i}(\cF),i\ge 2))$, the distribution of $(N_{a,\theta},\theta\!\ge\! 0)$ is that of a sum of an inhomogeneous Poisson process with intensity measure $L_a(\cF)\overline{H}_1(d\theta)$ and an independent counting process 
  $$\sum_{i\ge 2}\sum_{j=1}^{V_{a,i}(\cF)}1_{\{\Theta_{i,j}\le\theta\}},\quad\mbox{where $\Theta_{i,j}\sim H_i$, $1\le j\le V_{a,i}(\cF)$, $i\ge 2$, are independent.}$$
\end{lm}
\begin{pf} As $(N_{a,\theta}(\cX_n/\gamma_n),\theta\!\ge\! 0)$ is a simple counting process, we only need finite-dimensional convergence, see \cite[Proposition VI.3.37(b)]{JaS}. We use notation $L_a^{(n)}\!\!=\!L_{[\gamma_na]}(\cF_n(0))$ and
  $V_{a,i}^{(n)}\!\!=\!V_{a,i}(\cF_n(0)/\gamma_n)$, $i\!\ge\! 2$. Let $0\!=\!\theta_1\!<\!\theta_2\!<\!\cdots\!<\!\theta_k\!<\!\theta_{k+1}\!=\!\infty$. We set 
  $q_{i,j}^{(n)}\!=\!H^{(n)}_i((\theta_j,\theta_{j+1}])$, $i\!\ge\! 1$, $1\!\le\! j\!\le\! k$, and denote by $N_a^{(n)}$ the vector 
  with entries $N_{a,\theta_{j+1}}(\cX_n/\gamma_n)-N_{a,\theta_j}(\cX_n/\gamma_n)$, $1\le j\le k-1$. Let $m\ge 1$, 
  $m_i\ge 0$, $i\ge 2$ with $\sum_{i\ge 2}m_i\le m$ and $n_1\ge 0,\ldots,n_{k-1}\ge 0$. We also set $m_1:=m-\sum_{i\ge 2}m_i$ and
  $n_k:=m-n_1-\cdots-n_{k-1}$. Then
  $$\bP(N_a^{(n)}\!=\!(n_1,\ldots,n_{k-1})|L_a^{(n)}\!=\!m,V_{a,i}^{(n)}\!=\!m_i,i\!\ge\! 2)=\sum_{\genfrac{}{}{0pt}{}{r_{i,j}\in\bN,i\ge 2,1\le j\le k\colon}{\genfrac{}{}{0pt}{}{\sum_jr_{i,j}=m_i,i\ge 2;}{\sum_ir_{i,j}=n_j,1\le j\le k}}}m!\prod_{\genfrac{}{}{0pt}{}{i\ge 1\colon}{m_i\ge 1}}\prod_{j=1}^k\frac{(q_{i,j}^{(n)})^{r_{i,j}}}{r_{i,j}!},$$
  where the index set of the sum just captures formally all possible matchings of $m=n_1+\cdots+n_k$ pruning times
  distinguished only according to $k$ time intervals with $m=\sum_{i\ge 1}m_i$ vertices dinstinguished only according to their numbers of
  subtrees. By Lemma \ref{lm29}, \pagebreak[2] $(L_a^{(n)}/\gamma_n;V_{a,i}^{(n)},i\ge 2)$ converges in distribution.  
  As  $\gamma_nq_{1,j}^{(n)}\!=\!\gamma_nH_1^{(n)}((\theta_j,\theta_{j+1}])\rightarrow\overline{H}_1((\theta_j,\theta_{j+1}])\!=:\!\overline{q}_{1,j}$, $1\!\le\! j\!\le\! k\!-\!1$, and $q_{i,j}^{(n)}\rightarrow H_i((\theta_j,\theta_{j+1}])\!=:\!q_{i,j}$, $1\!\le\! j\!\le\! k$, $i\!\ge\! 2$, we find, 
  \begin{eqnarray*}&&\hspace{-0.5cm}\bE\left(e^{-\lambda L_a^{(n)}/\gamma_n}1_{\{V_{a,i}^{(n)}=m_i,i\ge 2\}}1_{\{\widetilde{N}_a^{(n)}=(n_1,\ldots,n_{k-1})\}}\right)\\
    &&\rightarrow\sum_{\genfrac{}{}{0pt}{}{r_{i,j}\in\bN,i\ge 2,1\le j\le k\colon}{\genfrac{}{}{0pt}{}{\sum_jr_{i,j}=m_i,i\ge 2;}{\sum_ir_{i,j}=n_j,1\le j\le k}}}\bE\left(e^{-\lambda L_a}1_{\{V_{a,i}=m_i,i\ge 2\}}\left(\prod_{j=1}^{k-1}\frac{(L_a\overline{q}_{1,j})^{r_{1,j}}}{r_{1,j}!}e^{-L_a\overline{q}_{1,j}}\right)\prod_{\genfrac{}{}{0pt}{}{i\ge 2\colon}{m_i\ge 1}}\prod_{j=1}^{k}\frac{(q_{i,j})^{r_{i,j}}}{r_{i,j}!}\right),
  \end{eqnarray*}
  as $n\rightarrow\infty$, as required for finite-dimensional convergence. 
\end{pf}

\begin{cor}\label{cor37} In the setting of the previous lemma, the sequence of distributions on $\bD([0,\infty),\bN)$ of $(M_{a,\theta}(\cX_n/\gamma_n),\theta\ge 0)$, $n\ge 1$, is tight.
\end{cor}
\begin{pf} The previous lemma yields tightness of the distributions of $(N_{a,\theta}(\cX_n/\gamma_n),\theta\ge 0)$, $n\ge 1$. Since
  $N_{a,\theta}(\cX_n/\gamma_n)-N_{a,\theta^\prime}(\cX_n/\gamma_n)\ge M_{a,\theta}(\cX_n/\gamma_n)-M_{a,\theta^\prime}(\cX_n/\gamma_n)$ for all $\theta\ge\theta^\prime\ge 0$, $n\ge 1$, we obtain the desired result from \cite[Proposition VI.3.35]{JaS}. 
\end{pf}

\begin{cor}\label{cor38} In the setting of Theorem \ref{thmconvH}, the sequence of distributions on $\bD([0,\infty),\bT)$ of $(\cF_n(\theta)/\gamma_n,\theta\ge 0)$, $n\ge 1$, is tight. 
\end{cor}
\begin{pf} Let $\cY_n=(\cF_n(\theta)/\gamma_n,\theta\!\ge\! 0)$. Recall that $(M_{a,\theta}(\cY_n),\theta\!\ge\! 0)$ 
  is the process counting times of pruning events below height $a$ of $\cY_n$, for each $n\ge 1$ and $a>0$. By Proposition \ref{oneH}, and 
  since ${\rm Blw}\colon\bT\times[0,\infty)\rightarrow\bT$ is continuous, the sequence of distributions on $\bT$ of 
  ${\rm Blw}(\cF_n(\theta)/\gamma_n,a)$, $n\ge 1$, is tight.
  Furthermore, the processes
  $({\rm Blw}(\cF_n(\theta)/\gamma_n,a),\theta\ge 0)$, $n\ge 1$, are $\bT$-valued pure jump processes. By Corollary \ref{cor37}, 
  the distributions of their jump counting 
  processes $(M_{a,\theta}(\cY_n),\theta\ge 0)$, $n\ge 1$, form a tight sequence. Finally note that 
  $$\sup_{x\in\bD([0,\infty),\bT)}\;\sup_{\theta\ge 0}\;\delta(x(\theta),{\rm Blw}(x(\theta),a))\le\int_a^\infty re^{-r}dr=(1+a)e^{-a}\rightarrow 0,\qquad\mbox{as $a\rightarrow\infty$.}$$
  It is now straightforward to see and stated in higher generality as an exercise problem in \cite[Problem 3.11.21]{EK86} that 
  this entails the required tightness of distributions on $\bD([0,\infty),\bT)$. 
\end{pf}

Combining Corollary \ref{cor38} with Proposition \ref{finiteH} completes the proof of Theorem \ref{thmconvH}.\hspace*{\fill} $\square$

\section{Proofs of Theorems \ref{thmprocconv} and \ref{thmprocconvedge}}

\subsection{Auxiliary convergence results for the proofs of Theorems \ref{thmprocconv} and \ref{thmprocconvedge}}\label{secaux}

We first collect some results not explicitly stated in \cite{DuWi12}, but that can be proved using similar arguments. The core condition in Lemma \ref{convgf}(i) is a version of the classical condition known to be necessary and sufficient for convergence of associated branching processes, see e.g. \cite{Gri74} or \cite[Chapter 3]{Li2011}. We strengthen this here to uniform convergence of all derivatives, as follows:
\begin{lm}\label{convgf}\label{lm8} In the setting of Theorem \ref{invprinc}, we have the following convergences.
  \begin{enumerate}\item[\rm(i)] $n\gamma_n\left(g_{\xi_n}(e^{-r_n/n})-e^{-r_n/n}\right)\rightarrow\psi(r)$ if $r_n\rightarrow r\ge 0$ for a sequence $r_n\ge 0$, $n\ge 0$.
    \item[\rm(ii)] $\gamma_n\left(1-g_{\xi_n}^\prime(e^{-r_n/n})\right)\rightarrow\psi^\prime(r)$ if $r_n\rightarrow r>q_0$ where
      $q_0$ is the largest root of $\psi$.
    \item[\rm(iii)] $\frac{\gamma_n}{n^{m-1}}g_{\xi_n}^{(m)}(e^{-r_n/n})\rightarrow (-1)^m\psi^{(m)}(r)$ if $r_n\rightarrow r>q_0$, for all $m\ge 2$.
    \item[\rm(iv)] $-n\log q_{\xi_n}\rightarrow q_0$, as $n\rightarrow\infty$, where $q_{\xi_n}=Q_{\xi_n}(\Gamma<\infty)$, $n\ge 1$.
    \item[\rm(v)] $n\eta_n(h_n\gamma_n)\rightarrow \eta(h)>q_0$ if $h_n\rightarrow h>0$.
 \item[\rm(vi)] $g_{\mu_n}(e^{-r_n/n})\rightarrow\int_{[0,\infty)}e^{-rx}\varrho(dx)$ if $r_n\rightarrow r\ge 0$.
  \end{enumerate}
\end{lm}
\begin{pf} (i) Expressing the $\nu_n$-convergence in (\ref{invass}) in terms of Laplace transforms, we find
  $$e^{\psi(r)}\leftarrow\left(e^{r/n}g_{\xi_n}(e^{-r/n})\right)^{\lfloor n\gamma_n\rfloor}
  =\left(1+\frac{\lfloor n\gamma_n\rfloor(e^{r/n}g_{\xi_n}(e^{-r/n})-1)}
  {\lfloor n\gamma_n\rfloor}\right)^{\lfloor n\gamma_n\rfloor},$$
  equivalently $\lfloor n\gamma_n\rfloor(e^{r/n}g_{\xi_n}(e^{-r/n})-1)\rightarrow\psi(r)$ or
  $f_n(r):=n\gamma_n(g_{\xi_n}(e^{-r/n})-e^{-r/n})\rightarrow\psi(r)$. Recall that $\psi$ is convex with
  $\psi(r)\rightarrow\infty$ as $r\rightarrow\infty$. Let $q$ be the unique position where $r\mapsto\psi(r)$ attains its minimum.
  Then for $r_n\rightarrow r\in(0,\infty)\setminus\{q\}$ and $\varepsilon>0$ small enough, $\psi$ is monotonic on
  $[r-3\varepsilon,r+3\varepsilon]$. Let $r>q$ (the case $r<q$, where applicable, is completely analogous). We have pointwise convergence
  $f_n(a)\rightarrow\psi(a)$ at $a=r-2\varepsilon$, $a=r-\varepsilon$ and $a=r$. In particular, since $\psi$ is strictly
  increasing on $[r-3\varepsilon,r+3\varepsilon]$, we will have $f_n(r-2\varepsilon)<f_n(r-\varepsilon)<f_n(r)$ for $n$ large
  enough, and since $f_n$ also has a unique minimum, $f_n$ is increasing on $[r-\varepsilon,r+\varepsilon]$ for $n$ large
  enough. The following basic result (known as Dini's second theorem) implies that the convergence is uniform and so 
  $f_n(r_n)\rightarrow\psi(r)$, as required:
  \begin{enumerate}\item[(R)] Let $f_n\colon[a,b]\rightarrow\bR$ be a sequence of monotonic functions converging pointwise to a
    continuous function $f$. Then the convergence is uniform.
  \end{enumerate}
  The cases $r=0$ and $r=q$ now follow easily via monotonicity and continuity of $\psi$.

  (iv) In the above argument for $r=q_0$, uniform convergence $f_n(a)\rightarrow\psi(a)$ on $[q_0-\varepsilon,q_0+\varepsilon]$ 
  yields that the largest root $r_{\xi_n}$ of $f_n$ lies in $[q_0-\varepsilon,q_0+\varepsilon]$ for $n$ sufficiently large. Hence, 
  $r_{\xi_n}\rightarrow q_0$ as $n\rightarrow\infty$. It is well-known that $q_{\xi_n}$ is the smallest root of $g_{\xi_n}(s)=s$. 
  Hence $q_{\xi_n}=e^{-r_{\xi_n}/n}$, and the result follows.

  (ii) Recall from Theorem \ref{invprinc} that $\cF_n/\gamma_n\convd\cF$. By continuity properties of $R^h$, we also have
  $R^{\gamma_nh}(\cF_n)/\gamma_n\convd R^h(\cF)$ and
  $R^{\lfloor\gamma_n h\rfloor}(\cF_n)/\gamma_n\convd R^h(\cF)$. By Lemma \ref{lmdiscrerasure},
  $R^{\lfloor\gamma_n h\rfloor}(\cF_n)$ is a Galton-Watson tree for each $n\ge 1$, and we can apply Theorem \ref{invprincdiscr} (c)$\Rightarrow$(a). In particular, we obtain the edge
  length parameter of $R^h(\cF)$ as a limit
  $$\psi^\prime(\eta(h))\leftarrow f_n(h):=\gamma_n(1-g_{\xi_n^{\lfloor\gamma_n h\rfloor}}^\prime(0))=\gamma_n(1-g_{\xi_n}^\prime(e^{-\eta_n(\lfloor\gamma_n h\rfloor)})).$$
  Now let $r_n\rightarrow r>q_0$ and $h$ so that $r=\eta(h)$. It is straightforward to check that $f_n$ and $\psi^\prime\circ \eta$ are
  monotonic decreasing, and
  by (R), this convergence is uniform on $[h-2\varepsilon,h+2\varepsilon]$, for any $\varepsilon<h/2$. Since $\psi\circ \eta$ is
  continuous, the range of $f_n$ becomes dense as $n\rightarrow\infty$, and we can find find $h_n^+$ and $h_n^-$ with
  $h_n^+-h_n^-\rightarrow 0$ such that
  $r\leftarrow n\eta_n(\lfloor\gamma_nh_n^+\rfloor)\le r_n\le n\eta_n(\lfloor\gamma_nh_n^-\rfloor)\rightarrow r$, so that
  $$\psi^\prime(r)\leftarrow f_n(h_n^-)\ge\lfloor\gamma_n\rfloor\left(1-g_{\xi_n}^\prime(e^{-r_n/n})\right)\ge f_n(h_n^+)\rightarrow\psi^\prime(r).$$
  (iii) Let $\varepsilon\in(0,r-q_0)$. Then $r_n-\varepsilon>q_0$ for $n$ sufficiently large. By the Mean Value Theorem, there are $u_{m,n}\in(r_n-\varepsilon,r_n)$ and $w_{m,n}\in(r_n,r_n+\varepsilon)$ such that
$$g_{\xi_n}^{(m-1)}(e^{-(r_n-\varepsilon)/n})-g_{\xi_n}^{(m-1)}(e^{-r_n/n})=g_{\xi_n}^{(m)}(e^{-u_{m,n}/n})(e^{\varepsilon/n}-1)e^{-r_n/n}$$
and
$$g_{\xi_n}^{(m-1)}(e^{-r_n/n})-g_{\xi_n}^{(m-1)}(e^{-(r_n+\varepsilon)/n})=g_{\xi_n}^{(m)}(e^{-w_{m,n}/n})(1-e^{-\varepsilon/n})e^{-r_n/n}.$$
Since $g_{\xi_n}^{(m)}\colon[0,1)\rightarrow[0,\infty)$ is increasing, this implies
\begin{eqnarray*}\frac{\frac{\gamma_n}{n^{m-2}}(g_{\xi_n}^{(m-1)}(e^{-r_n/n})-g_{\xi_n}^{(m-1)}(e^{-(r_n+\varepsilon)/n}))}{n(1-e^{-\varepsilon/n})e^{-r_n/n}}
&\le&\frac{\gamma_n}{n^{m-1}}g_{\xi_n}^{(m)}(e^{-r_n/n})\\
&\le&\frac{\frac{\gamma_n}{n^{m-2}}(g_{\xi_n}^{(m-1)}(e^{-(r_n-\varepsilon)/n})-g_{\xi_n}^{(m-1)}(e^{-r_n/n}))}{n(e^{\varepsilon/n}-1)e^{-r_n/n}}.\eeq
Now we proceed by induction on $m$. For $m=2$, the LHS tends to $(\psi^\prime(r+\varepsilon)-\psi^\prime(r))/\varepsilon$, while the RHS tends to $(\psi^\prime(r)-\psi^\prime(r-\varepsilon))/\varepsilon$, by (ii). Therefore, liminf and limsup are bounded by these quantities for all $\varepsilon$, and hence by their limit as $\varepsilon\downarrow 0$, which is $\psi^{\prime\prime}(r)$. Given 
the convergence for $m-1$, the same argument establishes the induction step for $m$, using the induction hypothesis instead of (ii). The factor $(-1)^m$ arises, because the limits of LHS (and similarly RHS) are now $(-1)^{m-1}(\psi^{(m-1)}(r)-\psi^{(m-1)}(r+\varepsilon))/\varepsilon$ and tend to $(-1)^m\psi^{(m)}(r)$.

   (v)-(vi) The results for $\mu_n$ and $\eta_n$ are easier and left to the reader. Let us just point out that $\eta(h)>q_0$ holds since $\int_{q_0}^\infty du/\psi(u)=\infty$ and $\int_{\eta(h)}^\infty du/\psi(u)=h$ for all $h\in(0,\infty)$.
\end{pf}

\subsection{Pruning at branch points and the proof of Theorem \ref{thmprocconv}}\label{secthm2}

Consider the setting of Theorem \ref{thmprocconv}. By Corollary \ref{cor25} it suffices to show convergence 
$$\cX_n^{h_n}:=\left(\cF_n^B(\theta/n)/\gamma_n\cap R^{h_n}(\cF_n^B(0)/\gamma_n),\theta\ge 0\right)\convd\left(\cF^{\rm AD}(\theta)\cap R^h(\cF^{\rm AD}(0)),\theta\ge 0\right)=:\cX^h,$$
for each $h>0$, where we choose $h_n:=\lfloor h\gamma_n\rfloor/\gamma_n$ and note that $h_n\gamma_n\in\bN$ and $h_n\rightarrow h$ as $n\rightarrow\infty$.
 
By hypothesis, $\cF_n^B(0)/\gamma_n\!\convd\!\cF^{\rm AD}(0)$. By Proposition \ref{prop4}(ii), we deduce  
$R^{h_n}(\cF_n^B(0)/\gamma_n)\!\convd\! R^h(\cF^{\rm AD}(0))$. By Lemma \ref{lmdiscrerasure}, $R^{h_n\gamma_n}(\cF_n^B(0))$ is a ${\rm GW}(\xi_n^{\gamma_nh_n};\mu_n^{\gamma_nh_n})$-real forest and by Definition \ref{Levyfor}, $R^h(\cF^{\rm AD}(0))$ is a ${\rm GW}(\xi^{h,\psi},c^{h,\psi};\mu^{h,\psi})$-real forest. Since $R^{h_n}(\cF_n^B(0)/\gamma_n)\convd R^h(\cF^{\rm AD}(0))$, we are in the framework of Theorem 
\ref{invprincdiscr}(c). By Proposition \ref{prop21}, $\cX^h$ is an $\overline{H}$-pruning process, with $\overline{H}\!=(\overline{H}_1;H_m,m\!\ge\! 2)$ as in the proposition. By Proposition \ref{propRhGWprun}, $\cX_n^{h_n}$ is an $H^{(n)}$-pruning process with
\begin{eqnarray*}H_m^{(n)}((\theta,\infty])=e^{-(m-1)\theta/n}\frac{g_{\xi_n}^{(m)}(e^{-\eta_n(\gamma_nh_n)}e^{-\theta/n})}{g_{\xi_n}^{(m)}(e^{-\eta_n(\gamma_nh_n)})}&\longrightarrow& \frac{\psi^{(m)}(\eta(h)+\theta)}{\psi^{(m)}(\eta(h))}=H_m((\theta,\infty])\\
  \gamma_nH_1^{(n)}([0,\theta])=\gamma_n\left(1-\frac{g_{\xi_n}^\prime(e^{-\eta_n(\gamma_nh_n)}e^{-\theta/n})}{g_{\xi_n}^\prime(e^{-\eta_n(\gamma_nh_n)})}\right)&\longrightarrow&\psi^\prime(\eta(h)+\theta)-\psi^\prime(\eta(h))=\overline{H}_1([0,\theta])
\end{eqnarray*}
by Lemma \ref{lm8}, as $n\rightarrow\infty$, since $g_{\xi_n}^\prime(e^{-\eta_{n}(\gamma_nh_n)})\rightarrow 1$, also by Lemma \ref{lm8}. Hence, the assumptions of Theorem \ref{thmconvH} are satisfied and the limit identified, so $\cX_n^{h_n}\convd\cX^h$ in $\bD([0,\infty),\bT)$, as required.\vspace{-0.3cm}
\begin{trivlist}\item[]\hspace*{\fill} $\square$\end{trivlist}

\subsection{Pruning at edges and the proof of Theorem \ref{thmprocconvedge}}\label{secthm3}

The structure of the proof is the same as for Theorem \ref{thmprocconv}, with $\cF^B_n(\theta/n)$ replaced by $\cF^E_n(\theta/\gamma_n)$, $n\ge 1$, and $\cF^{\rm AD}(\theta)$ by $\cF^{\rm AP}(\theta)$, $\theta\ge 0$. But here, $\cX^h$ is an Aldous-Pitman pruning process of a ${\rm GW}(\xi^{h,\psi},c^{h,\psi};\mu^{h,\psi})$-real forest, by Proposition \ref{edgeerase}, and $\cX^{h_n}_n$ is a ${\rm GW}(\xi_n^{h_n\gamma_n};\mu_n^{h_n\gamma_n})$-pruning process with pruning at edges, by Proposition \ref{edgeerasediscr}, with trees and pruning times scaled by $\gamma_n$. 

\begin{rem}\rm $\cX_n^{h_n}$ is almost an $H^{(n)}$-pruning process with $H^{(n)}_m\!=\!{\rm Exp}(\gamma_n)$, $m\!\ge\! 1$, but not
  quite, as $H^{(n)}$-pruning effectively prunes at the top of edges below branch points, while pruning at edges prunes 
  at the bottom of edges. A possible approach would be to couple $\cX_n^{h_n}$ to an $H^{(n)}$-pruning process and to show that
  the differences vanish in the limit. We do not pursue the coupling here, but we establish an invariance principle for pruning at branch points with ${\rm Exp}(1)$-pruning times in Theorem \ref{thmlast} below. We now provide relevant details of a direct approach to $\cX_n^{h_n}\convd\cX^h$. 
\end{rem}

\subsubsection{One-dimensional convergence for $\cX^{h_n}_n\convd\cX^h$}

It is straightforward to see that $\cF^E_n(\theta)\cap R^{h_n\gamma_n}(\cF^E_n(0))$ is a ${\rm GW}(\widehat{\xi}_n^{\theta,h};\widehat{\mu}_n^{\theta,h})$-real forest, where
$$g_{\widehat{\xi}_n^{\theta,h}}(s)\!=\!\frac{g_{\xi_n}(1\!-\!(1\!-\!p_n)e^{-\theta}\!+\!s(1\!-\!p_n)e^{-\theta})-p_n}{1-p_n}\ \mbox{and}\  g_{\widehat{\mu}_n^{\theta,h}}(s)\!=\!g_{\mu_n}(1\!-\!(1\!-\!p_n)e^{-\theta}\!+\!s(1\!-\!p_n)e^{-\theta})$$
with $p_n\!=\!e^{-\eta_n(h_n\gamma_n)}\!=\!Q_{\xi_n}(\Gamma\!\le\! h_n\gamma_n)$. We get one-dimensional convergence from Theorem \ref{invprincdiscr}, as
\begin{eqnarray*}\gamma_n\left(1-\widehat{\xi}_n^{\theta/\gamma_n,h}(1)\right)&\!\!\longrightarrow\!\!&c^{h,\psi}+\theta=\widehat{\psi}_\theta^\prime(\eta(h))=c_\theta^{h,\psi}\\
g_{\widetilde{\widehat{\xi}}_n^{\theta/\gamma_n,h}}(s)&\!\!\longrightarrow\!\!& s+\frac{\psi((1-s)\eta(h))+\theta(1-s)\eta(h)}{\eta(h)(\psi^\prime(\eta(h))+\theta)}=s+\frac{\widehat{\psi}_\theta((1-s)\eta(h))}{\eta(h)\widehat{\psi}^\prime_\theta(\eta(h))}=g_{\xi_\theta^{h,\psi}}(s)\\
g_{\widehat{\mu}_n^{\theta/\gamma_n,h}}(s)&\!\!\longrightarrow\!\!&\int_{[0,\infty)}e^{-x(1-s)\eta(h)}\varrho(dx)=g_{\mu^{h,\psi}}(s), 
\end{eqnarray*}
by straightforward calculations based on Lemma \ref{lm8}. By Proposition \ref{edgeerase} this is the required limit. 

\subsubsection{Finite-dimensional convergence for $\cX^{h_n}_n\convd\cX^h$}
\label{edgefd}

For finite-dimensional convergence, we adapt the proof of Proposition \ref{finiteH}, for scaled trees $\cS_n^{(j)}:=\cT_n(\theta_j/\gamma_n)/\gamma_n$, $1\le j\le k$, $n\ge 1$. The main difference is that $\cS_n^{(k)}=\{\rho\}$ with probability $e^{-\theta_k/\gamma_n}\rightarrow 0$, while conditionally given $\cS_n^{(k)}\neq\{\rho\}$, we have $(D(\cS_n^{(1)}),D(\cS_n^{(k)}))\sim(A_n,A_n\wedge B_n)$, where $\gamma_nA_n\sim{\rm geom}(1-\xi^{h_n\gamma_n}_n(1))$ and 
$\gamma_nB_n\sim{\rm geom}(1-e^{-\theta_k/\gamma_n})$ are independent, as in Lemma \ref{lm30}. This still allows us to show that the limiting distribution $Q$ along any convergent subsequence $(n(r))_{r\ge 1}$ satisfies for all $i_1=\cdots=i_j>i_{j+1}=\cdots=i_k=0$, $0\le j\le k$,
$$Q(g(\underline{D})1_{\{\underline{\fk}=(i_1,\ldots,i_k)\}}G(\underline{\vartheta}))=Q(g(\underline{D}))Q(\underline{\fk}=(i_1,\ldots,i_k))Q^{\circledast i_1}(G\circ\pi^{k,j}),$$
where, inductively, $Q^{\circledast i_1}(G\circ\pi^{k,j})=\widehat{Q}_{\theta_2,\ldots,\theta_j}$, $1\le j\le k-1$, while $\underline{D}\sim{\rm Exp}(c^{h,\psi}+\theta_k)$ under $Q$,
$$Q(\underline{\fk}=(i_1,\ldots,i_k))=\lim_{r\rightarrow\infty}\frac{\gamma_{n(r)}\xi_{n(r)}^{h_{n(r)}\gamma_{n(r)}}(1)\left(e^{-\theta_j/\gamma_{n(r)}}-e^{-\theta_{j+1}/\gamma_{n(r)}}\right)}{\gamma_{n(r)}\left(1-\xi_{n(r)}^{h_{n(r)}\gamma_{n(r)}}(1)e^{-\theta_k/\gamma_{n(r)}}\right)}=\frac{\theta_{j+1}-\theta_j}{c^{h,\psi}+\theta_k}$$
if $i_1=\cdots=i_j=1>i_{j+1}=\cdots=i_k=0$, $1\le j\le k-1$, and
\begin{eqnarray*}&&\hspace{-0.8cm}Q(\underline{\fk}=(i_1,\ldots,i_k))\\
  &&\hspace{-0.5cm}=\lim_{r\rightarrow\infty}\frac{\gamma_{n(r)}\left(1\!-\!\xi_{n(r)}^{h_{n(r)}\gamma_{n(r)}}(1)\right)\widetilde{\xi}_{n(r)}^{h_{n(r)}\gamma_{n(r)}}(i_1)}{\gamma_{n(r)}\left(1\!-\!\xi_{n(r)}^{h_{n(r)}\gamma_{n(r)}}(1)e^{-\theta_k/\gamma_{n(r)}}\right)}\frac{i_1!e^{-i_k\theta_k/\gamma_{n(r)}}}{i_k!}\prod_{j=1}^{k-1}\frac{\left(e^{-\theta_j/\gamma_{n(r)}}\!-\!e^{-\theta_{j+1}/\gamma_{n(r)}}\right)^{i_j-i_{j+1}}}{(i_j\!-\!i_{j+1})!}\\
   &&\hspace{-0.5cm}=\left\{\begin{array}{ll}\displaystyle\frac{c^{h,\psi}\xi^{h,\psi}(i_1)}{c^{h,\psi}+\theta_k}&\quad\mbox{if $1\neq i_1=\cdots=i_k$}\\ 0&\quad\mbox{if $1\neq i_1\ge\cdots\ge i_j>i_{j+1}\ge\cdots\ge i_k\ge 0$,}\end{array}\right.
 \end{eqnarray*}
since $e^{-\theta_k/\gamma_{n(r)}}\rightarrow 1$ and $(e^{-\theta_j/\gamma_{n(r)}}-e^{-\theta_{j+1}/\gamma_{n(r)}})\rightarrow 0$.
By Lemma \ref{lm1} and Proposition \ref{dfprunedgec}, $Q=\widehat{Q}_{\theta_2,\ldots,\theta_k}$, as required. The extension from single trees to forests is straightforward.

\subsubsection{Tightness of the family of distributions of $\cX^{h_n}_n$, $n\ge 1$}

Recall that we denote by $L_a(T)\in[0,\infty]$ the total length of a real tree $T\in\bT$ truncated at height $a$, i.e.\ the total length of ${\rm Blw}(T,a)$.

\begin{lm} In the setting of Theorem \ref{invprincdiscr}, the distribution of $(L_a/\gamma_n,a\ge 0)$ under $P_{\xi_n}^{\mu_n}(\cdot/\gamma_n)$ converges weakly in $\bD([0,\infty),[0,\infty))$ to the distribution of $(L_a,a\ge 0)$ under $P_{\xi,c}^\mu$, as $n\rightarrow\infty$.
\end{lm}
\begin{pf} For fixed $a\ge 0$, this is part of Lemma \ref{lm29}. Here, we give an independent proof of Skorohod convergence, as follows. For $T\in\bT$, let 
  $Z_a(T)=\#\{v\in T\colon d(\rho,v)=a\}\in[0,\infty]$, $a\ge 0$. The convergence of branching
  processes $(Z_{a+}(\cdot/\gamma_n),a\ge 0)$ under $P_{\xi_n}^{\mu_n}$ to $(Z_{a+},a\ge 0)$ under $P_{\xi,c}^\mu$ was obtained in 
  \cite[Theorem 3.24]{DuWi12}. Note that $L_a(T)=\int_0^a Z_{t+}(T)dt$.  Integration is a continuous function from 
  $\bD([0,\infty),[0,\infty))$ to $\bD([0,\infty),[0,\infty))$, see e.g. \cite[Problem 3.11.26]{EK86} or, in much higher 
  generality, \cite[Theorem VI.6.22]{JaS}, so the result follows.   
\end{pf}

Discrete pruning at edges is carried out for each edge (of unit length) at an independent identically distributed time. Recall that for a pruning process $\cX=(\cT(\theta),\theta\ge 0)$, we denote by $M_{a,\theta}(\cX)$ the number of pruning times that are jump times of the pruning process during time interval $(0,\theta]$ and below height $a$. Recall also that in the richer model that includes pruning times for all edges, we denote by $N_{a,\theta}(\cX)$ the total number of pruning times including for those edges already disconnected from the root. Although we will only require exponential pruning time distributions when pruning at edges, we can just as well consider more general pruning time distributions here. 

\begin{lm}\label{lm25} In the setting of Theorem \ref{invprincdiscr}, consider a sequence of non-atomic pruning time distributions $H_1^{(n)}\!$, $n\ge 1$, on 
  $(0,\infty)$, such that $\gamma_nH_1^{(n)}\!\rightarrow\!\overline{H}_1$ vaguely on $[0,\infty)$, as $n\rightarrow\infty$. Let 
  $\cX_n\!=\!(\cF_n^E(\theta),\theta\!\ge\! 0)$
  be a ${\rm GW}(\xi_n;\mu_n)$-pruning process with pruning at edges at independent identically $H_1^{(n)}$-distributed pruning times. 
  Then for each $a\ge 0$, as $n\rightarrow\infty$,
  $$(N_{a,\theta}(\cX_n/\gamma_n),\theta\ge 0)\rightarrow (N_{a,\theta},\theta\ge 0)\qquad\mbox{in distribution in $\bD([0,\infty),\bN)$,}$$
  where the distribution of $(N_{a,\theta},\theta\ge 0)$ is, as follows. For a ${\rm GW}(\xi,c;\mu)$-real forest $\cF$, conditionally given $L_a(\cF)$, the counting process $(N_{a,\theta}(\cF),\theta\ge 0)$ is an inhomogeneous Poisson
  process with intensity measure $L_a(\cF)\overline{H}_1(d\theta)$. 
\end{lm} 
\begin{pf} We adapt the proof of Lemma \ref{lm36}. Let $0=\theta_1<\theta_2<\cdots<\theta_k<\theta_{k+1}=\infty$ and set $q_j^{(n)}=H_1^{(n)}((\theta_j,\theta_{j+1}])$, $1\le j\le k$. We further simplify notation and set $L_a^{(n)}=L_a(\cF_n^E(0))$ and
$L_a=L_a(\cF^{\rm AP}(0))$, also $N_a^{(n)}$ for the vector of increments $N_{a,\theta_{j+1}}(\cX_n/\gamma_n)-N_{a,\theta_j}(\cX_n/\gamma_n)$, $1\le j\le k-1$. Then for all $m\ge 1$ and $(m_1,\ldots,m_{k-1})$ with $m_1+\cdots+m_{k-1}\le m$,
$$\bP(N_a^{(n)}=(m_1,\ldots,m_{k-1})|L_a^{(n)}=m)=m!\left(\prod_{j=1}^{k-1}\frac{(q_j^{(n)})^{m_j}}{m_j!}\right)\frac{(q_k^{(n)})^{m-m_1-\cdots-m_{k-1}}}{(m-m_1-\cdots-m_{k-1})!}.$$
Using $L_a^{(n)}/\gamma_n\rightarrow L_a$ in distribution, $\gamma_nq_j^{(n)}\rightarrow\overline{H}_1((\theta_j,\theta_{j+1}])=:\overline{q}_j$, $1\le j\le k-1$, and $(q_k^{(n)})^{\gamma_n}=(1-(\gamma_nq_1^{(n)}+\cdots+\gamma_nq_{k-1}^{(n)})/\gamma_n)^{\gamma_n}\rightarrow e^{-\overline{q}_1-\cdots-\overline{q}_{k-1}}$, we obtain
$$\bE\left(e^{-\lambda L_a^{(n)}/\gamma_n}1_{\{N_a^{(n)}=(m_1,\ldots,m_{k-1})\}}\right)\rightarrow\bE\left(e^{-\lambda L_a}\prod_{j=1}^{k-1}\frac{(L_a\overline{q}_j)^{m_j}}{m_j!}e^{-L_a\overline{q}_j}\right).\vspace{-0.5cm}$$
\end{pf}

\begin{cor}\label{counttight} In the setting of the previous lemma, the sequence of distributions on $\bD([0,\infty),\bN)$ of 
  $(M_{a,\theta}(\cX_n/\gamma_n),\theta\ge 0)$, $n\ge 1$, is tight.
\end{cor}
\begin{pf} The domination argument of Corollary \ref{cor37} applies again here.
\end{pf}

Tightness of the distributions of $\cX_n^{h_n\gamma_n}$, $n\ge 1$, on $\bD([0,\infty),\bT)$ follows as in Corollary \ref{cor38}. Together with finite-dimensional convergence established in Section \ref{edgefd}, this completes the proof of Theorem \ref{thmprocconvedge}.\hspace*{\fill} $\square$

\subsection{Invariance principle for equal-rate pruning at branch points}\label{secthmlast}

\begin{thm}\label{thmlast} In the setting of Theorem \ref{invprinc}, the associated $H$-pruning processes $(\cF_n^H(\theta),\theta\ge 0)$ with $H_m={\rm Exp}(1)$, $m\ge 1$, converge: 
$$(\cF_n^H(\theta/\gamma_n)/\gamma_n,\theta\ge 0)\convd(\cF^{\rm AP}(\theta),\theta\ge 0)\qquad \text{in }\bD([0,\infty),\bT),$$
where the limit is the Aldous-Pitman pruning process associated with a $(\psi;\varrho)$-L\'evy forest $\cF$.
\end{thm}
\begin{pf} We proceed as for the proof of Theorem \ref{thmprocconv}. Here, $\cX^h$ is an Aldous-Pitman pruning process of a ${\rm GW}(\xi^{h,\psi},c^{h,\psi};\mu^{h,\psi})$-real forest, by Proposition \ref{edgeerase}, and it is straightforward to see that $\gamma_n\cX^{h_n}_n$ is an $H^{(n)}$-pruning process of a ${\rm GW}(\xi_n^{h_n\gamma_n};\mu_n^{h_n\gamma_n})$-real forest, where $H^{(n)}_m={\rm Exp}(\gamma_n)$. Since $R^{h_n}(\cF_n^H(0)/\gamma_n)\convd R^h(\cF^{\rm AD}(0))$ is the same as in the proof of Theorem \ref{thmprocconv}, we are in the framework of Theorem \ref{invprincdiscr}(c) again. To apply Theorem \ref{thmconvH}, we check that for all $\theta\ge 0$
$$\gamma_nH_1^{(n)}((0,\theta])=\gamma_n(1-e^{-\gamma_n\theta})\rightarrow\theta\quad\mbox{and}\quad H_m^{(n)}((0,\theta])=H_1^{(n)}((0,\theta])\rightarrow 0.\vspace{-0.6cm}$$  
\end{pf}

\section{Applications}\label{ascsec}

Kesten \cite{Kes86} studied Galton-Watson trees conditioned on non-extinction. He showed that the resulting tree can be constructed by grafting onto an infinite half-line of vertices forests of Galton-Watson trees. We use this representation to define associated $\bT$-valued trees, which we call \em Kesten trees\em. See also the earlier Kallenberg \cite{Kal77}, where closely related structures are introduced as a tool to study cluster fields, Duquesne \cite{Duq-08} for an invariance principle for Kesten trees, also Athreya et al. \cite[Example 7.7]{ALW} for an application of the Brownian special case to walks on trees. See \cite{AP98,ADH12,AD12} for studies of pruning processes and ascension times separately in discrete and continuum settings. 

\subsection{Pruning of Kesten trees/forests and invariance principles}\label{ascsec1}

Let $\xi$ be an offspring distribution with $\xi(1)<1$ that is critical, i.e.\ $\sum_{i\ge 0}i\xi(i)=1$. Consider 
a random real tree $\cT$ obtained by grafting onto the infinite half-line $[0,\infty)$ at each $m\in\bN\setminus\{0\}$ an independent ${\rm GW}(\xi;\mu)$-real forest $\cF^{(m)}$, where $\mu(i)=(i+1)\xi(i+1)$, $i\ge 0$. We refer (to any random real tree isometric) to $\cT$ as a \em Kesten tree with offspring distribution $\xi$\em. We denote the distribution of the isometry class of $\cT$ in $\bT$ by $K_\xi$.

Let $\xi$ be an offspring distribution with $\xi(1)=0$ that is critical, and let $c\in(0,\infty)$. Let $S_m$, $m\ge 1$, be the times of a Poisson process of rate $c$. Consider a random real tree $\cT$ obtained by grafting onto the infinite half-line $[0,\infty)$ at each $S_m$, $m\ge 1$, an independent ${\rm GW}(\xi,c;\mu)$-real forest $\cF^{(m)}$, where $\mu(i)=(i+1)\xi(i+1)$, $i\ge 0$. We refer (to any random real tree isometric) to $\cT$ as a \em Kesten tree with offspring distribution $\xi$ and lifetime parameter $c$\em. We denote the distribution of the isometry class of $\cT$ in $\bT$ by $K_{\xi,c}$.  

\begin{prop}\begin{enumerate}\item[\rm(i)] Let $\xi$ be critical, $\xi(1)<1$. Then $K_\xi$ is the unique distribution $K$ on $\bT$  that satisfies \vspace{-0.3cm}
$$K(g(D)1_{\{\fk=i+1\}}G(\vartheta))=\sum_{m=1}^\infty g(m)\xi(1)^{m-1}(i+1)\xi(i+1)\left(K\circledast Q_\xi^{\circledast i}\right)(G)\vspace{-0.1cm}$$
for all $i\in\bN$ all nonnegative measurable functions $g$ on $[0,\infty)$ and $G$ on $\bT$.
\item[\rm(ii)] Let $c\in(0,\infty)$, $\xi$ critical and $\xi(1)=0$. Then $K_{\xi,c}$ is the unique distribution $K$ on $\bT$ with\vspace{-0.1cm}
$$K(g(D)1_{\{\fk=i+1\}}G(\vartheta))=\int_{0}^\infty g(x)ce^{-cx}dx(i+1)\xi(i+1)\left(K\circledast Q_{\xi,c}^{\circledast i}\right)(G).$$
\end{enumerate}
\end{prop}
Duquesne \cite{Duq-08} established invariance principles for (sub)critical Galton-Watson trees with immigration, which includes  Kesten trees. In the case that arises for Kesten trees, we can define here $\bT$-valued representations of Duquesne's limiting immigration L\'evy trees, as follows. Let $\psi$ be a branching mechanism of the form (\ref{brmech}) that is critical, i.e.\ $\psi^\prime(0)=0$. Let $\cP=\sum_{v\in I}\delta_{(v,T_v)}$ be a Poisson random measure on $[0,\infty)\times\bT$ with intensity measure $$\ell\times\left(2\beta\bN_\psi+\int_{(0,\infty)}xP_\psi^{\delta_x}\pi(dx)\right),$$
where $\ell$ is Lebesgue measure on $[0,\infty)$ and $\pi$ the L\'evy measure in (\ref{brmech}). 
Consider a random real tree $\cT$ obtained by grafting onto $[0,\infty)$ at each $v\in I$ (a representative of) the forest $T_v$. We refer (to any random real tree isometric) to $\cT$ as a \em $\psi$-Kesten-L\'evy tree\em. We denote the distribution of the isometry class of $\cT$ in $\bT$ by $K_\psi$. The following theorem is \cite[Theorem 1.5]{Duq-08}, restricted to the special case of Kesten trees and pushed forward from coding height functions to $\bT$. 

\begin{thm}[Duquesne \cite{Duq-08}]\label{invprinc2} Using notation from Theorem \ref{invprinc}, suppose that there is a positive sequence $\gamma_n\rightarrow\infty$ such that
\eq\nu_n(\ts n\,\cdot\,)^{*\lfloor n\gamma_n\rfloor}\rightarrow\nu\quad\mbox{weakly,}\qquad\mbox{and}\quad n\eta_n(\lfloor\ts\gamma_n\,\cdot\,\rfloor)\rightarrow \eta\quad\mbox{pointwise,}\label{invass2}
\en
where $\nu$ is such that $\int_{[0,\infty)}e^{-rx}\nu(dx)=e^{\psi(r)}$ for a (sub)critical branching mechanism {\rm(\ref{brmech})} satisfying {\rm(\ref{grey})}. 
Let $\cT_n$ be a Kesten tree with offspring distribution $\xi_n$, $n\ge 1$. Then, as $n\rightarrow\infty$,
$\cT_n/\gamma_n\convd\cT$ in $\bT$, for a $\psi$-Kesten-L\'evy tree $\cT$.
\end{thm}
We now take this result as the starting point to derive results analogous to Theorems \ref{thmprocconv} and \ref{thmprocconvedge}.

\begin{thm}\label{kesconv} In the setting of Theorem \ref{invprinc2}, the associated pruning processes  $(\cT_n^B(\theta),\theta\ge 0)$ with pruning at branch points converge, as $n\rightarrow\infty$: $$(\cT_n^B(\theta/n)/\gamma_n,\theta\ge 0)\convd(\cT^{\rm AD}(\theta),\theta\ge 0)\qquad \text{in }\bD([0,\infty),\bT),$$
where the limit is the Abraham-Delmas pruning process associated with a $\psi$-Kesten-L\'evy tree $\cT$. 
\end{thm}
\begin{thm}\label{kesconvedge} In the setting of Theorem \ref{invprinc2}, the associated pruning processes $(\cT_n^E(\theta),\theta\ge 0)$ with pruning at edges converge, as $n\rightarrow\infty$: 
$$(\cT_n^E(\theta/\gamma_n)/\gamma_n,\theta\ge 0)\convd(\cT^{\rm AP}(\theta),\theta\ge 0)\qquad \text{in }\bD([0,\infty),\bT),$$
where the limit is the Aldous-Pitman pruning process associated with a $\psi$-Kesten-L\'evy tree $\cT$.
\end{thm}
\begin{rem}\label{kestenforests}\rm From Theorems \ref{thmprocconv}, \ref{thmprocconvedge}, \ref{kesconv} and \ref{kesconvedge}, we can deduce generalisations to Kesten forests, which we define as concatenations of a single Kesten tree (or Kesten-L\'evy tree) and an independent forest of Galton-Watson trees (or L\'evy trees). Specifically, continuity of concatenation $\circledast\colon\bT^2\rightarrow\bT$ yields finite-dimensional convergence, and tightness reduces to the analogue of Lemmas
\ref{lm36} and \ref{lm25}, for which we add two independent convergent sequences of counting processes. 
\end{rem}
We leave the details of the proofs of Theorems \ref{kesconv} and \ref{kesconvedge} to the reader. Briefly, it suffices to prove that $h$-erasures converge. In representations of forests grafted onto the infinite half-line $[0,\infty)$, we can then apply Lemma \ref{convgf} for the convergence of the point process of numbers of trees and Theorems \ref{thmprocconv} and \ref{thmprocconvedge} to the convergence of the grafted forests themselves. This is straightforward, because Kesten trees behave nicely under $h$-erasure:

\begin{lm}\begin{enumerate}\item[\rm(i)] Let $h\!\in\!\bN$. Then $R_h$ under $K_\xi$ has distribution $K_{\xi^h}$, with $\xi^h$ from Lemma \ref{lmdiscrerasure}.
  \item[\rm(ii)] Let $h\!>\!0$. Then $R_h$ under $K_{\xi,c}$ has distribution $K_{\xi^{h,c},c^h}$, with $(\xi^{h,c},c^h)$ from Lemma \ref{lmdiscrerasure}.
  \item[\rm(iii)] Let $h\!>\!0$. Then $R_h$ under $K_\psi$ has distribution $K_{\xi^h,c^h}$, with $(\xi^h,c^h)$ from Definition \ref{Levyfor}.
  \end{enumerate}
\end{lm} 

\subsection{Pruning forests from their ascension time}\label{ascsec2}

We can use convergence results for Kesten trees to derive from the discrete setting of \cite{ADH12} a distributional identity relating pruning processes of L\'evy trees and Kesten-L\'evy trees due to Abraham and Delmas \cite{AD12}.  
Specifically, Aldous and Pitman \cite{AP98} studied pruning at edges of Kesten trees $\cT_*$ in the following context. For pruning at edges of a Galton-Watson tree $\cT$, they noted that pruning processes can be extended to $\theta\!<\!0$ for many offspring distributions including ${\rm Poi}(1)$, when $\cT^E(\theta)\!\sim\!{\rm Poi}(e^{-\theta})$, $\theta\!\in\!\bR$. In reverse time, as $\theta\!\rightarrow\!-\!\infty$, trees become more and more supercritical, and there is an \em ascension time \em $A\!:=\!\inf\{a\ge 0\colon\Gamma(\cT^E(-a))=\infty\}\!\in\!(0,\infty)$. In the Poisson case, $(\cT^E(\theta),\theta\!>\!-A)$ has the same distribution as $(\cT^E_*(W\!+\!\theta),\theta\!>\!-\Theta)$, where $W\!\sim\!{\rm Exp}(1)$ is independent and $\Theta\!:=\!\log(W/(1\!-\!e^{-W}))$. They point out the subtlety that while the left limit $\lim_{\varepsilon\downarrow 0}\cT^E(-A-\varepsilon)$ is infinite, $\lim_{\varepsilon\downarrow 0}\cT^E_*(-\Theta-\varepsilon)$ is finite, since $W-\Theta>0$ a.s. This study was generalised by Abraham et al. \cite{ADH12} to pruning at branch points for a wide class of offspring distributions, with $\log(W/(1-e^{-W}))$ replaced appropriately. 

\begin{thm}\label{thmasc} Let $\psi$ be a critical branching mechanism that is finite on 
  $(-\theta_0,\infty)$ for some $\theta_0\in(0,\infty]$, i.e.\ we require $\int_1^\infty e^{-\theta r}\pi(dr)<\infty$ for all $\theta\in(-\theta_0,\infty)$, where $\pi$ is the L\'evy measure in {\rm(\ref{brmech})}. Furthermore, suppose that  
  $\psi(\theta)\rightarrow\infty$ as $\theta\downarrow-\theta_0$. Let $(\cF^{\rm AD}(\theta),\theta>-\theta_0)$ be a (consistently 
  extended) Abraham-Delmas pruning process of a $(\psi,\delta_x)$-L\'evy forest and 
  $A=\inf\{a\ge 0\colon\Gamma(\cF^{AD}(-a))=\infty\}$ the ascension time of its time reversal. Then
  $$(\cF^{\rm AD}(\theta),\theta\ge-A)\ed(\cF_*^{\rm AD}(W/x+\theta),\theta\ge-\Theta),$$
  where, on the right hand side
  \begin{itemize}\item $(W,\Theta)$ and $(\cF_*^{\rm AD}(\theta),\theta\ge 0)$ are independent,
     \item $W\sim{\rm Exp}(1)$ and $\Theta=-q_0^{-1}(W/x)$, where $q_0^{-1}\colon[0,\infty)\rightarrow(-\theta_0,0]$ is the inverse of
       $\theta\mapsto q_0(\theta)$, where $q_0(\theta)$ is the largest $q\ge 0$ such that $\psi(\theta+q)=\psi(\theta)$,
     \item $(\cF_*^{\rm AD}(\theta),\theta\ge 0)$ is an Abraham-Delmas pruning process of a Kesten-L\'evy forest 
       $\cF^{\rm AD}_*(0)$ defined as the concatenation of $\cF^{\rm AD}(0)$ and an independent $\psi$-Kesten-L\'evy tree $\cT_*$.
  \end{itemize}\pagebreak[2]
\end{thm}
Before the proof, let us point out that this theorem is a forest version of \cite[Corollary 8.2]{AD12}, which can be deduced here using vague convergence $\frac{1}{x}P_\psi^{\delta_x}\rightarrow\bN_\psi$ on the space of $\sigma$-finite measures on $\bT$. The proof of Theorem \ref{thmasc} makes use of the following lemma, which is of some independent interest, since it demonstrates that the ``domains of attraction'' of $\psi$-L\'evy forests (in the sense of the invariance principle of Theorem \ref{invprinc}) are non-empty for all branching mechanisms. This is essential for us and other applications that use approximation of L\'evy forests by discrete Galton-Watson forests. Similar results for Galton-Watson forests with exponentially distributed edge lengths have been pointed out in \cite{DuLG02} and exploited in \cite{DuWi07,DuWi12} to construct L\'evy forests -- we use the same families of offspring distributions here, but combined with unit edge lengths.

\begin{lm}\label{lm52} Let $\psi$ be a branching mechanism satisfying {\rm(\ref{grey})}, and let $x>0$. Then a sequence of pairs
  $(\xi_n,\mu_n)$, $n\ge 1$, such that for all $n\ge 1$ with $\psi(n)>0$
  \begin{equation}g_{\xi_n}(s)=s+\frac{\psi(n(1-s))}{n\psi^\prime(n)},\quad \mu_n=\delta_{\lfloor nx\rfloor},\label{offsp}
  \end{equation}
  satisfies the hypotheses of Theorem \ref{invprinc} with $\gamma_n=\psi^\prime(n)$. 

  Furthermore, for all $\theta\in\bR$ such that $\psi(\theta-\varepsilon)<\infty$ for some $\varepsilon>0$, the sequences 
  $(\xi_{n,\theta/n}^0,\mu_n)$, $n\ge 1$, associated as in Proposition \ref{propRhGWprun}, also satisfy the hypotheses of Theorem 
  \ref{invprinc} with the same $\gamma_n=\psi^\prime(n)$, and with branching mechanism 
  $\psi_\theta(r)=\psi(\theta+r)-\psi(\theta)$. 
\end{lm}
\begin{pf} Fix $\theta\!\in\!\bR$ with $\psi(\theta\!-\!\varepsilon)\!<\!\infty$. For $n$ sufficiently large, $\psi(n)>0$ and $\psi(n)>\psi(\theta_n)$, where $\theta_n:=n(1-e^{-\theta/n})$, so that
  $\xi_n^\theta:=\xi_{n,\theta/n}^0$ is an offspring distribution with $\xi_n^\theta(0)>0$. Specifically,
  \begin{equation}g_{\xi_n^\theta}(s)=1+e^{\theta/n}(g_{\xi_n}(se^{-\theta/n})-g_{\xi_n}(e^{-\theta/n}))=s+\frac{\psi_{\theta_n}(ne^{-\theta/n}(1-s))}{ne^{-\theta/n}\psi^\prime(n)}\label{gxintheta}
  \end{equation}
  is of the same form as (\ref{offsp}) with $n$ replaced by $ne^{-\theta/n}$ and $\psi$ by $\psi_{\theta_n}$.
  Using Laplace transforms, we express the $\nu_n$-convergence in (\ref{invass}) as in the proof of Lemma \ref{convgf} and 
  check for all $r\ge 0$ that
  $$n\psi^\prime(n)(g_{\xi_n^\theta}(e^{-r/n})-e^{-r/n})=e^{\theta/n}\psi_{\theta_n}\left(ne^{-\theta/n}(1-e^{-r/n})\right)\rightarrow\psi_\theta(r),$$
  as $n\rightarrow\infty$. The $\eta_n$-convergence is just a tightness condition and  $\limsup_{n\rightarrow\infty}n\eta_n(\lfloor a\gamma_n\rfloor)<\infty$ is in fact sufficient; see e.g. \cite[Theorem 2.3.1]{DuLG02}. Now let $w^\theta_n(k)\!:=\!e^{-\eta^\theta_n(k)}\!:=\!Q_{\xi_n^\theta}(\Gamma\!\le\! k)$. Denote by $q_{\xi_n^\theta}$ the smallest $q\!\ge\! 0$ with $g_{\xi_n^\theta}(q)\!=\!q$. As $w_n^\theta(0)\!=\!0$ and $\xi_n^\theta(0)\!>\!0$, and $g_{\xi_n^\theta}$ is increasing on $[0,1]$, we have $w_n^\theta(k+1)\!=\!g_{\xi_n^\theta}(w_n^\theta(k))\!\in\!(0,q_{\xi_n^\theta})$. By (\ref{gxintheta}), this implies $\widetilde{\eta}_n^\theta(k)\!:=\!ne^{-\theta/n}(1\!-\!w_n^\theta(k))\!\in\!(q_0(\theta_n),\infty)$, where $q_0(\theta_n)$ is the largest $q\!\ge\! 0$ with $\psi_{\theta_n}(q)\!=\!0$. We obtain for all $k\!\ge\! 0$
  $$0<\left(\widetilde{\eta}_n^\theta(k)-\widetilde{\eta}_n^\theta(k+1)\right)\psi^\prime(n)
     =ne^{-\theta/n}\left(g_{\xi_n^\theta}(w_n^\theta(k))-w_n^\theta(k)\right)\psi^\prime(n)
      =\psi_{\theta_n}(\widetilde{\eta}_n^\theta(k)).
  $$
  Now let $a>0$. Since $\psi_{\theta_n}$ is positive increasing on $(q_0(\theta_n),\infty)$, we find
  $$\frac{\lfloor a\psi^\prime(n)\rfloor}{\psi^\prime(n)}=\sum_{k=1}^{\lfloor a\psi^\prime(n)\rfloor}\frac{\widetilde{\eta}_n^\theta(k-1)-\widetilde{\eta}_n^\theta(k)}{\psi_{\theta_n}(\widetilde{\eta}_n^\theta(k-1))}\le\int_{\widetilde{\eta}_n^\theta(\lfloor a\psi^\prime(n)\rfloor)}^{ne^{-\theta/n}}\frac{dx}{\psi_{\theta_n}(x)}.$$
  Then for $n$ and $x$ sufficiently large $\theta_n\ge\theta-\varepsilon$ and  $\psi^{\prime\prime}(x+\theta)\le 2\beta+1$, and as $x\rightarrow\infty$\vspace{-0.1cm}
 $$\frac{\psi_{\theta-\varepsilon}(x)}{\psi_\theta(x)}=\frac{\psi(x\!+\!\theta\!-\!\varepsilon)-\psi(\theta\!-\!\varepsilon)}{\psi(x\!+\!\theta)-\psi(\theta)}\sim\frac{\psi^\prime(x\!+\!\theta\!-\!\varepsilon)}{\psi^\prime(x\!+\!\theta)}=1-\frac{\int_0^\varepsilon\psi^{\prime\prime}(x\!+\!\theta\!+\!r)dr}{\psi^\prime(x\!+\!\theta)}\ge 1-\frac{\varepsilon(2\beta\!+\!1)}{\psi^\prime(x\!+\!\theta)}\rightarrow 1.\vspace{-0.1cm}$$
  In particular, for any $c\in(0,1)$ and $x$ and $n$ large enough, say $x\ge x_0$ and $n\ge n_0$, we can bound 
  $\psi_{\theta_n}(x)\ge\psi_{\theta-\varepsilon}(x)\ge c\psi_\theta(x)$.  
  Now suppose for contradiction that $\limsup_{n\rightarrow\infty}\widetilde{\eta}_n^\theta(\lfloor a\psi^\prime(n)\rfloor)=\infty$ and choose a subsequence $(n(k),k\!\ge\! 1)$ along which $\lim_{k\rightarrow\infty}\widetilde{\eta}_{n(k)}^\theta(\lfloor a\psi^\prime(n(k))\rfloor)\!=\!\infty$. Then there is $k_0\!\ge\! 1$  such that for $k\!\ge\! k_0$, we have $n(k)\!\ge\! n_0$ and $\widetilde{\eta}_{n(k)}^\theta(\lfloor a\psi^\prime(n(k))\rfloor)\!\ge\! x_0$ so that we get \vspace{-0.2cm}
  $$\int_{\eta^\theta(a)}^\infty\frac{dx}{\psi_\theta(x)}=a\le\liminf_{k\rightarrow\infty}\int_{\widetilde{\eta}_{n(k)}^\theta(\lfloor a\psi^\prime(n(k))\rfloor)}^{n(k)e^{-\theta/n(k)}}\frac{dx}{\psi_{\theta_{n(k)}}(x)}\le\frac{1}{c}\liminf_{k\rightarrow\infty}\int_{\widetilde{\eta}_{n(k)}^\theta(\lfloor a\psi^\prime(n(k))\rfloor)}^\infty\frac{dx}{\psi_\theta(x)}=0,\vspace{-0.1cm}$$
  which is the required contradiction. Hence $\limsup_{n\rightarrow\infty}\widetilde{\eta}_n^\theta(\lfloor a\psi^\prime(n)\rfloor)<\infty$, and this easily entails the tightness condition.
\end{pf}

\begin{pfofthmasc} Consider $(\xi_{n,\theta/n}^0,\mu_n)$ as in Lemma \ref{lm52} and recall notation 
  $\xi_n^\theta\!:=\!\xi_{n,\theta/n}^0$, $n\!\ge\! 1$. For $n$ sufficiently large, we can consider  
  $(\cF_n^B(\theta/n)/\gamma_n,\theta\!\ge\!-\theta_1)$. By Theorem \ref{thmprocconv}, we have   
  $(\cF_n^B(\theta/n)/\gamma_n,\theta\!\ge\!-\theta_1)\convd(\cF_n^{\rm AD}(\theta),\theta\!\ge\!-\theta_1)$ 
  in $\bD([-\theta_1,\infty),\bT)$ for all $\theta_1\!<\!\theta_0$, and hence\vspace{-0.1cm}
  $$(\cF_n^B(\theta/n)/\gamma_n,\theta>-\theta_0)\convd(\cF_n^{\rm AD}(\theta),\theta>-\theta_0)\qquad\mbox{in $\bD((-\theta_0,\infty),\bT)$.}\vspace{-0.1cm}$$
  By Theorem \ref{invprinc} and Lemma \ref{convgf}, we have for all $\theta>-\theta_0$ that, as $n\rightarrow\infty$, 
  $$\bP(-A_n\le\theta)=\bP(\Gamma(\cF_n^B(\theta/n))<\infty)=q_{\xi_n^\theta}^{\lfloor nx\rfloor}\rightarrow e^{-xq_0(\theta)}=\bP(\Gamma(\cF^{\rm AD}(\theta))<\infty)=\bP(-A\le\theta),$$
  where $q_0(\theta)$ is the largest $q$ with $\psi_\theta(q)\!=\!0$, $\theta\!<\!0$, and $q_0(\theta)\!=\!0$, $\theta\!\ge\! 0$. We deduce that  
  the distributions of $(A_n,\cF^B_n(\theta/n),\theta\!>\!-\theta_0)$ are tight. Let $n(k)$, $k\!\ge\! 1$, be a subsequence along which the 
  joint distributions converge. W.l.o.g., convergence holds almost surely, where $A_{n(k)}\!\rightarrow\! A_\infty$, as $k\!\rightarrow\!\infty$, a priori does not 
  mean that the limit is the ascension time $A$ of the limit process. However, as $A_\infty\ed A$ and $\Gamma\colon\bT\!\rightarrow\![0,\infty]$ is continuous, $A_\infty\!\le\!\liminf_{k\rightarrow\infty}A_{n(k)}\!=\!A$ gives $A_\infty\!=\!A$ a.s. This identifies the joint limiting distribution, which does not depend on the subsequence, and joint convergence in distribution follows. Since the limiting distribution of $A$ has no atoms, this implies\vspace{-0.2cm}
  \begin{equation}\label{c1}(\cF_n^B(\theta/n)/\gamma_n,\theta\ge-A_n)\convd(\cF^{\rm AD}(\theta),\theta\ge-A).
  \end{equation}
  To apply results from \cite{ADH12}, we check that \cite[Condition (4.10)]{ADH12} holds: $$\xi_n^\theta(0)=g_{\xi_n^\theta}(0)=\frac{\psi(n)-\psi(n(1-e^{-\theta/n}))}{n\psi^\prime(n)e^{-\theta/n}}=0\iff\psi(n)=\psi(n(1-e^{-\theta/n}))$$ 
  and since $\psi(n(1-e^{-\theta/n}))\rightarrow\infty$ as $\theta\downarrow-\theta_0$, there is $\overline{\theta}_n<0$ with $\xi_n^{\overline{\theta}_n}(0)=0$. By considering $A_n$ as the minimum ascension time of the $\lfloor nx\rfloor$ trees in $(\cF_n^B(\theta),\theta\ge 0)$, it 
  follows easily from \cite[Propositions 4.5-4.7]{ADH12} that\vspace{-0.1cm}
  \begin{equation}\label{d1}(\cF_n^B(\theta/n),\theta\ge-A_n)\ed(\cF_{n*}^B(W/\lfloor nx\rfloor+\theta/n),\theta\ge-\Theta_n),
  \end{equation}
  where on the right-hand side,
  \begin{itemize}\item $(W,\Theta_n)$ and $(\cF_{n*}^B(\theta),\theta\ge 0)$ are independent.
    \item $W\sim{\rm Exp}(1)$ and $\Theta_n=-nq_n^{-1}(e^{-W/\lfloor nx\rfloor})$, where $q_n^{-1}\colon[0,1)\rightarrow(\overline{\theta}_n,0]$ is the inverse function of $q_n(\theta)=q_{\xi_n^\theta}$, where $q_{\xi_n^\theta}$ is the smallest $q\ge 0$ with $g_{\xi_n^\theta}(q)=q$.
    \item $(\cF_{n*}^B(\theta),\theta\ge 0)$ is a pruning process with pruning at branch points of a Kesten forest $\cF_{n*}^B(0)$, obtained by concatenating one $\xi_n$-Kesten tree and $\lfloor nx\rfloor-1$ ${\rm GW}(\xi_n)$-real trees.
  \end{itemize}
  Finally, we check joint convergence $(nW/\lfloor nx\rfloor,nq_n^{-1}(e^{-W/\lfloor nx\rfloor}))\convd(W/x,q_0^{-1}(W/x))$ by noting
  that $g_n(z)=nq_n^{-1}(e^{-z/n})$ is the inverse of the strictly increasing function $\theta\mapsto-n\log q_n(\theta/n)$, which 
  we have seen converges to the continuous and strictly increasing function $\theta\mapsto q_0(\theta)$. This implies that the 
  convergence is uniform and that inverse functions also converge, as required. We can now conclude from Theorem \ref{kesconv} and Remark \ref{kestenforests} that
  $$(\cF_{n*}^B(W/\lfloor nx\rfloor+\theta/n)/\gamma_n,\theta\ge-\Theta_n)\convd(\cF_*^{\rm AD}(W/x+\theta),\theta\ge-\Theta),$$
  which completes the proof by uniqueness of limits, noting (\ref{d1}) and comparing with (\ref{c1}).
\end{pfofthmasc}

\section*{Acknowledgements}

This work was started during a research visit of the second author to Beijing Normal University. We would like to thank Beijing Normal University for support during this research visit. H.\ He is supported by the Fundamental Research Funds for the Central Universities (2013YB59) and NSFC (No.\ 11201030, 11371061).

\bibliographystyle{abbrv}
\bibliography{prun}

\end{document}